\RequirePackage{fix-cm}
\documentclass[final]{svjour3}
\usepackage{amsmath}
\usepackage{amssymb}
\smartqed

\usepackage{latexsym}
\usepackage{srcltx}
\usepackage{graphics}
\usepackage{graphicx}
\usepackage{epsfig}
\usepackage{GGGraphics}
\newtheorem{prop}{Proposition}
\newtheorem{cor}{Corollary}
\newtheorem{defi}{Definition}
\newtheorem{conj}{Conjecture}
\newtheorem{ex}{Example}

\newcommand{\co}{{\mathbb C}}
\newcommand{\rd}{{\mathbb R^d}}
\newcommand{\re}{{\mathbb R}}
\newcommand{\n}{{\mathbb N}}

\newcommand{\z}{{\mathbb Z}}

\newcommand{\cV}{{\mathcal{V}}}
\newcommand{\cU}{{\mathcal{U}}}

\newcommand{\cH}{{\mathcal{H}}}

\newcommand{\cM}{{\mathcal{M}}}

\newcommand{\cR}{{\mathcal{R}}}
\newcommand{\cT}{{\mathcal{T}}}

\newcommand{\cl}{{\rm c\ell}}

\begin{document}


\title{Exact computation of joint spectral characteristics of linear operators}

\author{Nicola Guglielmi\thanks{The research of the first author is supported by Italian M.I.U.R. (PRIN2007 project) and G.N.C.S.} \and
Vladimir Protasov\thanks{The research of the second author is supported by the RFBR grants No 11-01-00329 and  No 10-01-00293,
and by the grant of Dynasty foundation.}}
\institute{Nicola Guglielmi \at Department of Pure and Applied Mathematics, \\ University of L'Aquila, Italy \\ \email{guglielm@univaq.it} \and
Vladimir Protasov \at Department of Mechanics and Mathematics, Moscow State University, Vorobyovy Gory, Moscow, Russia 119992 \\ \email{v-protassov@yandex.ru}}

\date{}

\journalname{
}


%
%
%

\maketitle

\begin{abstract}

We address the problem of the exact computation of two joint spectral characteristics of a family of linear
operators, the joint spectral radius (in short JSR) and the lower spectral radius (in short LSR), which are
well-known different generalizations to a set of operators of the usual spectral radius of a linear operator.
In this article we develop a method which - under suitable assumptions - allows to compute the JSR and the
LSR of a finite family of matrices exactly.
We remark that so far no algorithm was available in the literature to compute the LSR exactly.

The paper presents necessary theoretical results on extremal norms (and on extremal antinorms)
of linear operators, which constitute the basic tools of our procedures, and a detailed description
of the corresponding algorithms for the computation of the JSR and LSR (the last one restricted to
families sharing an invariant cone).
The algorithms are easily implemented and their descriptions are short.

If the algorithms  terminate in finite time,
then they construct an extremal norm (in the JSR case) or antinorm (in the LSR case)
and find their exact values; 
otherwise they provide upper and lower bounds that both converge to the
exact values. A~theoretical criterion for termination in finite time is also derived.
According to numerical experiments, the algorithm for the JSR finds the exact value for the vast majority
of matrix families in dimensions~${\le 20}$. For nonnegative matrices it works faster and finds JSR in
dimensions of order~$100$ within a few iterations; the same is observed for the algorithm computing the LSR.
To illustrate the efficiency of the new method we are able to apply it in order to give
answers to several conjectures which have been recently stated in combinatorics, number theory,
and the theory of formal languages.


\begin{flushleft}
\noindent  \textbf{AMS 2010} {\em subject
classification: } 15A60, 15-04, 15A18, 90C90
\end{flushleft}

\end{abstract}

\section{ Introduction and background}

The joint spectral characteristics of linear operators are now applied in many areas,
from functional analysis and dynamical systems to discrete mathematics and number theory.
We focus on two characteristics: the joint spectral radius and the lower spectral radius, and elaborate
a method of their exact computation applicable even for relatively high dimensions.

The joint spectral radius of a set of matrices is a measure
identifying the highest possible rate of growth of the norm of
products of matrices (with no ordering and with repetition permitted)
of the set.
In contraposition, the lower spectral radius defines the lowest possible rate
of growth.
Both measures appear in several applications (see e.g. Strang \cite{Str}).

In this paper we consider the problem of the computation of both
joint spectral characteristics for a finite set of matrices.
In contrast to the fact that in the last twenty years much effort has
been devoted  to the computation of the joint spectral radius,
very little is known about computing the lower spectral radius
(to the best of our knowledge, the only available method of its approximate
computation was presented in~\cite{PJB}).


The joint spectral radius originated with Rota and Strang in 1960~\cite{RS}, and became extremely popular
after  Daubechies and Lagarias~\cite{DL} revealed its role in the study of refinement equations and
wavelets. Since then it has found applications in functional equations, approximation, probability,
combinatorics, etc. (see~\cite{J,PJB} for the extensive bibliography).  Let
$$
\cM = \{A_1, \ldots , A_m\}
$$
be a finite family of linear operators acting in~$\re^d$. We write
$$\, \cM^k \, = \, \bigl\{A_{d_k}\ldots A_{d_1} \ \bigl| \
d_j \, \in \, \{1, \ldots , m\}\, , \, j=1, \ldots , k\, \bigr\}$$
for the set of all~$m^k$ products of length~$k$
of  operators from~$\cM$. The {\em joint spectral radius} (JSR) of the family~$\cM$ is
\begin{equation}\label{jsr}
\widehat \rho (\cM)  \quad = \quad  \lim_{k \to \infty} \max_{B \in
\cM^k} \|B\|^{\, 1/k}\, .
\end{equation}
This limit exists for every family~$\mathcal{M}$ and does not depend on the norm in~$\re^d$~\cite{BW}.
Clearly, if $\cM$ consists of one operator~$A_1$, then $\widehat \rho (\cM)\, = \, \rho (A_1)$,
where $\rho (A_1)$ is the (usual) spectral radius of~$A_1$, which is the maximal modulus of its
eigenvalues.

For any family~$\cM$ there is a positive constant~$c_1$ such that
$$\max\limits_{d_1, \ldots , d_k}\|A_{d_k}\ldots A_{d_1}\|\, \ge \, c_1\, \widehat \rho^{\, k}$$
for every $k \in \n$.
The family is called {\em non-defective} if the inverse estimate holds, that is
there is a  constant~$c_2$ such that
$\max_{d_1, \ldots , d_k}\|A_{d_k}\ldots A_{d_1}\|\, \le \, c_2\, \widehat \rho^{\, k}$.
It appears that if a family~$\cM$ is irreducible, i.e., its operators do not share a common
nontrivial invariant subspace of~$\re^d$, then it is non-defective~\cite{P1}.
Thus, for an irreducible family one has $\, \max_{d_1, \ldots , d_k}\|A_{d_k}\ldots A_{d_1}\| \, \asymp \, \widehat \rho^{\, k}$,
where the symbol $\asymp$ denotes asymptotic equivalence
($\, a_k \, \asymp \, b_k\, $ if there are constants $c_1, c_2 > 0$ such that
$\, c_1 a_k \, \le \, b_k \, \le \, c_2 a_k\, $ for all $k$).
Whence,  the joint  spectral radius is the exponent of polynomial growth for the largest norm
of operator products of length~$k$. The geometric sense of JSR is the following:
$\widehat \rho < 1$ if and only if there exists a norm in $\re^d$ such that $\|A_j\| < 1$ for all
$j=1, \ldots , m$, where $\|A\| = \sup_{\|x\| \le 1}\|Ax\|$ is the corresponding operator norm.
In other words,  $\widehat \rho < 1$ precisely when there is a norm in $\re^d$ with respect to which all operators from~$\cM$
are contractive. So, it is natural to expect that each family of operators possesses some special norms
related to JSR. The following theorem established in 1988 by Barabanov~\cite{B}
shows that this is indeed the case, at least for irreducible families. A norm $\|\, \cdot\, \|$
in~$\re^d$ is called {\em invariant} for~$\cM$ if there is a number $\lambda \ge 0$ such that
$\max_{j=1, \ldots , m}\|A_jx\| \, = \, \lambda \, \|x\|\, $
for every $x \in \re^d$. It is shown easily that for every invariant norm one has $\lambda = \widehat \rho (\cM)$.
\begin{theorem}\label{th.bar}~\cite{B}
Every irreducible family~$\cM$ possesses an invariant norm.
\end{theorem}
 In practice it suffices to get a special norm with some weaker
requirements, the so-called {\em extremal norm}.
\begin{defi}\label{i.d2}
A norm $\| \, \cdot \, \|$ is called extremal
for $\cM\, $ if $\, \|A_j\, x\| \, \le \, \widehat \rho\, \|x\|\, $ for all $x \in \re^d$.
\end{defi}
Thus, the norm is extremal if and only if $\max_{j=1, \ldots , m}\|A_j\|\, = \, \widehat \rho $. Indeed,
from the definition it follows that $\max_{j=1, \ldots , m}\|A_j\|\, \le \, \widehat \rho $; on the other hand,
 the submultiplicativity of operator norms  yields $\max_{j=1, \ldots , m}\|A_j\|\, \ge \, \widehat \rho $.
Whence, this inequality becomes an equality precisely for extremal norms. This property justifies
the term ``extremal''.

  Clearly, any invariant norm is extremal, but not vice versa.
Let us remark that for every extremal norm one has
$\max_{d_1, \ldots , d_k}\|A_{d_k}\ldots A_{d_1}\|\, = \, \widehat \rho^{\, k}\, , \ k \in \n $, i.e.,  the asymptotic equality
becomes a sharp equality for all~$k$. In particular, for $k=1$ we have $\max_{j} \|A_j\|\, = \, \widehat \rho$.
Thus, if we know an extremal norm, then we have the exact value of JSR. The main idea of the approach presented in this paper
is to find JSR by constructing (in an iterative way) an extremal norm.

We first present an algorithm for general sets of matrices,
which under some suitable assumptions is able to check if a certain
product in the multiplicative semigroup is spectrum maximizing.
The algorithm is based on the computation of an extremal
norm whose unit ball is a balanced polytope, and we provide it by
a new criterion assuring a finite time termination and a new stopping
condition.

Then we analyze sets of matrices having an invariant cone (a most
important case is given by families of nonnegative matrices). For
such sets we refine the algorithm for the general case and exploit
the invariant property of the set in order to make the algorithm
faster.
Under this assumption we are also able to determine an
algorithm for the exact computation of the lower spectral radius,
which appears to be the first algorithm able to provide an exact
value of this important measure.

The algorithms compute respectively a bounded and an unbounded
polytope which represent the unit balls of, respectively, an extremal
norm and antinorm for the considered set. %


We write the formal routines of the  algorithms and illustrate their efficiency by suitable examples and by
numerical tests with randomly generated matrices. As we shall see, the algorithms find
the exact value of JSR for general families of matrices (under some minor restrictions) in dimensions up to~$20$.
For nonnegative matrices they work surprisingly fast even in dimension $d=100$ and higher.
Let us remark that our approach does not apply successfully to all families.
There are cases, in fact, where the algorithms we propose are not able to finitely compute the exact
values of the considered joint spectral characteristics, but only to approximate them.

In view of negative complexity results for the problem of JSR computation~\cite{BT},
this is unlikely that there are effective methods applicable for all families of operators.
Nevertheless, we claim that our approach works for the vast majority of families.
The results of many numerical tests with randomly generated matrices of dimension from $5$ to $100$
 (some of them are presented at the end of this paper)
 confirm this claim. In all the cases the algorithms found the exact values of JSR.
As a further confirmation of this, we are able to apply the new method to solve several open problems
in combinatorics and discrete mathematics.
\smallskip

In the literature there are several methods for the computation of the JSR. Some
of them work only for small dimensions~$d$,
but give either an exact or a very accurate value of ~$\widehat \rho$. For example, the method of polytope norms~\cite{P1,GZ1,BJP,GZ2,CGSZ};
see also special methods in~\cite{V,P3,HMR} elaborated for particular matrices.

Other methods aim to an approximate computation, such as the Kronecker lifting method~\cite{P2,BN},
ellipsoidal norm method~\cite{BNT}, Gripenberg's branch-and-bound method~\cite{DL,G} can work for bigger
dimensions (mostly, up to~$20$), but produce pretty rough estimations.
Recent approaches involving some modern tools of convex optimization (conic and semidefinite programming,
sum-of-squares approximation, etc.) have rather good accuracy for higher dimensions ($10$ or even bigger)~\cite{PJ,PJB}.
Most of those methods are actually based on the same simple idea.
For each $k$ we have
\begin{equation}\label{eq.radius}
\max_{B \in \cM^k}\, \bigl[ \, \rho(B) \, \bigr]^{1/k}\quad \le
\quad \widehat \rho(\cM)\quad \le \quad \max_{B \in \cM^k}\, \bigl\|\,
B \, \bigr\|^{1/k}\, .
\end{equation}
The right hand side of this inequality converges to $\widehat
\rho$ as $k \to \infty$, which follows from the definition. The upper limit of the left hand side
also equals to~$\widehat \rho$ \cite{BW}. So, choosing $k$ large enough, it is
possible to approximate JSR as close as we need. However, this
possibility is purely theoretical, because in most of practical cases
the number $k$ grows as $\frac{C}{\varepsilon}$, where $\,
\varepsilon > 0$ is the relative accuracy of the JSR approximation, and $C > 0$ is a
constant, which may be large  for high dimensions~$d$. That is why
the number of matrix products of~$\cM^k$ to look over becomes
enormous. The reason is that the norm~$\|\cdot \|$ in the
right hand side of~(\ref{eq.radius}) may not suit our
family~$\cM$, i.e., it may be far from the extremal norm of that
family. That is why, to achieve a good approximation of JSR one needs to find an appropriate
norm in~$\re^d$ for  the right hand side
of~(\ref{eq.radius}).  Actually, all the
methods of JSR computation use various techniques
to find such a special norm for a given family~$\cM$.
Those are, for instance, a polytope norm~\cite{P1,GZ1},  an ellipsoidal norm~\cite{BNT},
 a norm generated by a cone~\cite{PJB}, a norm defined by a
sum-of-squares polynomial~\cite{PJ},  etc. Sometimes this
idea leads even to finding the precise values of JSR. This happens when
 both the inequalities
in~(\ref{eq.radius}) become equalities. If we write
$$
\rho_k = \max\limits_{B \in \cM^s\, ,  \, s \le k}\, \bigl[  \rho(B)
\bigr]^{1/s},
$$
then for the extremal norm we have $\|A\|\, \le \,
\rho_k\, $ for all $\, A \in \cM$, and therefore $\rho_k = \widehat
\rho$. For example, if all the matrices of $\cM$ are symmetric, then
the Euclidean norm is extremal, if they are all column-stochastic,
then the $L_1$-norm is so. In some practical
cases people succeed in finding extremal norms for concrete
pairs of matrices arising in various applications: Gripenberg~\cite{G}
(matrices of Daubechies wavelets, dimensions from $4$ to $7$), Hechler, M\"o\ss ner, and Reif~\cite{HMR}
(matrices of the four-point subdivision schemes, $d = 4$),
Protasov~\cite{P4} (de Rham matrices, $d=2$), Guglielmi, Wirth, and Zennaro~\cite{GWZ}
(matrices of the Blondel-Theys-Vladimiov family, $d=2$), Protasov~\cite{P3}
(matrices of the binary partition function, $d=4, \ldots , 12$),
Villemoes~\cite{V} (matrices of refinement equations, $d=2$), Guglielmi, Manni and Vitale~\cite{GMV} (matrices
in Hermite subdivision schemes) etc.
The JSR computation in each case was  a nontrivial problem and
required special tricks applicable only for some narrow classes of
matrices.

The method of exact JSR computation presented in this paper
is related to previous works (see \cite{P1,GZ1,BJP,GZ2,CGSZ})
and aims to develop further ideas both for the general case, which
we are goind to recall, and for certain specific important cases,
like that of nonnegative matrices.
The method is applicable for all  families of matrices, under some general
assumptions. The main idea proposed in the above mentioned papers
is to build an extremal norm, whose
unit sphere is a polytope. At the first step we look over all
products of matrices from $\cM$ of length at most~$l$, and find a
product $\Pi$, for which the value $\bigl[\rho(\Pi )\bigr]^{1/n}$
is maximal ($n$ is the length of~$\Pi$). Then we denote this value by $\rho_l$ and try to prove that
$\widehat \rho (\cM)  \, = \,  \rho_l$.
\begin{defi}\label{i.d1}
A product $\Pi \in \cM^n$ is a spectrum maximizing product
(s.m.p.) if $$\, \bigl[\rho(\Pi )\bigr]^{1/n}\, = \, \widehat \rho (\cM). $$
\end{defi}
To prove that $\Pi$ is an s.m.p. it suffices to have an extremal
norm $\|\cdot \|$ in $\re^d$, for which $\|A_j\| \, \le \, \rho_l\, , \
j=1, \ldots , m$. By~(\ref{eq.radius}) in this case we indeed have $\rho_l \, = \, \widehat \rho$.
We try to build a {\em polytope extremal
norm},  whose unit sphere is some polytope~$P$. Such a polytope will
also be called {\em extremal}. It is characterized by the property
$\, A_jP \, \subset \, \rho_l P\, , \ j=1, \ldots , m$. The polytope is
constructed successively: its first vertices are the {\em leading
eigenvector}~$v_1$ of~$\Pi$ (i.e., the eigenvector corresponding
to the largest by modulo eigenvalue, which is assumed to be real
for the moment), the leading eigenvectors~$v_i$ of the~$(n-1)$
cyclic permutations of~$\Pi$, and the same vectors taken with minus, i.e., $\, -v_i$.
 We call an eigenvalue $\lambda$ of an operator $A$ {\em leading} if $|\lambda| = \rho(A)$.

Then we consider their images
$(\rho_l)^{-1}A_jv_i\, , \ j = 1, \ldots , m$ and remove those are
in the convex hull of the previous ones, etc., until we obtain a
set of points $\cV$ such that
$$ (\rho_l)^{-1}A_j \cV\, \subset \,{\rm co}_s\, (\cV)\, , \ j=1, \ldots , m.$$
By ${\rm co}_s (\cV)$ we denote the {\em symmetrized convex hull}:
${\rm co}_s\, (\cV) \, = \, {\rm co}\, \bigl( \cV \cup (-\cV) \bigr)$, where ${\rm co}\, (\cdot )$ is the (usual) convex hull.
 Then the polytope~$P \, = \, {\rm co}_s
(\cV)$ possesses the desired property: $(\rho_l)^{-1}A_j P
\subset P$, so $P$ is an extremal polytope. This implies $\widehat
\rho = \rho_l$. The algorithm involves standard tools of linear programming.

In case the leading eigenvalue of $\Pi$ is complex, one has to replace
polytopes by the so-called {\em complex polytopes} (see e.g. \cite{GZ3}).

Our goal is to develop this approach for general families of matrices in higher dimensions,
to analyze  the structure of extremal polytopes and to derive  the conditions of
convergence of this algorithm.
\smallskip

 Let us now emphasize the shortcomings of our approach. First
 of all, not every family of matrices has an s.m.p. Moreover, even if a non-defective family~$\cM$ possesses an s.m.p.,
 it may not have extremal polytopes (neither real nor complex~\cite{JP}).
 For such families our method apparently does not work. Another
 disadvantage appears, when the s.m.p. is not unique, up to
 cyclic permutations. In  this case an extremal polytope, even if
 it exists, in general cannot be found by our method.
The first two cases are rather pathological. It required constructing
  special nontrivial examples to show that they are
 possible~\cite{BTV,JP}. The third case of multiple s.m.p., in contrast,
  being also quite rare in general, nevertheless,
 appears in practical applications.

 We believe that our method can be extended to this case as well,
 which may be a challenging problem for further research.


 To work with those ``bad cases'', we apply our approach also to approximate
 computation of JSR. The algorithm constructs a polytope, which is either extremal or not.
 If it is, then the JSR is found. Otherwise, we stop the algorithm after a certain iteration,
 say the $N$-th, and use the obtained polytope as  a unit ball of the corresponding norm in
 estimations~(\ref{eq.radius}). In most cases this
 gives very sharp bounds for JSR. Thus, for an arbitrary family~$\cM$ the algorithm either
produces an extremal polytope, or a polytope norm that gives good upper and lower bounds
 for JSR. Proposition~\ref{p2.r} guarantees that both those bounds converge to $\widehat \rho(\cM)$
 as $N \to \infty$.
 \smallskip

The second part of the paper deals with the {\em lower spectral
radius} (LSR) defined as follows:
\begin{equation}\label{lsr}
\check \rho (\cM)  \quad = \quad  \lim_{k \to \infty} \min_{B \in
\cM^k} \|B\|^{\, 1/k}\, .
\end{equation}
Thus,  LSR is the exponent of asymptotic growth of the {\em
minimal} product of operators from the family~$\cM$. This notion
defined in~\cite{Gu}  have been applied in problems of dynamical
systems,  functional analysis, coding theory, combinatorics,
number theory,  etc. (see~\cite{J} for many references). The limit in~(\ref{lsr}) always exists and
does not depend on the norm. A simple observation is that LSR can
be estimated by the usual spectral radii as follows:
\begin{equation}\label{lsr.est}
\check \rho (\cM)  \quad \le \quad   \min_{B \in \cM^k}
\bigl[ \rho (B)\bigr]^{\, 1/k}\quad \le \quad  \min_{B \in \cM^k}
\|B\|^{\, 1/k}\, .
\end{equation}
In contrast to inequality~(\ref{eq.radius}) for JSR,
estimation~(\ref{lsr.est}) gives only  upper bounds. In fact,
there is no effective lower bounds for LSR, and this causes the
main difficulty for its computation. Basically, the lower spectral
radius is still harder to compute or to estimate than the joint
spectral radius (see, for instance,~\cite{TB} for the corresponding
complexity results). The notions of invariant and extremal norms
cannot be directly extended to LSR. The reason is that the
operation of taking minimum of several functions, in contrast to
the maximum, does not obey convexity. This means that the
pointwise minimum of several convex functions may not be convex.
Hence, the function~$f(x)\, = \, \min_{j=1, \ldots , m}\,
\|A_jx\|\, , \ x \in \re^d$, in general, is not a norm in~$\re^d$.
To overcome this difficulty, we use in Section~\ref{sec:LSR}  a notion of {\em
antinorm} defined on a convex cone $K \subset \re^d$
(Definition~\ref{d.antinorm}). This notion originated in~\cite{P6} to
study the Lyapunov exponents of linear operators. As we shall see,
it can also be applied  to analyze  the lower spectral radius. We
prove that every family of operators that share a common invariant
cone~$K$ possesses an extremal antinorm on that cone
(Theorem~\ref{th20}). This allows us to extend the new approach
to the LSR computation, replacing norms by  antinorms, and polytopes by {\em
infinite polytopes}, i.e. the sets of the type ${\rm co}\, (\cV)
\, + \, K$, where $\cV\subset \re^d$ is a finite set, and
$K\subset \re^d$ is a cone. In particular, this approach can be
used for nonnegative matrices, since the corresponding operators preserve the cone $K =
\re^d_+$. This yields  an algorithm of exact computation of LSR for
nonnegative matrices. In numerical examples we show that the
algorithm works well for rather big dimensions (like $d=100$). Let
us note that the problem of LSR computation for nonnegative
matrices arise naturally in combinatorics,
discrete mathematics, and number theory~\cite{C,P3,JPB2}. Some of those
applications will be considered in detail in Section~\ref{sec:numer}. See also~\cite{MS,FV} for
 applications to  the problem of stabilization of switched linear systems.
\smallskip

The  main results of the paper can be summarized as follows:
\smallskip

\begin{itemize}

\item[(i) ] we analyze the considered algorithm for the JSR
computation of an arbitrary family and improve it by
elaborating a stopping criterion that indicates whether a chosen
product~$\Pi$ can be an s.m.p. or not.
If our initial guess is wrong, and~$\Pi$ is not an s.m.p., then
the criterion determines it (usually, after a few iterations) and
suggests a new candidate for s.m.p. with a bigger spectral radius.
\smallskip

\item[(ii) ] Theorem~\ref{th.cond-r} in Section~\ref{sec:CFT} gives a criterion for a family~$\cM$ insuring that
the algorithm terminates within finite time, i.e.,  produces an extremal polytope.
\smallskip

\item[(iii) ] we improve the considered algorithm when applied to
nonnegative matrices; the new algorithm finds the exact values of
JSR in much higher dimensions (up to~$d = 100$);
\smallskip

\item[(iv) ] we obtain a new algorithm which is able to exactly
compute the LSR for families of nonnegative matrices, by computing
a polytope extremal antinorm;
\smallskip

\item[(v) ] as examples we compute the exact values of JSR for special
families of matrices (of dimensions up to~$40$) from well-known
problems of combinatorics and number theory. This, in particular,
allows us to solve three open problems. We discuss this aspect
below in more detail.
\smallskip

\item[(vi) ] we provide numerical tests with randomly generated matrices
(both arbitrary and nonnegative), showing that for all considered cases the
algorithms produce extremal polytopes and, consequently, the exact value of
the JSR (LSR).
\smallskip

\end{itemize}

The structure of the paper is the following. We describe the
algorithm for JSR computation in three possible cases, which will
be considered separately and called~\textbf{(R)}, \textbf{(C)} and
\textbf{(P)}. The case~\textbf{(R)}, when the leading eigenvalue
of the product~$\Pi \in \cM^n$ (a candidate for s.m.p.) is real,
is recalled and further analyzed
in Section~\ref{sec:JSR}. We discuss an algorithm for constructing
an extremal polytope and for computing JSR, give necessary
explanations and proofs, and establish two efficiency results: on
the stopping criterion (to indicate within finite
time, whether the chosen product $\Pi$ is s.m.p. or not) and on
the estimation for JSR. Thus,
Algorithm~\textbf{(R)} either terminates within finite time, in
which case JSR is found, or produces lower and upper bounds
converging to JSR. According to our numerical experiments with
randomly generated matrices (Section \ref{sec:numer}), for almost all matrix
families Algorithm~\textbf{(R)} finds the exact value of JSR, and
works efficiently for dimensions up to $20$.

In Section~\ref{sec:C} we briefly consider the
case~\textbf{(C)}, when the leading eigenvalue of~$\Pi$ is
complex, where we refer to \cite{GWZ,GZ2,GZ3}.
The algorithm and all the efficiency results are very
similar, but with complex polytopes.
By the numerical results in
Section~\ref{sec:numer}, it works slower than Algorithm~\textbf{(R)}, and works
in smaller dimensions.

In Section~\ref{sec:P} we consider
the case~\textbf{(P)}, when all matrices are nonnegative. In this
case the corresponding  Algorithm~\textbf{(P)} works faster and
much more efficiently. Of course,~\textbf{(P)} is a special case
of~\textbf{(R)}, which is, in turn, a special case
of~\textbf{(C)}. In fact all the three algorithms are very
similar and differ in a few  key details. Nevertheless, we describe them
separately and independently of each other for convenience of the
reader. Besides, their practical efficiency is very different, and
it would be non-reasonable to compute JSR of nonnegative matrices
by Algorithm~\textbf{(R)} or by~\textbf{(C)}.

In Section~\ref{sec:CFT} we formulate one of the main results of the paper. This is a
criterion of terminating of Algorithms~\textbf{(R)}, \textbf{(C)},
and~\textbf{(P)} within finite time (Theorem~\ref{th.cond-r}). It
shows that an algorithm produces and extremal polytope and finds
the precise values of JSR if and only if the family~$\cM$ has a
{\em dominant product} (Definition~\ref{d.simple}). In
particular, if the algorithm terminates within finite time, then~$\Pi$
is a dominant product for~$\cM$.

In Section~\ref{sec:LSR} we extend
our method to  the lower spectral radius computation. To this end
we first define an antinorm, prove several theoretical results
about it, and then describe Algorithm~\textbf{(L)} for the exact
computation of LSR of nonnegative matrices. Its practical
efficiency for randomly generated matrices (Section~\ref{sec:numer}) is approximately
 the same as for Algorithm~\textbf{(P)}.

Section~\ref{sec:IE} presents two detailed examples in dimension~$2$ to illustrate
the algorithms.

In Section~\ref{sec:Appl} we consider
applications to several problems of combinatorics, coding theory and number theory.
In the problem of asymptotic growth of the number of overlap-free words (\S \ref{subsec:OFW})
we compute precise values of exponents of the upper and lower growth. This proves two
conjectures stated in 2008~\cite{JPB2}. Then in \S  \ref{subsec:PR} we do the same for the problem of density of ones
in the Pascal rhombus, and disprove one previously know conjecture. In \S  \ref{subsec:EPF} \, and \, \S  \ref{subsec:ETF}
we find precise values of the lower and upper growths of the Euler partition functions for some values
of the parameters.

Section~\ref{sec:numer} presents the results of numerical tests for JSR and LSR computation for randomly generated matrices
of dimensions from $5$ to $100$. In all the cases the algorithms find the exact values
of JSR and LSR, which suggests that our approach generically has finite convergence.
\smallskip

In the sequel we assume that the basis $\{e_i\}_{i=1}^d$ of the space~$\re^d$ is fixed
 and do not distinguish between  operators and the corresponding matrices.
An eigenvalue $\lambda$ is simple, if it is of multiplicity~$1$. The largest by modulo
eigenvalue of an operator $B$ is called {\em leading} and denoted by $\lambda_{\max}$
(if there are several such eigenvalues, then each of them is leading). We use the following notation:
$B^*$ is the operator adjoint to $B$, $\, {\rm int} \, M$ is the interior of a set~$M$, $\, {\rm co} \, (M)$
is the convex hull of~$M$. We use the short abbreviation ``LP'' for linear programming problems.

\section{ Computing of the joint spectral radius: the case of real leading eigenvectors \textbf{(R)}}
\label{sec:JSR}

In this section we present Algorithm~\textbf{(R)} for JSR computation.

We consider an irreducible  family $\cM = \{A_1, \ldots ,
A_m\}$. For some (as large as possible)~$l$ we look over all
products $\Pi$ of length $\le l$ and take one with the biggest
value~$[\rho (\Pi)]^{1/n}$, where $n$ is the length of the
product. We denote it as~$\Pi= A_{d_n}\cdots A_{d_1}$.

Let~$\widetilde M \, = \, \{\widetilde A_1,
\ldots , \widetilde A_m\}$ be the normalized family, where
 $\widetilde A_i \, = \, [\rho(\Pi)]^{\, -1/n}\, A_i$.
For the product $\widetilde \Pi = \widetilde A_{d_n}\cdots \widetilde A_{d_1}$
we have~$\rho(\widetilde \Pi)\, = \, 1$ which implies $\widehat \rho (\widetilde \cM) \ge 1$.
\smallskip

Define, for an arbitrary nonzero vector $v\in \re^d$ the set
\begin{equation}
{\Omega}(v) = \bigcup\limits_{k \ge 0}
\Big\{ \Gamma\,v \ \mid \Gamma \in \widetilde\cM^k \Big\},
\label{eq:traj}
\end{equation}
(where $\widetilde \cM^0 = Id$, the identity matrix),
i.e. the set obtained by joining $v$ to all vectors obtained by applying
the products of the semigroup of $\widetilde {\mathcal M}$ to  $v$.
The following theorem (see \cite{P1} and \cite{GZ1}) relates the set
${\Omega}(v)$ and an extremal norm for $\widetilde \cM$.

\begin{theorem}
Let~$\widetilde M \, = \, \{\widetilde A_1,\ldots , \widetilde A_m\}$ be irreducible and
such that $\widehat \rho (\widetilde \cM) \ge 1$ and let ${\Omega}(v)$ (for a given $v \ne 0$)
be a bounded subset of $\re^d$ spanning $\re^d$.
Then $\widehat \rho (\cM) = 1$.
Furthermore the set
\begin{equation}
\overline{{\rm co}_s\left({\Omega}(v) \right)} =
\overline{{\rm co}\left({\Omega}(v) \cup -{\Omega}(v) \right)}
\label{eq:ball}
\end{equation}
is the unit ball of an extremal norm $\|\cdot \|$ for $\widetilde M$ (and for $\cM$).
\label{th:extnorm}
\end{theorem}

The main idea of the algorithm we present is to finitely compute the set (\ref{eq:ball})
whenever it is a polytope. Let us clarify this key point.

We say that a bounded set $P\subset \re^d$ is a {\em balanced real polytope} (b.r.p.) if there
exists a finite set of vectors ${\mathcal V}=\{v_{i}\}_{1\leq i \leq p}$ (with $p \ge d$)
such that ${\rm span}({\mathcal V})=\re^d$ and
\begin{equation}
P = {\rm co}_s (\cV) = {\rm co}(\cV,-\cV).
\label{eq:rpn}
\end{equation}
Therefore
$$
P = \Big\{ z = \sum\limits_{x \in \cV}\, t_x\, x \ \ \ {\rm with} \ \ \ -q_x \le t_x \le q_x, \quad q_x \ge  0 \ \forall x \in \cV
\ \ \ {\rm and} \ \ \ \sum\limits_{x \in \cV}\, q_x \le 1 \Big\}.
$$
The set $P$ is the unit ball of a norm $\|\cdot \|_{P}$ on
$\re^d$, which we call a {\em real polytope norm}.

Assume that the hypotheses of Theorem \ref{th:extnorm} hold. The
possibility of actually determining an extremal polytope norm, if
any, crucially relies on the search of the initial vector $v$, which we will
address later in Theorem \ref{th.cond-r} that suggests to choose $v$
as a leading eigenvector of $\Pi$ (although a different choice would be
admissible). We will assume here $v$ to be real.

The idea is that of computing the set ${\Omega}(v)$ by
applying recursively the family $\widetilde M$ to a finite set
of vectors (which in the beginning is simply the vector $v$),
checking at every iteration $h$ whether $\widetilde M$ maps the symmetrized convex
hull (${\rm co}_s ( {\Omega}^{h-1}(v) )$) of the computed set of
vectors
$$
{\Omega}^{h-1}(v) = \bigcup\limits_{0 \le k \le h-1}
\Big\{ \Gamma\,v \ \mid \Gamma \in \widetilde\cM^k \Big\},
$$
into itself.

Algorithm \textbf{(R)} we are going to present is similar to the one described in \cite{GZ1}; the main differences
are that a new vertex is included even if it lies on the boundary of the current polytope,
all the leading eigenvectors are considered as starting vertices of the
searched extremal polytope and a new and efficient stopping criterion is
added.

We start by an auxiliary result and then describe the algorithm.

\begin{lemma}\label{l.5}
Let an operator $B$ have a unique simple  leading eigenvalue
$\lambda \in \re$ with the leading eigenvector $v$; let also $v^*$
be the leading eigenvector of $B^*$ such that $(v^*, v) = 1 $. If
for some operator $C$ one has~$\bigl| (v^*, Cv)\bigr|\, > \, 1$,
then for sufficiently large $r$ the  operator $B^rC$ has a unique
simple  leading eigenvalue, which is real and bigger
than~$\lambda$ by modulo.
\end{lemma}
{\tt Proof}. Without loss of generality, after a suitable normalization, it can be assumed
that~$\lambda = 1$. Since all other eigenvalues of $B$ are smaller
by modulo than~$1$, it follows that $B^r$ converges to the
one-rank operator $B_{\infty}\, x \, = \, (v^*, x)\, v\, $ as $\,
r \to \infty$. Hence $B^rC$ converges to the operator
$B_{\infty}C$, whose unique simple leading eigenvalue is $(v^*,
Cv)$, which exceeds~$1$ by modulo.
{\hfill $\Box$}
\smallskip

\begin{remark}\label{r.5}
{\rm If $v$ and $v^*$ are the leading eigenvectors of $B$ and $B^*$ respectively,
then these vectors cannot be orthogonal, otherwise the leading eigenvalue is not simple.
So, $(v^*, v)\, \ne \, 0$, and hence, after a suitable normalization it can always be assumed that
$(v^*, v) = 1$.
}
\end{remark}

\smallskip




It is well-known that the problem of JSR computation has to be considered
only for irreducible families of matrices, which do not possess common invariant
linear subspaces.  Otherwise this family is factorable in a suitable basis in~$\re^d$: all the matrices~$A_j$
get a block upper-triangular form, and $\widehat \rho (\cM)$ equals to the maximal JSR of the blocks.
This reduces the problem of JSR computation to several problems in smaller dimensions.
Therefore, in the sequel of this section we assume that~$\cM$ is irreducible.

\smallskip

\subsection{Algorithm~\textbf{(R)}}
\label{subsec:R}

{\tt Initialization.} Given the irreducible  family $\cM = \{A_1, \ldots ,
A_m\}$ we look over all products $\Pi$ of length $\le l$ and consider
the shortest product $\Pi$ such that $[\rho (\Pi)]^{1/n}$ is maximal,
where $n$ is the length of the product. We denote it as~$\Pi= A_{d_n}\cdots A_{d_1}$
and consider the main assumption:
\smallskip

{\em (i) The product~$\Pi$ has a real nonzero leading eigenvalue.}
\smallskip

\noindent We assume that the leading eigenvalue $\lambda_{\,
\max}$ of~$\Pi$ is positive; the case of negative eigenvalue is
considered in the same way. Let~$\widetilde M \, = \, \{\widetilde A_1,
\ldots , \widetilde A_m\}$ be the normalized family, where
 $\widetilde A_i \, = \, [\rho(\Pi)]^{\, -1/n}\, A_i$.
For the product $\widetilde \Pi = \widetilde A_{d_n}\cdots \widetilde A_{d_1}$
we have~$\lambda_{\max}\, = \, 1$.
\smallskip

Let $\widetilde \Pi_1 = \widetilde \Pi\, , \, \widetilde \Pi_i = \widetilde
A_{d_{i-1}} \cdots \widetilde A_{d_1}\widetilde A_{d_n}\cdots \widetilde
A_{d_i}$ be a cyclic permutation of $\widetilde \Pi_1$, $\, i = 2,
\ldots , n$. We denote by $v_1$ the leading eigenvector of~$\widetilde
\Pi_1$, for which $\widetilde \Pi_1v_1 = v_1$. If it is not unique (in
which case $\lambda_{\, \max}$ is multiple) we take any of them.
Then for every $i\ge 2$ we set
$$
v_i \quad = \quad \widetilde A_{d_{i-1}}\cdots \widetilde A_{d_{1}}v_1 \ .
$$
Thus, $v_i$ is a leading eigenvector of $\widetilde \Pi_i$.
\smallskip

In case $\widetilde \Pi$ has a unique simple eigenvalue, we also need
the corresponding dual system of vectors:
$v_1^*$ the leading eigenvector of the conjugate operator~$\widetilde \Pi_1^*$
normalized by the condition $(v_1^*, v_1) = 1$ (see Remark~\ref{r.5}), and
$$
v_i^* \quad = \quad \widetilde A_{d_{i}}^*\cdots \widetilde
A_{d_{n}}^*v_1^*\,
$$
for $i = 2, \ldots , n$. Thus, $v_i^*$ is the leading eigenvector of $\widetilde \Pi_i$,
and $(v_i^*, v_i) = 1$.
If the leading eigenvalue of $\widetilde \Pi$ is multiple or not unique, then we
do not need the conjugate system.
\medskip

{\tt Set $k=0$}.  We set $\cV_0\, = \, \cU_{\, 0}\, = \, \{v_1, \ldots
, v_n \}\, $ and $\, \cR_0\, = \, \bigl\{ \, ( v_i\, , \, \widetilde
A_p) \bigl|\ i = 1, \ldots , n\, ; \ p=1, \ldots , m \, , \  p\ne
d_i \bigr\}$.
\medskip

{\tt Main loop}
\smallskip

{\tt For $k\, \ge \, 1$}. We have  finite sets $\cV_{k-1}\,
\subset \, \re^d\, ,\,  \cU_{k-1} \, \subset \, \cV_{k-1}$, and
$\cR_{k-1}\, \subset \, \cU_{k-1}\times \widetilde \cM$. Set
$\cV_{k}\, = \, \cV_{k-1}\, $ and $\, \cU_k \, = \, \emptyset$.
Take an arbitrary pair $(v, \widetilde A)\, \in \, \cR_{k-1}$ and
compute the norm whose unit ball is the polytope with vertices
$\cV_{k-1}$ and $-\cV_{k-1}$, of the corresponding vector
$z = \widetilde A v$. This is done by
solving the following LP problem with variables $\{t_x\}_{x \in
\cV_k}\, , \, \{q_x\}_{x \in \cV_k}$  and  $t_0$ (which represents the
reciprocal of the value of the norm):
\begin{equation}\label{LP.r}
\left\{
\begin{array}{rcl}
\max & & t_0
\\[0.25cm]
{\rm subject \ to}
     & & t_0\,  z = \sum\limits_{x \in \cV_{k}}\, t_x\, x
\\[0.3cm]
     & & -q_x \le t_x \le q_x, \quad q_x \ge  0 \qquad \forall x \in \cV_k
\\[0.3cm]
{\rm and}
     & & \sum\limits_{x \in \cV_{k}}\, q_x \le 1,
\end{array}
\right.
\end{equation}
The value of the problem, i.e., the value $\, \max \, t_{0}$ will be denoted by
$t_{\{v, \widetilde A\}}$. Thus, for a given pair $(v, \widetilde A) \,
\in \ \cR_{k-1} \, $ we obtain  the value $\, t_{\{v, \widetilde
A\}}$.
\smallskip

{\tt If} $\, t_{\{v, \widetilde A\}} \, > \, 1$,
then we leave the sets $\cV_k  $ and $\cU_k$ as they are, take the next
pair~$(v, \widetilde A)\, \in \, \cR_{\, k-1}$
and consider problem~(\ref{LP.r}) for it.
\smallskip

\ {\tt If}  $\, t_{\{v, \widetilde A\}} \, \le \, 1$,  then we distinguish between two cases
\medskip

\ \ {\tt If} \ the leading eigenvalue of~$\Pi$
is  unique and simple, we apply the following
\smallskip

\ \ {\tt Stopping criterion:}
\smallskip

\ \ For a given pair $(v, \widetilde A)$ we check the condition
\begin{equation}\label{cond.r}
\bigl| \, \bigl(v_j^*\, , \, \widetilde A\, v \bigr)\, \bigr|\quad \le \quad 1\ , \qquad
j \, = \, 1, \ldots ,n\, .
\end{equation}

\ \ \ {\tt If} (\ref{cond.r}) is satisfied, then we set $\,  \cV_k \, = \,
\cV_k\cup \{\, \widetilde A\, v\}\, ,
 \cU_k \, = \, \cU_k\cup \{ \widetilde A\, v\}$,
 take the next pair~$(v, \widetilde A)\, \in \, \cR_{k-1}$
and consider problem~(\ref{LP.r}) for it.

\ \ \ {\tt Otherwise If} (\ref{cond.r})  is not satisfied, then $\Pi$ is not an s.m.p. for $\cM$, and
$\widehat \rho (\widetilde \cM)> 1$ (Lemma~\ref{l.5}).
We stop the algorithm and go either to the {\tt Final step}, or back to the
{\tt Initialization.} In the latter case we need to find another candidate s.m.p.
The first option is to increase $l$ and to look over all products of a bigger length.
Lemma~\ref{l.5} provides also a different approach.  We take an index $j$, for which
$\bigl| \, (v_j^*\, , \, \widetilde A\, v )\, \bigr|\, > \,  1$. Applying Lemma~\ref{l.5} for the
vectors $v_j^*$ and $ -\, v_j^*$, we conclude that there is~$r$ such that
$\lambda_{\max}(\widetilde \Pi_j^r\, \widetilde A_{s_q}\cdots \widetilde A_{s_1})\, > \, 1$,
where $\widetilde A_{s_q}\cdots \widetilde A_{s_1} v_j \, = \, \widetilde A v$. We take the new initial product
$\Pi = \Pi_j^r\,  A_{s_q}\cdots  A_{s_1} $ and restart the
algorithm.
\smallskip

\ \ \ {\tt End If}
\smallskip

\ \ {\tt Otherwise If}\  the leading eigenvalue of $\Pi$ is not unique or
multiple, then we do not apply the stopping criterion, and
 set $\,  \cV_k \, = \,  \cV_k\cup \{ \, \widetilde A\, v\}\, ,
 \cU_k \, = \, \cU_k\cup \{ \widetilde A\, v\}$,
 take the next pair~$(v, \widetilde A)\, \in \, \cR_{k-1}$
and consider problem~(\ref{LP.r}) for it.

\smallskip
\ \ {\tt End If}
\medskip

{\em The $k$th step is over when the whole set $\cR_{k-1}$ is
exhausted.}
\medskip

{\tt If} $\, \cU_k = \emptyset$,
then $\widehat \rho (\widetilde \cM) = 1$, and so $\widehat \rho (\cM)\, = \,
[\rho(\Pi)]^{1/n}$. The extremal polytope is $P_{k-1} = {\rm
co}_{s}\, (\cV_{k-1})$, and the s.m.p. for $ \cM$ is $ \Pi$. The
algorithm terminates having performed~$k$ steps.
\smallskip

{\tt Otherwise If} $\cU_k \ne \emptyset$, then we set $\cR_{k}\, = \, \cU_k \times \widetilde \cM$ and continue.
\smallskip

{\tt End If}

\smallskip
{\tt End For}

\bigskip

{\tt Final step.} If the algorithm has not terminated, then we
stop it after some $N$ steps, denote $t_N\, = \, \,
\min\limits_{(v, \widetilde A)\in \cR_{N-1}}\, t_{\{v, \widetilde A\}}$,
where $t_{\{v, \widetilde A\}}$ is the solution of LP
problem~(\ref{LP.r}) for the last step, i.e., for $k=N$, and have
the following estimate for the joint spectral radius of the family
$\cM$:
\begin{equation}\label{final.r}
[\rho (\Pi)]^{1/n}\quad \le \quad \widehat \rho (\cM) \quad \le \quad
t_N^{-1}\, [\rho (\Pi)]^{1/n}  \,  .
\end{equation}
\medskip

{\tt End of Algorithm \textbf{(R)}}.
\medskip

\begin{remark} \label{rem:diff}
{\rm An important difference with respect to previous similar algorithms
is that if at step $k$ a new vector $v$ lies on the boundary of the
polytope $P_{k-1} = {\rm co}_{s}\, (\cV_{k-1})$ (which means
$t_{\{v, \widetilde A_p\}} = 1$ for some $p$)
then we include the vector as a new vertex.
Clearly this condition is non generic and requires - to be tested
in floating point arithmetics - the use of a suitable error
tolerance.}
\end{remark}

\noindent Before we give the proofs, let us explain the general
scheme of the algorithm.
\smallskip

\subsection{The cyclic tree structure of the algorithm}
\label{subsec:CT}

Consider a combinatorial {\em cyclic tree} $\cT$ defined as
follows. The root is formed by a cycle $\mathbf{B}$ of $n$
nodes~$v_1, \ldots , v_n$. They are, by definition, the nodes of
zero level. For every $i\le n$ an edge (all edges are directed)
goes from $v_i$ to $v_{i+1}$, where we set $v_{i+1} = v_1$. At each
node of the root $m-1$ edges start to nodes of the first level.
So, there are $n(m-1)$ different nodes on the first level. The
sequel is by induction: there are $n(m-1)m^{k-1}$ nodes of the
$k$th level, $k\ge 1$, from each of them $m$ edges (``children'')
go to $m$ different nodes of the $(k+1)$st level.

Consider now an arbitrary  word $\mathbf{b} = d_n\ldots d_1$ of
length $n\ge 1$, where each $d_j$ belongs to the alphabet $\{1,
\ldots , m\}$. The product of several words is their
concatenation. We assume that $\mathbf{b}$ is irreducible, i.e.,
is not a power of a shorter word. To every edge of the tree~$\cT$ we
associate a letter $d$ as follows: the edge $v_iv_{i+1}$
corresponds to $d_i\, , \, i=1, \ldots , n$; at each node $m$
edges start associated to $m$ different letters. To a given  word
$q_k\ldots q_1$ we associate  the node, which is the end of the
path from $v_1$ along the edges $q_1, q_2, \ldots , q_k$. For
example, the empty word corresponds to $v_1$, the
word~$\mathbf{b}$ also corresponds to $v_1$, the word $d_2d_1$
corresponds to $v_3$, the word $d_2$ corresponds to either $v_2$,
if $d_2 = d_1$, or to a child of $v_1$ from the first level,
otherwise. This tree is said to be generated by the word
$\mathbf{b}$, or by the cycle $\mathbf{B}$.

For a family of operators $\widetilde \cM \, = \, \{\widetilde A_1, \ldots
, \widetilde A_m\}$ and for some product $\widetilde \Pi \, = \, \widetilde
A_{d_n}\cdots \widetilde A_{d_1}$ with an eigenvalue~$1$ we associate
the cyclic tree $\cT$ generated by the word $d_n\ldots d_1$. The
node $v_1$ corresponds to an  eigenvector with the eigenvalue~$1$;
to a given node $v \in \cT$ we associate a point $\widetilde A_{q_k}\ldots
\widetilde A_{q_1}v_1$, where the word $q_k\ldots q_1$ corresponds to
the node $v_1$.

When we start the algorithm, we take the set $\mathbf{B} \, = \,
\{v_1, \ldots, v_n\}$ as the root of the tree. At the first step
we take any node $v_i$ and consider successively  its $(m-1)$
children from the first level. For each neighbor  $\, u \, = \,
\widetilde A v_i$, where $\widetilde A \in \widetilde M \setminus \{\widetilde
A_{d_i}\}$ we solve LP problem~(\ref{LP.r}) and determine, whether
$u$ belongs to the interior of the set ${\rm co}_s\, (\cV_1)$, where
${\rm co}_{s}(M)\, = \, {\rm co} \{M, -M\}$ is the symmetrized convex hull.
If it does, then $u$ is a ``dead leaf'' generating a ``dead branch'':
we will never come back to $u$, nor to nodes of the branch
starting at $u$ (so, this branch is cut off). If it does not, then
$u$ is an ``alive leaf'', and we add this element~$u$ to the set $\cV_1$
and to the set $\cU_1$. After the first step all alive
leaves of the first level form the set~$\cU_1$.  At the second
step we deal with the leaves from $\cU_1$  only and obtain the next set
of alive leaves of the second level~$\cU_2$, etc. Thus, after the
$k$th step we have a family $\cU_k$ of alive leaves from the $k$th
level, and a set $\cV_k , = \, \cup_{j=0}^{k}\, \cU_j$. A node~$u$
belongs to~$\cV_k$ iff its level does not exceed~$k$ and it
belongs to an alive branch starting from the root. The polytope
$P_k$ is the symmetrized convex hull ${\rm co}_s \, (\cV_k)$.
The polytope $P_{k-1}$ is extremal iff $\cU_{k}\, = \, \emptyset$,
i.e., the $k$th step produces no alive leaves (only dead ones).
This means that there are no alive paths of length~$k$ from
the root. Therefore $P_{k} = P_{k-1}$. Otherwise, if $\cU_{k}$ is nonempty, we make the next
step and go to the $(k+1)$st level: take children of each element
of $\cU_k$, determine whether they are alive or dead and proceed.

\smallskip

\subsection{Explanations and proofs}
\label{subsec:EP}

 The algorithm produces a
sequence of embedded polytopes $P_1 \subset
 P_2 \subset \ldots $
such that $P_{j+1} \, = \, {\rm co}_s\, \bigl\{ \widetilde A_1P_j,
\ldots , \widetilde A_mP_j\bigr\}\, $ for every~$j$. If the algorithm
terminates after the $k$th step, then $P_{k} = P_{k-1}$. The $k$th
step is actually needed only to ensure that the polytope $P_{k-1}$
is extremal, i.e., $\widetilde A_j\, P_{k-1} \subset P_{k-1}\, , \
j=1, \ldots , m$. In this case $P_{k-1}$ possesses an interior of nonzero
measure, otherwise its linear span is a common invariant
nontrivial subspace of the family~$\widetilde M$, which contradicts
the irreducibility assumption. Moreover, $P_{k-1}$ is
centrally-symmetric, hence $0 \in {\rm int}\, P_{k-1}$. This, in
particular, yields that if for some $v \in \re^d$ and $\, t > 1$
one has $\, t\, v \, \in \, P_{k-1}$, then $\, v \, \in \, {\rm
int}\, P_{k-1}$. Thus, if the value $t_{\, \{v, \widetilde A
\}}\, $ of LP problem~(\ref{LP.r}) is bigger than~$1$, then
$\widetilde Av\, \in \, {\rm int}\, P_{k-1}$. Thus all dead
leaves removed by the algorithm are internal points
for~$P_{k-1}$.

In the Minkowski norm $\|\cdot \|_{k-1}$ whose unit ball is given by $P_{k-1}$,
one has $\|\widetilde A\|_{k-1} \le 1$ for all $\, \widetilde A\in \widetilde M$,
therefore $\widehat \rho (\widetilde \cM) \le 1$.
On the other hand, $\widehat \rho (\widetilde M) \,\ge \, \rho(\widetilde \Pi)^{1/n}\, = \, 1$,
hence  $\widehat \rho (\widetilde \cM) \, = \, \rho(\widetilde \Pi)^{1/n}\, = \, 1$, an so
$\widehat \rho ( \cM) \, = \, \rho(\Pi)^{1/n}$. Thus, if the algorithm terminates within
finite time, then the s.m.p. and the exact value of JSR are
found.
\smallskip

Suppose the algorithm does not terminate within finitely many
steps. After the final step we take the polytope $P_{N-1} = {\rm
co}_s\, (\cV_{N-1})$ as a unit ball of the new norm $\|\, \cdot \,
\|_{N-1}$ in $\, \re^d$. Then $\, \max\limits_{A \in \cM}\,
\|A\|_{N-1}\, \ge \, \widehat \rho(\cM)$. We have
$$\, \max\limits_{A \in
\cM}\, \|A\|_{N-1}\quad = \quad \bigl[
\rho(\Pi)^{1/n}\bigr]\max\limits_{\widetilde A \in \widetilde \cM}\,
\|\widetilde A\|_{N-1} \quad = \quad \bigl[ \rho(\Pi)^{1/n}\bigr]\cdot
 \max\limits_{v \in \cU_{N-1}} (t_{\{v, \widetilde A\}})^{-1}\, ,
$$
where $t_{\{v, \widetilde A\}}$ is the value of LP
problem~(\ref{LP.r}) for $k=N$. Therefore, $\widehat \rho (\widetilde \cM)
\, \le \, (t_N)^{-1}$, and after multiplying by $
[\rho(\Pi)]^{1/n}$ we arrive at~(\ref{final.r}).
\bigskip

\begin{remark}\label{r290}
{\rm
Although we perform operations numerically, the obtained results
have to be considered exact since apart from the vertices
of~$P_{k-1}$, all other vectors obtained by applying the scaled
matrices $\widetilde A_j$ to the vertices are either vertices or
internal points to the polytope~$P_{k-1}$.}
\end{remark}
\smallskip

\smallskip
\begin{remark}\label{r300}
{\rm By the construction of the algorithm, each  vertex  of the polytope~$P_k$
belongs either to  $\cV_k$ or to $-\cV_k$. However, not all elements of the set~$\cV_k$
are actually vertices: some of them may lie in the convex hull of the others.
This means that in general the set $\cV_k$ is not an essential system of
vertices.
Nevertheless, for the sake of simplicity we call all elements of~$\cV_k$ vertices.
}
\end{remark}

\begin{remark}\label{r310}
{\rm Actually Algorithm~\textbf{(R)} can be applied to a reducible family~$\cM$ as well.
If the algorithm terminates after $k$th iteration, and the set $\cV_k \subset \re^d$ does not lie in
a linear subspace of a smaller dimension (i.e., the system of equations $(x, v)=0\, , \ v \in \cV_k$ has only trivial solution
$x=0$), then $P_{k-1}$ is an extremal polytope, and $\Pi$ is an s.m.p. Thus, one can apply
Algorithm~\textbf{(R)} without preliminary checking of irreducibility of~$\cM$. Nevertheless,
if the family $\cM$ is reducible, then it is always better  to factorize~$\cM$ before applying Algorithm~\textbf{(R)},
because this reduces the dimension of matrices.
}
\end{remark}
\smallskip

\subsection{Efficiency results for Algorithm~\textbf{(R)}}
\label{subsec:ER}

If the algorithm terminates, then it finds the exact value
of JSR, otherwise estimate~(\ref{final.r}) gives its approximate
value with the relative error $\, \varepsilon \, = \, (t_N)^{-1} -
1$. This error depends on two integer parameters:
the maximal length~$l$ of the products, among which we choose an s.m.p. $\Pi$, and
the number of iterations~$N$ of the algorithm. Let us show that the error~$\varepsilon$
tends to zero as both these parameters increase:
\begin{prop}\label{p2.r}
For an arbitrary irreducible family $\cM$ we have  $t_N \to 1$ as  $l \to
\infty$ and $N \to \infty$.
\end{prop}
Thus, both sides of inequality~(\ref{final.r}) tend to $\widehat \rho (\cM)$ as $l \to \infty\, , \, N \to \infty$.
Algorithm~\textbf{(R)} either finds the value of JSR or provides lower and upper bounds for it; those bounds are arbitrarily
close to each other, whenever both $l$ and $N$ are large enough.

In the proof we use Dini's theorem on monotone convergence:  if a sequence of continuous real-valued functions defined on a compact metric space~$Q$ is monotone and converges pointwise to a continuous function,  then this convergence is uniform on~$Q$ (see~\cite[theorem 7.13]{Ru}). We use the  Minkowski norm $\|\cdot\|_{\, D}$ associated to a given symmetric convex body
$D \subset \re^d$ as follows:  $\, \|\cdot\|_{\, D} \, = \, \inf \, \bigl\{ t^{-1}  \ \bigl| \ t >0, \  t x \in D  \bigr\}$.

\noindent {\tt Proof of Proposition~\ref{p2.r}}. Assume first that $\widehat \rho (\widetilde \cM)= 1$, i.e., that $\Pi$ is an s.m.p.
The algorithm produces the polytopes
$\{P_k\}_{k \in \n}$ such that
$$P_{k} \, \subset \, P_{k+1} \, = \,  {\rm co}_{s}\, \bigl\{ \widetilde A_1P_k, \ldots ,  \widetilde A_mP_k\bigr\}.$$
 Since the family~$\widetilde \cM$ is irreducible, there is $p \ge 1$ such that
 all the polytopes $P_k$ have nonempty interior
 for $k \ge p$. Hence, for $k\ge p$ the polytope $P_k$ generates the Minkowski norm
 $f_k(\cdot )\, = \, \|\cdot \|_{\, P_k}$. For each $x \in  \re^d$ the sequence $\{f_k(x)\}_{k \ge p}$
 is non-increasing. Moreover, it is uniformly bounded below by a positive constant, because all the polytopes $\{P_k\}_{k \ge p}$
 are contained in some ball, since the family $\widetilde \cM$ is non-defective. Therefore, the sequence
 $f_k(x)$ converges pointwise to a function $f(x)$, which is also a norm in~$\re^d$. By Dini's theorem,
 this convergence is uniform on any compact subset of~$\re^d$. In, particular, it is on the unit sphere
 $S = \{ x \in \re^d \ | \ f(x) = 1\}$ of the norm~$f$. Thus, $f_k(x) \to 1$ uniformly for $x \in S$, as $k \to \infty$. Hence, there is $N_{\varepsilon}$ such that $f_{N-1}(x) \le 1+\varepsilon$ for all $x \in S$, whenever $N \ge N_{\varepsilon}$.
 Consequently, $(t_N)^{-1} \, = \, \sup_{x \in P_N}f_{N-1}(x) \, \le \, \sup_{x \in S}f_{N-1}(x)\, \le \, 1+\varepsilon$,
 which completes the proof for the case $\widehat \rho = 1$. Consider now the general case. We have
 $[\rho (\Pi)]^{1/n} \to \widehat \rho (\cM)$ as $l \to \infty$, where, let us remember, $n = n(l)$ is the length
  of $\Pi$. Hence, for every $\delta > 0$ there is $l_{\delta}$ such that
 $\widehat \rho (\widetilde \cM) < 1 + \delta$. Since each polytope $P_k$ continuously depends on the family $\widetilde \cM  =
 [\rho (\Pi)]^{-1/n} \cM$, for all sufficiently small $\delta$ one has $(t_{N_{\varepsilon}})^{-1} <   1 + \varepsilon$.
 This inequality holds for all $l \ge l_{\delta}$. It remains to note that for every family $\cM$ the value $t_N$
 is non-decreasing in~$N$. Indeed, $t_{N}P_N \subset P_{N-1}$, hence $t_{N}\widetilde A_jP_{N} \subset \widetilde A_j P_{N-1}$
 for every $j = 1, \ldots , m$, and so $t_{N}P_{N+1} \subset P_{N}$. Therefore,
 $t_{N+1} = \sup \, \{t > 0 \ | \ t\, P_{N+1} \subset P_{N}\}\, \ge \, t_{N}$. We see that
 $(t_{N})^{-1}  < 1+ \varepsilon$,  provided $N \ge N_{\varepsilon}$. Thus, for every $\varepsilon > 0$
 there are $l_{\delta}$ and $N_{\varepsilon}$ such that
 $(t_{N})^{-1} - 1 < \varepsilon$, whenever $N \ge N_{\varepsilon}$ and $l \ge l_{\delta}$.

{\hfill $\Box$}
\smallskip

\bigskip

Let us now show the efficiency of the stopping criterion.
\begin{prop}\label{p1.r}
Assume the leading eigenvalue
of~$\Pi$ is real, unique and simple.  If the assumption of the
algorithm is wrong (i.e., $\Pi$ is not an s.m.p.)  then for every
$j = 1, \ldots , n$ condition~(\ref{cond.r}) is violated at some
step. Conversely, if condition~(\ref{cond.r}) is violated at some
step for some~$j$, then $\Pi$ is not an s.m.p.
\end{prop}
{\tt Proof.} The sufficiency follows from Lemma~\ref{l.5}. To prove the necessity
suppose $\widehat \rho (\widetilde M) > 1$; then for every point $v \ne 0$ and for every number
$R > 0$ there is a product $C$ of operators of the family $\widetilde M$ such that $\|Cv\|\ge R$~\cite[theorem 1]{P1}.
On the other hand, since the family $\widetilde M$ is irreducible, it follows that there is $\gamma > 0$
such that for every point $y \ne 0$ the set $K_d(y)\, = \,
{\rm co}_s\, \{\, A y\ | \  A \in \widetilde M^d\}$ contains a ball of radius
$\gamma \|y\|$ (see~\cite{K2} for the proof), where $\gamma > 0$ is a constant.
Applying these results to the points $v =v_1$ and $y = Cv_1$
and using the fact that the polytope $P_{s+d}$ contains $K_d(y)$, where $s$ is the length of the product $C$,
we see that  the polytope $P_{s+d}$ contains a ball of radius $\gamma R$ centered at the
origin.  Therefore,
$$
\sup_{x \in P_{s+d}}|(v_j^*, x)| \ \ge \ \gamma\, R\, \|v_j^*\|.
$$
Since this supremum is attained at some
vertex of~$P_{s+d}$, which is produced by the algorithm, we see that condition~(\ref{cond.r})
will fail by the $(s+d)$-th step, whenever $R > 1/(\gamma \|v_j\|^*)$.

{\hfill $\Box$}
\smallskip

Thus, if the chosen product $\Pi$ is not an s.m.p., then the stopping criterion always
determines this in finite time.
In Section~\ref{sec:CFT} we formulate Theorem~\ref{th.cond-r} that gives a
sharp criterion for the algorithm to terminate in finitely many
steps (and, respectively, to produce an extremal polytope).


\section{ Computing of the joint spectral radius: the case of complex leading eigenvectors \textbf{(C)}}
\label{sec:C}

For the theoretical results and the algorithms relevant to this case we mainly address the reader to the papers
\cite{GWZ,GZ2,GZ3}.

We recall from \cite{GZ3,VZ} the definition of a balanced complex polytope, which generalizes to the complex case
a centrally symmetric real polytope.

Let ${\mathcal V}=\{v_{i}\}_{1\leq i\leq p}$ be a finite set of vectors,
then
\begin{equation}
{\rm absco}({\mathcal V})=\Big\{z\in \co^d \ \Big| \ z=\sum_{x \in \cV}
t_x\,x \quad {\rm with}
\quad \sum_{x \in \cV} |t_x|\leq 1\Big\}.
\label{160}
\end{equation}

\begin{definition}
A set ${\mathcal P}\subset \co^d$ is a balanced complex polytope (b.c.p.)
if there exists a finite set of vectors ${\mathcal V}=\{v_{i}\}_{1\leq i \leq p}$
such that

\begin{equation}
{\rm span}({\mathcal V})=\co^d \qquad \mbox{and} \qquad {\mathcal P}={\rm absco}({\mathcal V}).
\label{103}
\end{equation}
Moreover, if ${\rm absco}({\mathcal V}') \varsubsetneq {\rm
absco}({\mathcal V})$ for all ${\mathcal V}' \varsubsetneq {\mathcal V}$,
we say that ${\mathcal V}$ is an essential system of vertices
for ${\mathcal P}$. Every vector $u\,v_{i}$ with $u \in \co$,
$|u|=1$, is called a vertex of ${\mathcal P}$.
\label{definition13}
\end{definition}

Note that geometrically a b.c.p. ${\mathcal P}$ is not a
classical polytope (see \cite{GZ3}).

A polytope norm can be defined in a natural way.
\begin{lemma}
Any b.c.p. ${\mathcal P}$ is the unit ball of a norm $\|\cdot \|_{{\mathcal P}}$ on
$\co^d$.
\label{lemma1}
\end{lemma}
The proof is immediate (see e.g. \cite{GZ3})
\begin{definition}
We shall call {\em complex polytope norm} any norm $\|\cdot \|_{{\mathcal P}}$
whose unit ball is a b.c.p. ${\mathcal P}$.
\label{definition21}
\end{definition}
The corresponding vector norm is characterized by the following Lemma
(for a proof see \cite{GZ3}).
\begin{lemma}
Let ${\mathcal P}$ be a b.c.p. and let $\|\cdot \|_{{\mathcal P}}$ be the
corresponding complex polytope norm. Then, for any $z\in \co^d$, it holds
that
\begin{equation}
\|z\|_{{\mathcal P}}= \Big\{ \max t_0  \ \Big| \
t_0 z = \sum_{x \in \cV} t_x \,x, \quad \sum_{x \in \cV} |t_x| \le 1\Big\}, \label{eq:polyv}
\end{equation}
where ${\mathcal V}=\{v_{i}\}_{1\leq i\leq p}$ is an essential system
of vertices for ${\mathcal P}$.
\end{lemma}

Complex polytope norms are dense in the set of all norms defined
on $\co^d$ and consequently the corresponding set of induced matrix
complex polytope norms is dense in the set of all induced
$d\times d$-matrix norms (see \cite{GZ3}). This implies the following
important property:
$$
\widehat\rho(\cM) = \inf\limits_{\| \cdot\|_{\mathcal P}} \max\limits_{A \in \cM} \| A \|_{\mathcal P}
$$
where $\| \cdot \|_{\mathcal P}$ denotes the set of polytope norms.

From an algorithmic point of view, the above property has the consequence that although an
extremal polytope norm may not exist, it is possible to compute a polytope norm which is
$\varepsilon$-close to an extremal one, for any $\varepsilon > 0$.

\smallskip

\subsection{Main differences between Algorithm~\textbf(C) and Algorithm~\textbf(R).}
\label{subsec:MD}

Algorithm~\textbf{(R)} extends to the complex case in a direct way, as well
as the convergence and approximation results (see \cite{GWZ}, \cite{GZ2}).
The only (important) difference lies in the computation of the polytope
norm of a vector. We obtain this by rewriting (\ref{eq:polyv}) as a real
optimization problem.

Let ${\mathcal P}={\rm absco}({\mathcal V})$ (with ${\mathcal V} = \{ v_1,v_2,\ldots,v_p \}$)
be a b.c.p. and $\|\cdot \|_{{\mathcal P}}$ the associated norm.
For any $z\in \co^d$, we write (\ref{eq:polyv}) (with
$t_x = \alpha_x + {\rm i} \beta_x$) in the following way:
\begin{equation}\label{LP.c}
\left\{
\begin{array}{rcl}
\max & & t_0
\\[0.25cm]
{\rm subject \ to}
     & & \sum\limits_{x \in \cV} \alpha_x\,{\rm Re}(x) - \beta_x\,{\rm Im}(x)  = t_0 {\rm Re}{(z)}
\\[0.35cm]
     & & \sum\limits_{x \in \cV} \alpha_x\,{\rm Im}(x) + \beta_x\,{\rm Re}(x)  = t_0 {\rm Im}{(z)}
\\[0.35cm]
{\rm and}
     & & \sum\limits_{x \in \cV} \sqrt{\alpha_x^2 + \beta_x^2} \le 1
\end{array}
\right.
\end{equation}


This problem can be efficiently solved in the framework of
the conic quadratic programming, by the interior point method on Lorentz cones,
see~\cite{AG,ART}. The corresponding pocket of programs can be found in
http://www.mosek.com .


The second difference with respect to Algorithm~\textbf{(R)}  is concerned with the stopping criterion.
In particular (\ref{cond.r}) has to be replaced by the following condition (\ref{cond.c}).
For a given pair $(v, \widetilde A)$ we have to check in fact the condition
\begin{equation}\label{cond.c}
\bigl| {\rm Re} \, \bigl(v_j^*\, , \, \widetilde A\, v \bigr)\, \bigr|\quad \le \quad 1\ , \qquad
j \, = \, 1, \ldots ,p\, .
\end{equation}



\section{ Computing of the joint spectral radius: the case of nonnegative matrices \textbf{(P)}}
\label{sec:P}

Algorithm~\textbf{(R)} can be modified for families of nonnegative matrices
to improve significantly its efficiency. The corresponding Algorithm~\textbf{(P)}
has a very similar structure, but differs from Algorithm~\textbf{(R)} in
several key points. Before describing the algorithm we need to establish several auxiliary
results on operators with an invariant cone.

\smallskip

\subsection{Operators with invariant cones. Monotone extremal norms}
\label{subsec:IC}

Let $K$ be a convex closed pointed nondegenerate cone with the
apex at the origin. In the sequel we write $K^*$ for the dual
cone: $K^* \, = \, \bigl\{u \in \re^d \ \bigl| \ \inf_{x \in K}(u,
x)\ge 0\, \bigr\}$. According to the Perron-Frobenius theorem, every operator~$B$
that leaves a cone~$K$ invariant has a positive leading eigenvalue $\lambda_{\max} = \rho(B)$
and $K$ contains a leading eigenvector corresponding to this eigenvalue.
Any leading eigenvector of $B$ that belongs to the cone~$K$ will be referred as Perron-Frobenius
eigenvector.

If all operators of the family~$\cM$ share a common
invariant cone, then Theorem~\ref{th.bar} on the existence of invariant  norms can
be slightly sharpened. First, the irreducibility condition can be
relaxed; second, an invariant  norm can  always be chosen to be
monotone with respect to the invariant cone.
 A function $g$ is {\em monotone} on a cone~$K$ if
$g(x)\ge g(y)$, whenever $(x -y)\,  \in \, K$. If $g$ is a
monotone norm defined on the cone~$K$, then it is extended
 onto~$\re^d$ in a standard way:
the unit ball of that norm is
\begin{equation}\label{ext}
\bigl\{x \in \re^d \ \bigr| \ \|x\| \le 1\, \bigr\}\quad = \quad
{\rm co}_s\,  \bigl\{\, x \in K \, , \, g(x)\, \le \, 1\,
\bigr\}\, .
\end{equation}
All extreme points of the ball defined by~(\ref{ext}) are in the
cones~$K$ and $-K$. Since the norm of any operator~$A$ is attained
at an extreme point of the unit ball, we see that if $A$ leaves
$K$ invariant, it attains its norm in the cone~$K$. Thus,
$$\|A\| = \max\limits_{x\in K, \ g(x)\le 1}g(Ax).$$
In particular, if $g$ is
an extremal norm for a family $\cM$, i.e., $\, \max_{i= 1, \ldots , m} \|A_i\|\, = \,
\widehat \rho $, then its
extension defined by~(\ref{ext}) is extremal as well. Thus, for
families with a common invariant cone it suffices to construct an
extremal monotone norm $g$ on that cone.

We are going to show that there exists not only extremal, but invariant
monotone norm on~$K$.
Recall, that a norm in $K$ is invariant for $\cM$ if
$$\max \{\|A_1 x\|, \ldots , \|A_mx\|\} = \widehat \rho \, \|x\|\, , \
x \in K.$$
To formulate the main result we need some further notation.
A hyperplane $L \subset \re^d$ is called {\em a plane of support} of a cone~$K$
if $L \cap K \, \ne \, \{0\}\, $ and $\, L \cap {\rm int}\, K \, = \, \emptyset$.
A {\em face} of a cone is its intersection with some plane of support.
For example, a spherical cone has only one-dimensional faces (rays);
the faces of the cone $K = \re^d_+$ are coordinate planes:
$F_{i_1 \ldots  i_r} \, = \, \{x \in \re^d_+ \ | \ x_{i_1} = \cdots x_{i_r} = 0\}\, , r = 1, \ldots , d-1$.
A face $F$ of a cone $K$ is {\em invariant} for an operator~$A$ if $AF \subset F$.
\begin{theorem}\label{th10}
If operators of a family $\cM = \{A_1, \ldots , A_m\}$ share an invariant cone $K$
and do not have common invariant faces of that cone, then $\cM$ possesses a monotone invariant norm on~$K$.
\end{theorem}
In the proof of Theorem \ref{th10} we use the following lemma from~\cite[section 4]{P3}:

\begin{lemma}\label{l300} \cite{P3}
For any cone $K$ and for any norm on this cone there is a homogeneous continuous function $\gamma (x)$ positive on the
interior of $K$ such that for every operator $B$ leaving the cone invariant we have
$\|Bx\| \ge \gamma (x) \|B\|$.
\end{lemma}

{\tt Proof.} The  case $\widehat \rho = 0$ is impossible, because in this case
the operators must have a common invariant face. This fact is simple, and we omit its proof.
If
$\widehat \rho > 0$, then after normalization it can be assumed that
$\widehat \rho = 1$. Let us first get an extremal norm for~$\cM$. Take any $\, e^* \in {\rm int}\, K^*$ and for each
$n \ge 1$ consider the function
\begin{equation}\label{norm}
 g_n(x)\quad = \quad \sup_{k \, \ge\,
n}\max_{A \, \in \, \cM^k}\ \bigl(\, e^*\, , \, Ax\, \bigr)\ , \qquad x \, \in \, K\, .
\end{equation}
Note the following properties of these functions.
\smallskip

\begin{itemize}
\item[1. ] For very $n$ we have $g_{n+1}(x) \le g_n(x)$, so the sequence $\{g_n\}_{n  \in \n}$ is monotone.
\smallskip

\item[2. ] We have
$$\sup\limits_{x \in K\,, (e^*, x) = 1} g_1(x) < \infty,$$
hence, by the monotonicity of the sequence $\{g_n\}$, all
$g_n$ are uniformly bounded on the unit sphere. To prove that $g_1$ is bounded
observe that the set $\, L \, = \, \bigl\{x \in K \ \bigl| \ g_1(x)\, < \, +\infty \, \bigr\}$
is either $\{0\}$, or $K$, or a common  invariant face of~$K$.
The latter contradicts the assumption. Assume $L = \{0\}$.
Then consider the compact set  $S = \{x \in K\ | \ (e^*, x)= 1\, \}$. For
each $j \ge 1$ we define the set
$$
V_j \quad = \quad \Bigl\{\ x \in S \ \Bigl| \
\ \max_{A \in \cM^j} \, (e^*, Ax) \ > \ 2\,  \
\Bigr\}.
$$
Thus $V_j$  consists of vectors for which some product~$A$ of length~$j$
increases the value $(e^*, x)$ more than twice. If $L = \{0\}$, then
 $\, \cup_{j \ge 1}V_j \, = \, S$,
and, since all $V_j$ are open in $S$, from the compactness of~$S$ it follows that
$\, \cup_{j = 1}^{N} V_j \, = \, S$ for some~$N$.
This means that
for every $x \in K$ there is a product~$A$ of length at most~$N$ increasing  the value $(e^*, x)$ at least twice.
Applying this argument successively~$k$ times, we conclude  that for every $x \in K$ there is a product~$\Pi_k$
of length~$l_k \le kN$ such that $(e^*, \Pi_k x)\, \ge \, 2^{k}\, (e^*, x)$,
and hence $\|\Pi_k\|\, \ge \, 2^{k}\, (e^*, x)\, \|e^*\|^{-1}$ (the norm of $\Pi_k$ is the Euclidean).
Taking the power $1/l_k$ and the limit as $k \to \infty$, we obtain $\, \widehat \rho\, \ge \, 2^{1/N}$,
which is a contradiction. Thus, $L \ne \{0\}$, and hence $L = K$. Thus, each function $g_n$
is bounded on~$S$.

\smallskip

\item[3. ] For each $n$ the function $g_n$ is homogeneous, positive (as a supremum of positive values)
and convex on~$K$, as a pointwise supremum of linear functionals. Thus, $g_n(\cdot )$ is a norm on~$K$.
\smallskip

\item[4. ] For every $x \in K$ we have    $g_n(A_ix)\le g_n(x)\, , \ A_i\in
\cM$. Hence for each $n$ the norm $g_n$ is extremal. Thus, we have established the
existence of a monotone sequence of extremal norms.
\smallskip

\item[5. ] For every operator $A$ leaving $K$ invariant we denote by
$$\|A\|_{e^*} = \sup_{x \in S}\, (e^*, Ax)$$
the operator norm corresponding to the norm $\|x\|_{e^*}  = (e^* , x)$.
Since $\widehat \rho = 1$, it follows that $$\max\limits_{A \in \cM^k} \, \|A\|_{e^*} \, \ge \, 1.$$
Now using Lemma~\ref{l300}, we obtain $\max\limits_{A \in \cM^k} \, \|Ax\|_{e^*} \, \ge \, \gamma (x)\|x\|_{e^*}$.
This holds for each $k$, therefore, $g_n(x) \ge \gamma (x)\|x\|_{e^*}$ for every $x \in K$.

\item[6. ] Since the sequence $\{g_n(x)\}_{n \in N}$ is non-increasing and bounded below, it converges
to some limit function $g(x)$. For every $x \in {\rm int}\, K$ we have $g(x) \ge \gamma (x)\|x\|_{e^*} > 0$.
Thus, the function $g$ is convex, positively homogeneous, and invariant, i.e., possesses the property
$$
g(x) \ = \ \max_{j=1, \ldots , m} \ g(A_jx)\, .
$$
\end{itemize}
It remains to show that $g$ is positive on $K$, in such case it constitutes an invariant norm.
Since $g(x) > 0$ for $ x \in {\rm int}\, K$, we see that the set $L = \{x \in K \ | \ g(x) = 0\}$ lies
on the boundary of $K$. This set is obviously convex, hence it is contained on a face of~$K$.
Let $V$ be the minimal (by inclusion) face containing $L$.
Since $A_j L \subset L\, , \ j = 1, \ldots m$, it follows that $A_j V \subset V\, , \ j = 1, \ldots m$,
which contradicts the assumption. Thus, $g$ is an invariant norm, which completes the proof.

 {\hfill $\Box$}
\smallskip

Now we focus on the case of nonnegative operators, i.e., operators
defined by nonnegative matrices (which means, with nonnegative entries).
Each family of nonnegative operators share an invariant cone~$\re^d_+$.
A family of nonnegative operators is called {\em positively-irreducible}
if they do not have common invariant faces among the coordinate planes.
Applying Theorem~\ref{th10} to the case $K = \re^d_+$, we obtain:
\begin{cor}\label{c10}
A family of nonnegative positively-irreducible operators possesses a monotone invariant norm on~$\re^d_+$.
\end{cor}
\begin{remark}\label{r10}
{\rm The assumption of Corollary~\ref{c10}
is not restrictive, because the general case of nonnegative matrices
is reduced to the case of matrices without invariant coordinate planes.
If the matrices possess common invariant planes, then after a suitable permutation
of the basis vectors all the matrices get a block upper-triangular form.
The joint spectral radius of the matrices equals to the largest joint spectral radius
of the blocks. Thus, the problem of JSR computation comes to several similar problems with nonnegative
matrices of smaller dimensions.
A fast polynomial procedure to realize this reduction can be found in~\cite[section 2]{JPB1}.
}
\end{remark}

\medskip

Before describing Algorithm~\textbf{(P)} we formulate an analogue of Lemma~\ref{l.5} for
nonnegative operators.

\begin{lemma}\label{l.5p}
Let a nonnegative  operator $B$ have a unique simple Perron-Frobenius eigenvalue
$\lambda$  with an eigenvector $v$; let also $v^*$
be the Perron-Frobenius eigenvector of $B^*$ such that $(v^*, v) = 1 $. If
for some nonnegative operator $C$ one has~$ (v^*, Cv)\, > \, 1$,
then for sufficiently large $r$ the  operator $B^rC$ has a unique
simple Perron-Frobenius eigenvalue bigger
than~$\lambda$.
\end{lemma}
The proof is literally the same as for Lemma~\ref{l.5}.
{\hfill $\Box$}
\smallskip

If a cone $K\subset \re^d$ is fixed, then for a given set $Q \subset \re^d$
we denote
$$
 {\rm co}_{-}(Q)\ = \ \bigl( \, {\rm co}(Q)\, - \, K\, \bigr)\, \cap \, K\ = \
 \bigl\{ \, x \, \in \, K\ \bigl| \ x = y \, - \, z\, , \ y \, \in \,
{\rm co}\, (Q)\, , \   z \in K\bigr\}\, .
$$
If the cone~$K$ is not specified, we always assume  $K = \re^d_+$.

Everywhere below in this section the family~$\cM$ is assumed to be positively irreducible
(see Remark~\ref{r10}).
\smallskip

\subsection{{\em Algorithm~(P)} versus {\em Algorithm~(R)}}
\label{subsec:PvR}

Algorithm~(P) ia similar to Algorithm~(R) but has some peculiar differences which
we remark in the sequel. The second has a major computational importance.

\begin{itemize}
\item[(i) ] By the Perron-Frobenius theorem the candidate s.m.p. $\Pi$ has a nonnegative
leading eigenvalue~$\lambda_{\max} = \rho(\Pi)$, and the corresponding eigenvector belongs
to~$\re^d_+$. We assume that ${\lambda_{\max} > 0}$. Hence the main assumption for
Algorithm~(R) holds true automatically.
\item[(ii) ] The LP problem performed in the loop at step $k$ should be replaced by the following:
\begin{equation}\label{LP.p}
\left\{
\begin{array}{rcl}
\max & & t_0
\\[0.25cm]
{\rm subject \ to}
     & & t_0\,  z \le \sum\limits_{x \in \cV_{k}}\, t_x\, x
\\[0.3cm]
{\rm and}
     & & \sum\limits_{x \in \cV_{k}}\, t_x \le 1, \qquad t_x \ge 0 \quad \forall x \in \cV_k
\end{array}
\right.
\end{equation}
where we recall that $z = \widetilde A v$.

\item[(iii) ] In contrast to the
case (R), we work now only with nonnegative vectors, and do not include the vectors $-v_i$ to the polytope~$P_k$.
We do not construct a symmetric polytope ${\rm co}_s\, (\cV_k)$, but a positive polytope ${\rm co}_{-}\, (\cV_k)$.
So, our norm will have a unit ball ${\rm co}_{-}\, (\cV_k)$. This explains the differences between LP-problems
(\ref{LP.r}) and (\ref{LP.p}).

\item[(iv) ]

Condition (\ref{cond.r}) in the {\tt Stopping criterion} should be replaced by the following:
\begin{equation}\label{cond.p}
 \bigl(v_j^*\, , \, \widetilde A\, v \bigr)\, \quad \le \quad 1\ , \qquad
j \, = \, 1, \ldots ,n\, .
\end{equation}
\end{itemize}

\smallskip

\subsection{Explanations and proofs}
\label{subsec:EPP}

The theoretical base of Algorithm~\textbf{(P)} is actually the same as
for Algorithm~\textbf{(R)}. Let us only stress the distinctions. First of
all, we construct a monotone extremal norm on the positive
orthant~$\re^d_+$, so the polytopes $P_k$ are in the orthant, and
are not centrally-symmetric. That is why we do not need to
symmetrize the convex hull of $\cV_k$. Consequently,
LP-problem~(\ref{LP.p}), in contrast to LP-problem~(\ref{LP.r}),
does not have extra variables $q_x\, , x \in \cV_k$. The other
difference is that the equality constraint $t_{0} \widetilde A\, v \,
= \, \sum_{x \in \cV_{k}}\, t_x\, x\, $ becomes an inequality.
This means that the polytope $P_k$ is not a convex hull of
$\cV_k$, but $\, {\rm co}_{-}(\cV_k)$. Thus  there are two
advantages of Algorithm~\textbf{(P)}: 1) the number of variables
and the number of constraints in the LP-problem is  a half of the
LP-problem in Algorithm~\textbf{(R)}; 2) the polytope~${\rm
co}_{-}(\cV_k)$ is larger than~${\rm co}\, (\cV_k)$, therefore
this algorithm sorts out more vertices (``dead branches'') at each
step, which leads to a lower complexity. In practice
Algorithm~\textbf{(P)} works much faster than
Algorithm~\textbf{(R)} (see Section~\ref{sec:IE} and Section~\ref{sec:numer}).


In the worst case, if the algorithm has not terminated, one gets
an approximate value of JSR from inequality
\begin{equation}\label{final.p}
[\rho (\Pi)]^{1/n}\quad \le \quad \widehat \rho (\cM) \quad \le \quad
(t_N)^{-1}[\rho (\Pi)]^{1/n}\, .
\end{equation}
\begin{prop}\label{p2.p}
For an arbitrary  positively-irreducible family~$\cM$
we have $\, t_N \to 1\, $ as~$l \to
\infty$ and $N \to \infty$ in estimate~(\ref{final.p}).
\end{prop}
The proof is the same as for  Proposition~\ref{p2.r} with the use
of Theorem~\ref{th10} instead of Theorem~\ref{th.bar}.
The proof  of the efficiency of the stopping criterion (Proposition \ref{p1.p})
is also the same as for Proposition~\ref{p1.r}:
\begin{prop}\label{p1.p}
Assume  the
Perron-Frobenius eigenvalue of~$\Pi$ is unique and simple.  If the
assumption of the algorithm is wrong (i.e., $\Pi$ is not an
s.m.p.)  then for every $j = 1, \ldots , n$
condition~(\ref{cond.p}) is violated at some step. Conversely, if
condition~(\ref{cond.p}) is violated at some step for some~$j$,
then $\Pi$ is not an s.m.p.
\end{prop}

\begin{remark}\label{r312}
{\rm In Algorithm~\textbf{(P)} the family~$\cM$ is assumed to be positively irreducible.
Actually, this was done for the sake of simplicity. The algorithm can be applied to arbitrary nonnegative families.
If the algorithm terminates after $k$th iteration, and the set $\cV_k$ does not lie in
a coordinate plane of a smaller dimension, then $P_{k-1}$ is an extremal polytope, and $\Pi$ is an s.m.p.
That condition means that for each $i = 1, \ldots , d$ there is a vector from~$\cV_k$ with strictly positive
$i$th coordinate. This simple condition allows us to apply
Algorithm~\textbf{(R)} to arbitrary family, without preliminary checking its positive irreducibility. Nevertheless,
if the family $\cM$ is reducible, then it is always advisable to factorize it before starting the
algorithm,
because this significantly reduces the dimension (see Remark~\ref{r10}). Especially as the factorization is realized by
a fast polynomial routine~\cite{JPB1}.
}
\end{remark}


\section{ The criterion for finite termination of Algorithms~(R), (C) and (P).}
\label{sec:CFT}

Algorithms~\textbf{(R)},~\textbf{(C)}, and~\textbf{(P)}  compute JSR by step-by-step
constructing a polytope norm. If the algorithm terminates within finite time,
then it produces an extremal polytope, and, hence, proves that the chosen product~$\Pi$ is an s.m.p.
If it does not terminate, then it gives upper and lower bounds for JSR
that converge to the exact value. An important issue is the following: what are the conditions
for the family~$\cM$ and for the product~$\Pi$, under which the algorithm terminates
and produces the extremal polytope~? Conditions guaranteeing the convergence of the algorithm to
an extremal polytope norm have been discussed in \cite{GWZ,GZ2}. We give here a further
result which is related to those obtained in the mentioned papers.

Certainly, the product $\Pi$ must be spectral maximizing for that.
This condition, however, does not guarantee the convergence of the algorithm. It appears that a bit stronger condition
solves the problem completely: it is both sufficient and necessary. The product $\Pi$ has to be
not just maximizing but {\em dominant}. To formulate the criterion we need some further notation.

Let $\cM = \{A_1, \ldots , A_m\}$ be a given family of operators,
$\Pi = A_{d_n}\cdots A_{d_1}$ be some product, which is not a
power of a shorter product, $n \ge 1$. We denote $\, \widetilde A_i \,
= \, [\rho(\Pi)]^{-1/n}A_i\, , \ \widetilde M = \{\widetilde A_1, \ldots ,
\widetilde A_m\}\, $ and $\, \widetilde \Pi \, = \, \widetilde A_{d_n}\cdots
\widetilde A_{d_1}$. Clearly, the spectral radius of any power of
$\widetilde \Pi$ or of any power of its cyclic permutation is~$1$.
\begin{defi}\label{d.simple}
A product $\Pi\in \cM^n$ is called dominant for the family
$\cM$  if there is $q< 1$ such that the spectral radius of every
product of operators of the normalized family $\widetilde M$, that is
not a power of~$\widetilde \Pi$ nor a power of its cyclic
permutations, is smaller than $\, q$.
\end{defi}
Obviously, the dominant product along with all its cyclic permutations are all s.m.p.,
but vice versa.

The following theorem  gives a sharp criterion on the family~$\cM$ ensuring that
our algorithm produces an extremal polytope.
\begin{theorem}\label{th.cond-r}
For each of Algorithms \textbf{(R)}, \textbf{(C)} and \textbf{(P)}
the following holds:

the algorithm terminates within finitely many iterations if and
only if $\Pi$ is  dominant for $\cM$, and its leading eigenvalue
is unique and simple.
\end{theorem}
The proof is in Appendix.

\begin{cor}\label{c200}
If a family $\cM$ possesses a dominant product, whose leading
eigenvalue is unique and simple, then it has an extremal polytope.
\end{cor}

\begin{remark}\label{r20}
{\rm Theorem~\ref{th.cond-r} remains true even if we do not apply
the stopping criterion in the algorithms.}
\end{remark}

\begin{remark}\label{r23}
{\rm Theorem~\ref{th.cond-r} implies  that if the algorithm
terminates within finite time, then the family~$\cM$ possesses a
dominant product. Actually the algorithm ensures that a chosen product~$\Pi$
is dominant. In numerical examples from applications
(Sections~\ref{sec:Appl}) and from randomly generated matrices (Section~\ref{sec:numer})
 most of matrix families possess dominant products.}
\end{remark}


The assumption on a dominant product allows to exclude a limit spectrum maximizing product,
that is a matrix in the closure of the multiplicative semigroup of $\widetilde {\mathcal M}$
with spectral radus equal to $1$ (see \cite{GZ4}).

Consider in fact the following example.
Let ${\cM}=\{ A_1,A_2 \}$:
\begin{eqnarray*}
\displaystyle{ A_1 \, = \, \left(
\begin{array}{rr}
1 & 1 \\
0 & 1
\end{array}
\right) \qquad \mbox{and} \qquad A_2 \, = \, \frac{4}{5}\,  \left(
\begin{array}{rr}
1 & 0 \\
1 & 1
\end{array} \right). } &&
\end{eqnarray*}
We can prove that $\widehat\rho(\cM) = 1+\frac{1}{\sqrt{5}}$ and there is a unique {\em finite} spectrum maximizing product $\Pi=A_1 A_2$
such that $\rho( \widetilde P) = 1$ (apart from its cyclic permutation $A_2 A_1$ and their powers). Nevertheless the product is
not dominant.
In fact the sequence of products $\widetilde Q_k = \widetilde A_1 \left( \widetilde A_1 \widetilde A_2 \right)^k$ is convergent and such that
\[
\lim\limits_{k \rightarrow \infty} \widetilde Q_k =
\left(
\begin{array}{cc}
 \frac{\sqrt{5}+1}{4} & \frac{1}{2} \\[0.2cm]
 \frac{\sqrt{5}-1}{4} & \frac{3-\sqrt{5}}{4}
\end{array}
\right) := \widetilde Q_\infty
\]
which is such that $\rho(\widetilde Q_\infty) = 1$.

This implies that $P$ is not dominant; the matrix $\widetilde Q_\infty$ is indeed a limit spectrum maximizing product
of the normalized family $\widetilde \cM$.


Indeeed the algorithms \textbf{(R)} and \textbf{(P)}  do not converge
when applied to this example. However, if we modify the algorithms and
remove a vector when it lies on the boundary of the polytope $P_{k-1}$,
then we obtain a finite convergence also in this case.


\section{ Computing the lower spectral radius. Algorithm~(L).}
\label{sec:LSR}

In this section we describe a method for the exact computation of the lower spectral radius
of a finite family of matrices.

\subsection{Antinorms on convex cones}
\label{subsec:antinorm}

To extend our approach to computing the lower spectral radius, first of all we need
the notion of extremal norm for this case. One can define it
by the inequality $\, \min_{A_i \in \cM}\|A_ix\|\, \ge \, \check \rho \, \|x\|\, , \ x \in \re^d$.
However, simple examples show that such a norm may not exist even for very ``good''
families~$\cM$ (for instance, irreducible families of positive matrices). The reason is that
the function $\, x \, \mapsto \, \min_{A_i \in \cM}\|A_ix\|$ may not be convex,
in which case it is not a norm (in contrast to the situation with JSR, when
the function $\, x \, \mapsto \, \max_{A_i \in \cM}\|A_ix\|$ is always a norm).
One of the ways to generalize the notion of extremal norm for the lower spectral radius
is to consider concave positive homogeneous functionals on $\rd$ instead of convex ones (i.e., instead of norms).
However, such functionals do not exist.
Indeed, if $f$ is concave, then $f(x) + f(-x) \, \le \, f(0)\, = \, 0$
(homogeneity failure), hence $f$ cannot be positive. Nevertheless, extremal concave ``norms''
can be defined, provided all operators of the family~$\cM$ share an invariant cone.
In particular, this can be done for families of nonnegative matrices.
As in the previous section,  $K$ is a convex closed pointed nondegenerate cone with an
apex at the origin.
\begin{defi}\label{d.antinorm}
An antinorm is a continuous nonnegative nontrivial (not identical
zero) concave positively-homogeneous function on a cone~$K$.
\end{defi}
From the concavity it easily follows that an antinorm can vanish
only on the boundary of $K$. An antinorm $f$ is called {\em positive}
if $f(x) > 0$ for all $x \in K\setminus \{0\}$. So, $f$ is positive, whenever
it is positive on the boundary. This is well known that a concave
function is continuous at each interior point of its domain. Hence
the continuity condition for antinorms can be relaxed to the
continuity on the boundary.

\begin{figure}[ht]
\centering \global\def\path{#1}\input{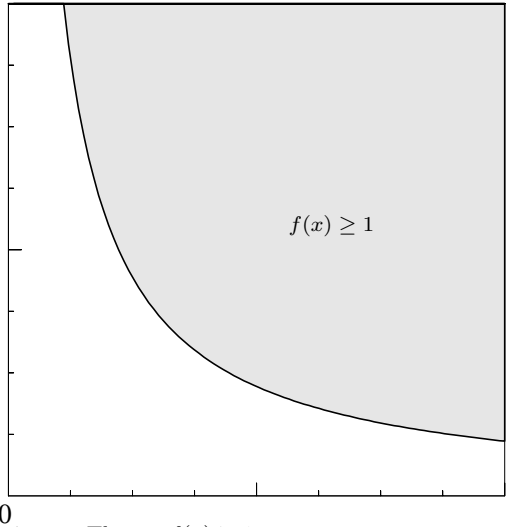}
 \caption{Example of antinorm. The set $f(x) \ge 1$. \label{fig:anorm}}
\end{figure}

Let us observe some basic properties of antinorms. First of all,
every antinorm is asymptotically bounded above by every norm in~$\re^d$:
\begin{lemma}\label{l10}
For any antinorm $f$ and for any norm $\|\cdot \|$ there is a constant $C$ such that
\linebreak $f(x)\, \le \, C\, \|x\|\, , \ x \in K$.
\end{lemma}
{\tt Proof.} Since $f$ is continuous, the value $C\, = \,
\sup\limits_{x \in K, \|x\|=1}f(x)\, $
 is finite.  Now by the
homogeneity the lemma follows.

{\hfill $\Box$}
\smallskip

Consider  now a family of operators $\cM = \{A_1, \ldots ,
A_m\}$ that share an invariant cone $K$.
\begin{prop}\label{p20}
If for some antinorm $f$ and for a constant $\lambda$ we have
$f(A_ix)\, \ge \, \lambda\, f(x)\, , \ x \in K\, , \ A_i \in \cM$,
then $\check \rho \, \ge \, \lambda$.
\end{prop}
{\tt Proof.} Applying Proposition~\ref{p20} for an arbitrary point
$e \in {\rm int}\, K\, , \|e\|=1$, we get $\, \|A_{d_k}\cdots
A_{d_1}e\| \, \ge \, C^{-1}f(A_{d_k}\cdots A_{d_1}e)\, \ge \,
C^{-1}\, \lambda^k\, f(e)$. Thus,
$$\min\limits_{d_k, \ldots ,d_1}\|A_{d_k}\cdots A_{d_1}\| \, \ge \,
C^{-1}\, f(e)\,\lambda^k.
$$
Taking the power $1/k$ and the limit as $k \to
\infty$, we conclude the proof.

{\hfill $\Box$}
\smallskip

\begin{defi}\label{d.ext-antinorm}
An antinorm is called extremal if $\ f(A_ix)\, \ge \, \check
\rho\, f(x)\, , \ x \in K\, , \ A_i \in \cM$.
\end{defi}
Similar to monotone norms, an antinorm is called {\em monotone} if
$f(x)\ge f(y)$, whenever $\, x - y\, \in\, K$.
\begin{theorem}\label{th20}
For every family of matrices with a common invariant cone $K$ there
exists  a monotone  extremal antinorm on~$K$.
\end{theorem}
{\tt Proof.} If $\check \rho = 0$, then any antinorm suffices. If
$\check \rho > 0$, then after normalization it can be assumed that
$\check \rho = 1$. Take any $\, e^* \in {\rm int}\, K^*$ and
consider the function
\begin{equation}\label{anorm}
 f(x)\quad = \quad \inf_{k \, \ge\,
0}\min_{A \, \in \, \cM^k}\, \bigl(\, e^*\, , \, Ax\, \bigr) \ , \qquad x \, \in \, K\, .
\end{equation}
 This function is concave,
homogeneous  and monotone, and $f(A_ix)\ge f(x)$ for each $A_i\in
\cM$. It remains to show that $f$ is not an identical zero.
Consider the compact set $S = \{x \in K\ | \ (e^*, x)= 1\, \}$ and for
each $k \ge 1$ define the set
$$
U_k \quad = \quad \Bigl\{\ x \in S \ \Bigl| \
\ \min_{A \in \cM^k} \, (e^*, Ax) \ < \ \frac12\,  \
\Bigr\}.
$$
If $f \equiv 0$, then $\, \bigcup_{k \ge 1}U_k \, = \, S$,
and, since all $U_k$ are open in $S$, from the compactness of~$S$ it follows that
$\, \bigcup_{k = 1}^N U_k \, = \, S$ for some~$N$.
This means that
 for every $x \in K$ there is a product~$A$ of length at most~$N$ reducing  the value $(e^*, x)$ at least twice.
 Applying this argument $k$ times, we obtain that for every $x \in K$ there is a product~$\Pi_k$
 of length~$l_k\le kN$ such that $(e^*, \Pi_k x)\, \le \, 2^{-k}\, (e^*, x)$.
 Note that if $x \in {\rm int}\, K$, then there is a constant $C$ that depends on $e^*$ and on $x$, and such that
  $\|B\|\, \le \, C\, (e^*, B x)$ for any operator~$B$ that leaves~$K$ invariant~(see, for instance~\cite{P3}).
  Thus, $\|\Pi_k\|\, \le C \, 2^{-k}\, (e^*, x)$ for each~$k$. Taking the power $1/l_k$
  and the limit as $k \to \infty$, we get $\, \check \rho\, \le \, 2^{-1/N}$,
  which is a contradiction.

{\hfill $\Box$}
\smallskip

Applying this theorem for the case $K = \re^d_+$, we obtain:
\begin{cor}\label{c20}
For an arbitrary family of nonnegative matrices there is a  monotone extremal antinorm on~$\re^d_+$.
\end{cor}

In the algorithm we need the following analogue of Lemma~\ref{l.5p}, whose proof is the same.
\begin{lemma}\label{l.5min}
Under the assumptions of Lemma~\ref{l.5p}, if for some nonnegative operator $C$ one has~$(v^*, Cv)\, < \, 1$
then for sufficiently large $n$ the operator $B^nC$ has
a unique simple Perron-Frobenius eigenvalue smaller than~$\lambda$.
\end{lemma}

\medskip

Now we are ready to describe Algorithm~\textbf{(L)} of LSR
computation for nonnegative matrices. We use the notation $\, {\rm
co}_{+}(X)\, = \, {\rm co}(X) \, + \, K \, = \, \bigl\{ \, x + h\ \bigl| \ x \in {\rm co}\,
(X)\,  \ h \, \in \, K\, \bigr\}$, where $X$ is a subset
of~$\re^d$ and $K \subset \re^d$ is a cone. If $X$ if finite, then
${\rm co}_{+}(X)$ will be referred as an {\em infinite polytope}.
Thus, an  infinite polytope is a set $P+K$, were $K$ is a cone and $P$
is a polytope. In the algorithm we always assume~$K = \re^d_+$.

A product $\Pi \in \cM^n$ is called {\em the spectral lowest product (s.l.p.)} if $[\rho(\Pi)]^{1/n} = \, \check \rho
(\cM)$. We shall also call it a {\em spectrum  minimizing product}, but always use the abbreviation s.l.p. to avoid
confusion with the spectral maximizing product (s.m.p.).
  By inequality~(\ref{lsr.est}) a product $\Pi$ is an s.l.p. iff $[\rho(\Pi)]^{1/n} \le \, \check \rho
(\cM)$. Proposition~\ref{p20} implies that if there is an  antinorm $f: \re^d_+\, \to \, \re_+\, $
such that $\, f(A_jx) \, \ge \, [\rho(\Pi)]^{1/n} f(x)$, \,
$x \in \re^d_+$, then $[\rho(\Pi)]^{1/n} = \, \check \rho (\cM)$, and $f$ is extremal.
The main idea of the algorithm is to select a candidate $\Pi$ for s.l.p. (by a reasonable exhaustion)
and then to prove that it is actually an s.l.p. The proof is by step-by-step constructing an extremal infinite polytope,
which  generates an extremal antinorm.

\smallskip

\subsection{Algorithm~(L)}
\label{subsec:AL}

{\tt Initialization.} We have an arbitrary family $\cM = \{A_1, \ldots ,
A_m\}$. For some (as large as possible)~$l$ we look over all
products $\Pi$ of length $\le l$ and take one with the smallest
value~$[\rho(\Pi)]^{1/n}$, where $n$ is the length of the product.
We take the shortest product possessing this property and denote
it as~$\Pi= A_{d_n}\cdots A_{d_1}$. If $\rho(\Pi)=0$, then $\check
\rho (\cM)= 0$, and the algorithm terminates. So, we assume~$\rho(\Pi)
> 0$. By the Perron-Frobenius theorem $\Pi$ has a positive leading
eigenvalue, and the corresponding eigenvector belongs
to~$\re^d_+$.
 We normalize the family $\cM$ as
 $\widetilde A_i = [\rho(\Pi_1)]^{\, -1/n}\, A_i\, , \, \widetilde M = \{\widetilde A_i\}_{i=1}^m$.
The leading eigenvalue of the operator $\widetilde \Pi = \widetilde
A_{d_n}\cdots \widetilde A_{d_1}$ equals to~$1$.

Let $\widetilde \Pi_1 = \widetilde \Pi\, , \
\widetilde \Pi_i = \widetilde A_{d_{i-1}} \cdots \widetilde A_{d_1}\widetilde A_{d_n}\cdots
\widetilde A_{d_i}$ be a cyclic permutation of $\widetilde \Pi_1$,
$\, i = 2, \ldots , n$. We take a Perron-Frobenius eigenvector~$v_1$ of~$\widetilde \Pi_1$
(if it is not unique, take any of them), and
$$
v_i \quad = \quad \widetilde A_{d_{i-1}}\cdots \widetilde A_{d_{1}}v_1 \, .
$$
Thus, $v_i$ is a Perron-Frobenius eigenvector of~$\widetilde \Pi_i$
with the eigenvalue~$1$.

In case the leading eigenvalue $\lambda_{\, \max}\, = \, 1$ is unique and simple,
we also need a dual system of vectors: $v_1^*$ the leading eigenvector of~$\widetilde \Pi_1^*$
normalized as $(v_1^*, v_1) = 1$ (Remark~\ref{r.5}), and
$$
v_i^* \quad = \quad \widetilde A_{d_{1}}^*\cdots \widetilde A_{d_{i-1}}^*v_1^*\ , \qquad  i \ = \ 2, \ldots, n\, .
$$
Thus, $v_i$ and $v_i^*$ are the Perron-Frobenius eigenvectors of $\Pi_i$ and $\Pi_i^*$ respectively,
and $(v_i^*, v_i) = 1$.

\smallskip

{\tt Set $k=0$}. We set $\cV_0\, = \, \cU_{\, 0} \, = \, \{v_1,  \ldots , v_n \}$
 and $\, \cR_0\, = \, \bigl\{ \, ( v_i\, , \, \widetilde A_p) \bigl|\
i = 1, \ldots , n\, ; \ p=1, \ldots , m \, , \  p\ne d_i \bigr\}$.
\medskip

{\tt Main loop}
\smallskip

{\tt For $k \, \ge \, 1$}. We have  finite sets $\cV_{k-1}\, \subset \, \re^d\, ,\,  \cU_{k-1} \,
\subset \, \cV_{k-1}$, and  $\cR_{k-1}\, \subset \, \cU_{k-1}\times \widetilde \cM$.
Put $\cV_{k}\, = \, \cV_{k-1}\, $ and $\, \cU_k \, = \, \emptyset$.

We successively take all pairs $(v, \widetilde A)\in \cR_{k-1}$.
If for a given pair $\widetilde A \, v\, = \, 0$,
then we stop the algorithm, it is inapplicable for this case.
If $z = \widetilde A \, v\, \ne  \, 0$, then we  solve the following LP problem with variables  $t_0$ and
$\{t_x\}_{x \in \cV_k}$:
\begin{equation}\label{LP.pmin}
\left\{
\begin{array}{rcl}
\min & & t_0
\\[0.25cm]
{\rm subject \ to}
     & & t_0\,  z \ge \sum\limits_{x \in \cV_{k}}\, t_x\, x
\\[0.3cm]
{\rm and}
     & & \sum\limits_{x \in \cV_{k}}\, t_x \ge 1, \qquad t_x \ge 0 \quad \forall x \in \cV_k.
\end{array}
\right.
\end{equation}
The value of the problem, i.e., $\, \min \, t_{0}\, $ will be denoted by $t_{\{v, \widetilde A\}}$.
Thus, for every pair $(v, \widetilde A) \, \in \ \cR_{k-1} \, $ we have a nonnegative
number $\, t_{\{v, \widetilde A\}}$,   which may take value $+\infty$, when
the system of inequality constraints has no solution.
\smallskip

{\tt If} $\, t_{\{v, \widetilde A\}} \, < \, 1$,
then leave the sets $\cV_k$ and $ \cU_k$ as they are, take the next pair~$(v, \widetilde A)\, \in \, \cR_{k-1}$
and consider problem~(\ref{LP.pmin}) for it.
\smallskip

{\tt Otherwise If} $\, 1 \, \le \, t_{\{v, \widetilde A\}}$, then
\smallskip

\ {\tt If} the leading eigenvalue of $\Pi$
 is  unique and simple, we apply the following stopping criterion:
\smallskip

\ \ {\tt Stopping criterion}
\smallskip

We check the condition
\begin{equation}\label{cond.pmin}
 \bigl(v_j^*\, , \, \widetilde A\, v \bigr)\, \quad \ge \quad 1\ , \qquad
j \, = \, 1, \ldots ,n\, .
\end{equation}

\ \ \ {\tt If} (\ref{cond.pmin}) is satisfied, then we set $\,
\cV_k \, = \,  \cV_k\cup \{ \widetilde A\, v\}\, , \cU_k \, = \,
\cU_k\cup \{ \widetilde A\, v\}$, take the next pair~$(v, \widetilde A)\,
\in \, \cR_{k-1}$ and consider problem~(\ref{LP.pmin}) for it.

\ \ {\tt Otherwise If} (\ref{cond.pmin})  is not satisfied, then   our
assumption is wrong, $\Pi$ is not an s.l.p., and $\check \rho
(\widetilde M)< 1$ (Lemma~\ref{l.5min}). We stop the algorithm and go
either to the {\tt Final step}, or back to the {\tt
Initialization.} In the latter case we need to find another
pretender to s.l.p. The first opportunity is to increase~$l$ and
to look over all products of a bigger length. Lemma~\ref{l.5min}
provides also a different approach.  We take an index $j$, for
which $ \, (v_j^*\, , \, \widetilde A\, v )\, \, < \,  1$. Applying
Lemma~\ref{l.5min}, we conclude that there is~$r$ such that
$\lambda_{\max}(\widetilde \Pi_j^r\, \widetilde A_{s_q}\cdots \widetilde
A_{s_1})\, < \, 1$, where $\widetilde A_{s_q}\cdots \widetilde A_{s_1} v_j
\, = \, \widetilde A v$. We take the new initial product $\Pi =
\Pi_j^r\,  A_{s_q}\cdots  A_{s_1} $ and restart the
algorithm.
\smallskip

\ \ \ {\tt End If}
\medskip

\ {\tt Otherwise If} the leading eigenvalue of $\Pi$ is not unique or multiple, then we do not apply the stopping criterion, and
set $\,  \cV_k \, = \,  \cV_k\cup \{ \, \widetilde A\, v\}\, ,
\  \cU_k \, = \, \cU_k\cup \{ \widetilde A\, v\}$, take the next pair~$(v, \widetilde A)\, \in \, \cR_{k-1}$
and consider problem~(\ref{LP.pmin}) for it.
\smallskip

\ {\tt End If}
\medskip

{\em The $k$th step is over, when all pairs $(v, \widetilde A)\, \in
\, \cR_{k-1}$ are exhausted.}
\smallskip

{\tt If} $\, \cU_k = \emptyset$,
then $\check \rho (\widetilde M) = 1$, and so $\check \rho (\cM)\, =
\, [\rho(\Pi)]^{1/n}$. The extremal infinite polytope is $P_{k-1}
= {\rm co}_{+}\, (\cV_k)\, $, and the product $ \Pi$ is an s.l.p. for $ \cM$.
The algorithm terminates after the $k$th step.

\smallskip

{\tt Otherwise If} $\cU_k \ne \emptyset$, then we set  $\cR_k \, = \, \cU_k \times \widetilde \cM$
 and go to the $(k+1)$st step.

\smallskip
{\tt End If}

\smallskip
{\tt End For}
\medskip

\noindent {\tt Final step.} If the algorithm has not terminated,
then we stop it after some $N$ steps, denote $t_N\, = \, \,
\max\limits_{(v, \widetilde A)\in \cR_{N-1}}\, t_{\{v, \widetilde A\}}$,
and have the following estimate for the lower spectral radius:
\begin{equation}\label{final.pmin}
(t_N)^{-1} [\rho (\Pi)]^{1/n}\quad \le \quad \check \rho (\cM)
\quad \le \quad  [\rho (\Pi)]^{1/n}\,  .
\end{equation}
\medskip

{\tt End of Algorithm \textbf{(L)}}
\smallskip

\subsection{Explanations and proofs}
\label{subsec:EPL}

The algorithm produces a sequence of embedded infinite polytopes
$P_1 \subset P_2  \subset \ldots $ such that $P_{j+1} \, = \, {\rm
co}_{+}\, \bigl\{ \widetilde A_1P_j, \ldots , \widetilde A_mP_j\bigr\}\, $
and $P_{j} \subset P_{j+1} $ for every~$j$. If the algorithm
terminates after the $k$th step, then $P_{k} = P_{k-1}$. The $k$th
step is actually needed only to ensure that the infinite polytope
$P_{k-1}$ is extremal, i.e., $\widetilde A_j\, P_{k-1} \subset
P_{k-1}\, , j=1, \ldots , m$. Since $0\notin P_{k-1}$ (otherwise
at some step we had $\widetilde A v = 0$, in which case the algorithm
would be stopped by the end of the~$k$th step), the antinorm
$f_{k-1}(x)\, = \, \sup\, \bigl\{ t^{-1} \ \bigl| \ t\, x \, \in
\, P_{k-1} \bigr\}\, $ is well-defined on~$\re^d$. Since
$f_{k-1}(\widetilde A_j x)\, \ge \, f_{k-1}(x)\, $ for every $x$, we see
that $\check \rho (\widetilde \cM) \, \ge \, 1$. On the other hand,
$\check \rho (\widetilde \cM)\, \le \, [\rho (\widetilde \Pi)]^{1/n}\, =
\, 1$. Thus $\check \rho (\widetilde \cM)\, = \, 1$ and so $\,
\check \rho (\cM)\, = \, [\rho (\Pi)]^{1/n}$.
 Thus, {\rm if the algorithm terminates within finite time,
then the s.l.p. and the exact value of LSR are found.} In this
case $P_{k-1}$ is an extremal infinite polytope and $f_{k-1}$ is an
extremal antinorm.

Although the extremal antinorm is obtained numerically, the
results are actually exact, because the algorithm  removes only
those points~$v$, for which the strict inequality $f_i(v)> 1$
holds, where $f_i$ is the antinorm generated by the current
polytope~$P_i\, = \, {\rm co}_{+}\, (\cV_i)$.
\smallskip

If the algorithm does not terminate within finitely many steps,
then we have estimate~(\ref{final.pmin}) to get an approximate
value of the LSR. The right hand side inequality is obvious, the
left hand side is equal to $\min\limits_{j=1, \ldots ,
m}\inf\limits_{f_{N-1}(x)\, = \, 1}f_{N-1}(A_jx)\, \le \, \check
\rho (\cM)$, from which the estimate follows.
\smallskip

\subsection{Efficiency results for Algorithm~\textbf{(L)}}
\label{subsec:ERL}

Let us start with the stopping criterion. If the stopping
criterion is applicable (i.e., the leading eigenvalue of $\Pi$ is
unique and simple), then it always determines, whether $\Pi$ is an
s.l.p. or not.
\begin{prop}\label{p1.pmin}
Assume  the Perron-Frobenius eigenvalue of~$\Pi$ is unique and
simple. If the assumption of the algorithm is wrong (i.e., $\Pi$
is not an s.l.p.)  then for every $j = 1, \ldots , n$
condition~(\ref{cond.pmin}) is violated at some step. Conversely,
if condition~(\ref{cond.pmin}) is violated at some step for
some~$j$, then $\Pi$ is not an s.l.p.
\end{prop}
{\tt Proof.} The sufficiency follows from Lemma~\ref{l.5min}. To
show the necessity, we assume that $\check \rho (\widetilde \cM)\, <
\, 1$. Then there is a product $C \in \widetilde \cM^s$ such that
$\rho(C) < 1$. This yields $C^{\, r} \to 0$, and hence $ C^{\, r}v_1\to 0$
as $r\to \infty$. Consequently, for every $j$ one has $(v_j^*\, ,
\, C^{\, r}v_1)  < 1$, whenever $r$ is large enough. Since the point
$C^{\, r}v_1$ belongs to the infinite polytope~$P_{rs}$, we see that
$\inf\limits_{x \in P_{\, rs}}(v_j^*, C^{\, r}v_1)< 1$. This infimum is
attained at some vertex $v$ of $P_{\, rs}$, hence that vertex
violates condition~(\ref{cond.pmin}).

{\hfill $\Box$}
\smallskip

Now let us analyze estimate~(\ref{final.pmin}). In contrast to the
algorithms for JSR computation, the lower bound $(t_N)^{-1}[\rho
(\Pi)]^{1/n}$ may not converge to $\check \rho$ at all, even if
the family $\cM$ is positively irreducible. There are simple
examples already in the dimension~$d=2$. The reason is that some
product $A \in \cM^l$ may have the leading eigenvector~$v$ on the
boundary of the invariant cone~$\re^d_+$ i.e., have some zero
entries. If the corresponding leading eigenvalue is unique and
simple, then for each $j=1, \ldots , n$ the sequence $A^rv_j$
converges to $t_jv$ as $r\to \infty$, where $t_j\ge 0$ depends
only on~$j$. This means that some vertices of $P_k$ approach
closer and closer to the boundary as $k \to \infty$. In this case
Algorithm~\textbf{(L)} is useless: it gives neither an extremal infinite
polytope nor good lower bound for $\check \rho$. We suggest two
methods to avoid this situation. The first one is to impose a {\em
second invariant cone} assumption. Second, is to enlarge the
invariant cone~$\re^d_+$ by adding {\em extra directions}. Let us
begin with the first one.
\smallskip

\subsection{Convergence of the algorithm. Case 1; the second invariant cone.}
\label{subsec:SIC}

A cone $\widetilde K$ is {\em embedded} to a cone $K$ if $(\widetilde
K\setminus\{0\}) \, \subset \, {\rm int}\, K$. The pair $(\widetilde K,
K)$ will be refereed as an {\em embedded pair}.
\begin{defi}\label{d.pmin1}
An embedded pair is invariant for a family $\cM$ if both its cones are invariant for this family.
\end{defi}
An embedded invariant pair of cones for a given family will be
called an {\em invariant pair}. For a family~$\cM$ of
nonnegative operators (i.e., operators defined by nonnegative
matrices) we say that~$\widetilde K$ is a {\em second invariant cone}
if $(\widetilde K, \re^d_+)$ is an invariant pair. So, the cone
$\widetilde K$ is embedded in $\re^d_+$ and invariant for~$\cM$. A
simple sufficient condition for the existence of a second
invariant cone is the so-called {\em eventual positivity} of a
matrix family.
\begin{defi}\label{d.pmin2}
A nonnegative  family $\cM$ is called eventually positive if there
is $k$ such that all matrices of the family $\cM^r$ are positive
for all $r \ge k$.
\end{defi}
In particular, if all matrices of~$\cM$ are positive, then $\cM$
is eventually positive. The following trivial fact clarifies the
notion of eventual positivity:
\begin{lemma}\label{l.pmin8}
A family $\cM$ is eventually positive iff its matrices have
neither zero columns nor zero rows, and there is $k$ such that all
matrices of $\cM^k$ are positive.
\end{lemma}
\begin{lemma}\label{l.pmin10}
Every eventually positive family possesses a second invariant cone.
\end{lemma}
{\tt Proof.} A conic hull of the set $\cup_{A \in \cM^k}AK$, where
$K = \re^d_+$ is the second invariant cone for~$\cM$.

{\hfill $\Box$}
\smallskip

The key property of embedded pairs is formulated in the following lemma.
\begin{lemma}\label{l.pmin15}
If $(\widetilde K, K)$ is an embedded pair, then every antinorm on $K$ is continuous and strictly positive on~$\widetilde K$.
\end{lemma}
{\tt Proof.} Any nonnegative concave function,
which is not an identical zero, is continuous and positive at any internal point of its domain.

{\hfill $\Box$}
\smallskip

Now we can prove the existence of {\em invariant antinorms} in the interior cone.
\begin{defi}\label{d.pmin3}
Let a family~$\cM$ has a common invariant cone $K$. An antinorm $f: K \to \re_+$ is called invariant, if
there is a constant $\lambda \ge 0$ such that
$$\min\limits_{j=1, \ldots , m}f(A_jx)\, = \, \lambda \, f(x)\, , \ x \in K.$$
\end{defi}
\begin{theorem}\label{th20.5}
If the family $\cM$ possesses an embedded pair~$(\widetilde K, K)$, then it has
a positive monotone invariant antinorm on~$\widetilde K$. For any invariant antinorm
on~$\widetilde K$ we have~$\lambda = \check \rho(\cM)$.
\end{theorem}
Let us recall that for extremal antinorms we have 
$\min_{j=1, \ldots , m}f(A_jx)$\, $\ge \, \check \rho \, x.$
It becomes an equality if the antinorm is invariant.


By Theorem~\ref{th20} an extremal antinorm always exists, whenever the operators
share an invariant cone. For the invariant antinorm this is not the case. There are simple examples of irreducible
pairs of nonnegative $2\times 2$-matrices that do not have an invariant antinorm. So, the embedded pair assumption is essential
in Theorem~\ref{th20.5}. In the proof of Theorem~\ref{th20.5} we use the following simple fact from~\cite[section 4]{P3}.
\begin{lemma}\label{l.tree3p}
For every pair of embedded cones $(\widetilde K, K)$ and for every norm
in~$\re^d$ there is a constant~$\gamma$ such that for any
operator~$B$ with these invariant cones and for each $x \in \widetilde
K$ one has $\|Bx\| \, \ge \, \gamma \, \|B\|\, \|x\|$.
\end{lemma}

{\tt Proof of Theorem~\ref{th20.5}.} Let $f$ be an invariant antinorm. For each $k$ we have
$\min_{A \in \cM^k}f(Ax)\, = \, \lambda^k f(x)$. On the other hand, combining  Lemmas~\ref{l10} and \ref{l.tree3p}
we see that for any operator $B$ that preserves the cones $K$ and $\widetilde K$ one has
$$C_0\|B\|\, \|x\| \, \ge \, f(Bx)\, \ge \, C_0^{-1}\|B\|\, \|x\|\, , \ x \in \widetilde K,$$
where $C_0$ does neither depend on $B$ nor on $x$.


A  nontrivial concave nonnegative function is strictly positive on the
interior of its domain. Therefore, $f(x) > 0$ for all $x \in \widetilde K \setminus \{0\}$.
By the compactness argument it follows that $f(z) \ge C_1 \|z\|$ for any $z \in \widetilde K$, where
$C_1 > 0$ does not depend on $z$. Hence, applying Lemma~\ref{l.tree3p}, we conclude that
$f(Bx) \ge C_1 \|Bx\| \ge C_1 \gamma \|B\| \, \|x\|$.

Therefore $\min_{A \in \cM^k}\|Ax\|\, \asymp \, \lambda^k$, and therefore
$\check \rho = \lambda$. Now let us prove the existence of an invariant antinorm. It suffices to consider the case
$\check \rho (\cM) = 1$. By Theorem~\ref{th20} there is a monotone  extremal antinorm $f_0$, for which
$\min_{A \in \cM}f_0(Ax)\, \ge  \,  f_0( x)\, , \, x \in \widetilde K$. Let $f_{j}(x)\, = \, \min_{A \in \cM}f_{j-1}(Ax)\, , \
j \in \n$. This is a nondecreasing sequence of antinorms. If for some $\, x \in \widetilde K\, $ we have
$\, f_j(x)\, \to \, +\infty\, $ as $\, j \to \infty$, then this holds for all nonzero $x \in \widetilde K$, and hence
$\check \rho (\cM) >1$. Thus, the sequence $\{f_j\}_{j\in \n}$ is bounded, hence it converges pointwise
to some function~$f$, which is a monotone invariant antinorm.

{\hfill $\Box$}
\smallskip

Applying now invariant antinorms we can prove the convergence results for Algorithm~\textbf{(L)}.
We start with inequality~(\ref{final.pmin}).

\begin{prop}\label{p2.pmin}
If the family $\cM$ possesses a second invariant cone~$\widetilde K
\subset \re^d_+$, and $v_1 \in \widetilde K$, then for
estimate~(\ref{final.pmin}) we have $t_N \to 1$ as $l \to
\infty$ and $N \to \infty$.
\end{prop}
Thus, if $\cM$ has a second invariant cone, then
Algorithm~\textbf{(L)} is always applicable for, at least, approximate
computation of~LSR. It either finds the value of LSR or provides lower and upper bounds for it; those bounds are arbitrarily
close to each other, whenever both $l$ and $N$ are large enough.

In the proof we use Dini's theorem (see~\cite[theorem 7.13]{Ru}) and the following  analogue  of the  Minkowski norm
for concave functionals. We call a convex closed set $D \subset \re^d_+$ {\em admissible} if it does not contain the origin, and if with every point $x \in D$ it contains all points $y \ge x$. In particular, all infinite polytopes not containing the origin
are admissible.
The {\em Minkowski antinorm}  associated to an admissible set $D \subset \re^d_+$ is defined as
  $\, f_D(x) \, = \, \sup \, \bigl\{ t^{-1}  \ \bigl| \ t >0, \  t x \in D  \bigr\}\, , \ x \in \re^d_+$. For any admissible set $D$ the function $f_D$ is a positive monotone antinorm on~$\re^d_+$.

\noindent {\tt Proof of Proposition~\ref{p2.pmin}}. Assume first that $\check \rho (\widetilde \cM)= 1$, i.e., that $\Pi$ is an s.l.p.
The algorithm produces the infinite polytopes
$\{P_k\}_{k \in \n}$ such that $P_{k} \, \subset \, P_{k+1} \, =
\,  {\rm co}_{+}\, \bigl\{
 \widetilde A_1P_k, \ldots ,  \widetilde A_mP_k\bigr\}$. All their vertices are in~$\widetilde K$, because $v_1 \in \widetilde K$.
By Theorem~\ref{th20.5} there is a positive invariant antinorm $f_0$ on~$\widetilde K$, for which
$$f_0(x) = \min_{j=1, \ldots , m}f_0(\widetilde A_jx)\, , \, x \in \widetilde K.$$ Therefore,
for all vertices $v$ of the polytopes $P_k$ one has $f_0(v) \ge f_0(v_1)$. Since $f_0$ is positive on~$\widetilde K$,
by the compactness argument it follows that $f_0(x) \ge C\|x\|\, , \, x \in \widetilde K$, where
$C > 0$ is a constant. Whence, all polytopes $P_k$ are uniformly separated from zero: they do not intersect the ball of radius $C\|v_1\|$ centered at the origin. Consequently, the sequence $\{f_k\}_{k \in \n}$ of Minkowski antinorms
generated by the infinite polytopes $\{P_k\}_{k \in \n}$ is non-decreasing and bounded. So, it converges pointwise
to a positive monotone antinorm~$f$.
 By Dini's theorem~~\cite[theorem 7.13]{Ru},
 this convergence is uniform on the set
 $S = \{ x \in \widetilde K \ | \ f(x) = 1\}$. Thus, $f_k(x) \to 1$ uniformly for $x \in S$, as $k \to \infty$. Hence, there is $N_{\varepsilon}$ such that $f_{N-1}(x) \ge 1-\varepsilon$ for all $x \in S$, whenever $N \ge N_{\varepsilon}$.
 Hence, $(t_N)^{-1} \, = \, \inf_{x \in P_N}f_{N-1}(x) \, \ge \, \inf_{x \in S}f_{N-1}(x)\, \ge \, 1-\varepsilon$,
 which completes the proof for the case $\check \rho (\widetilde M)= 1$. The transfer to the general case is realized in the same way as in the proof of Proposition~\ref{p2.r}.

{\hfill $\Box$}
\smallskip

\begin{cor}\label{c210}
Each of the following conditions is sufficient for the convergence
$t_N \to 1$ as $l \to \infty$ and $N \to \infty$:

1) the family~$\cM$ is eventually positive;

2) the family~$\cM$ has a second invariant cone, and the leading
eigenvalue of $\Pi$ is simple.
\end{cor}
{\tt Proof.}
Observe that if a family is eventually positive, then every product $\Pi$
of its matrices has a unique simple largest by modulo eigenvalue. This eigenvalue is positive, and
the corresponding eigenvector is strictly positive. To see this note that the matrix $\Pi$ is obviously
primitive (i.e., it is nonnegative and some power $\Pi^k$ is strictly positive). A primitive matrix
always has a unique simple largest by modulo eigenvalue, which is positive, and the corresponding
eigenvector is strictly positive~\cite[chapter 8]{HJ}.
If $\cM$ is eventually positive, then for each $A \in
\cM^k$ and every $v \in \re^d_+$ the vector  $Av$ belongs to the
interior cone $\widetilde K$, which is a conic hull of the
set~$\bigcup_{d_k, \ldots , d_1}\, A_{d_k}\cdots A_{d_1}(\re^d_+)$.
Hence, $\widetilde K$ contains the leading eigenvector of any
product~$\Pi$ of matrices from~$\cM$.


If~$\cM$ has a second invariant cone~$\widetilde K$, then, by the
Perron-Frobenius theorem, $\widetilde K$ contains some of the leading
eigenvectors of~$\Pi$. If the leading eigenvalue is simple,
then~$v_1 \in \widetilde K$.

{\hfill $\Box$}
\smallskip

\subsection{The criterion for finite termination of Algorithm~\textbf{(L)}.}
\label{sub:FTL}

Now we are ready to prove  a sharp criterion ensuring that
Algorithm~\textbf{(L)} produces an extremal infinite polytope. It looks similar to Theorem~\ref{th.cond-r}
and use the notion of under-dominant product.

\begin{defi}\label{d.simple-min}
A product $\Pi\in \cM^n$ is called under-dominant for the family $\cM$
if there is $p> 1$ such that the spectral radius of every product of
operators of $\widetilde M$, that is not a power of $\widetilde \Pi$ nor a
power of its cyclic permutations,  is bigger than $\, p$.
\end{defi}

\begin{theorem}\label{th.cond-min}
Assume the family $\cM$ is eventually positive.
Algorithm~{\rm \textbf{(L)}} terminates within finitely many iterations
if and only if $\Pi$ is under-dominant for $\cM$.
\end{theorem}
The proof of Theorem~\ref{th.cond-min} is in Appendix \ref{sec:app}.
\begin{remark}\label{r.30-min}
{\rm If the family~$\cM$ is eventually positive, then every
product $\Pi \in \cM^n$ is a {\em primitive matrix}, i.e., some of
its powers is positive.


This is well known that the leading
eigenvalue of a primitive matrix is always unique and simple (see,
for instance,~\cite[chapter 8]{HJ}). That is why in Theorem~\ref{th.cond-min} we do
not need the uniqueness and simplicity of the leading eigenvalue
assumption, in contrast to Theorem~\ref{th.cond-r}. }
\end{remark}
\begin{remark}\label{r.31-min}
{\rm Let us stress again that Algorithm~\textbf{(L)} can produce
extremal polytopes for nonnegative families that have no second
invariant cone or not eventually positive. We will see some
examples in Section~\ref{sec:IE} and~\ref{sec:Appl}.
 The only difference is that we have not
succeeded in finding a reasonable criterion for that case. As for
Theorem~\ref{th.cond-min}, the eventual positivity  assumption is
essential and cannot be omitted.}
\end{remark}

\begin{cor}\label{c200-min}
If an eventually positive  family $\cM$ possesses  an under-dominant
product,  then it has an extremal infinite polytope.
\end{cor}

\subsection{Case 2. Modification of Algorithm~\textbf{(L)}}
\label{subsec:ML}

In the previous subsection  we showed that if the family~$\cM$ has a second
invariant cone (in particular, if this family is eventually
positive), then the algorithm is always applicable. It may
converge within finite time, in which case it produces the
extremal infinite polytope and finds the exact value of LSR. By
Theorem~\ref{th.cond-min} this happens precisely when the
product~$\Pi$ is under-dominant. Otherwise, if it does not converge,
it produces  upper and lower bounds in~(\ref{final.pmin}) that both
tend to $\check \rho (\cM)$ as $N\to \infty$
(Proposition~\ref{p2.pmin}). If~$\cM$ does not have the second
invariant cone, then the algorithm can be applied  as well, but in
some cases it may not converge to the value of LSR. This happens,
for instance, when there is a product $A \in \cM^r$, whose leading
eigenvalue $\lambda_{\max}$ is unique and simple, and the
corresponding eigenvector~$v$ has some zero entries. If
$[\lambda_{\max}]^{1/r} > [\rho(\Pi)]^{1/n}$ then, for each point
$v_i$ produced by the algorithm we have $A^sv_i \to t_iv$ as $s\to
\infty$, where the sequence $t_i\ge 0$ diverges.


Whence, some vertices of the polytopes
$P_k$ converge to the boundary of~$\re^d_+$ as $k \to \infty$. In
this case the algorithm does not terminate within finite time,
since new vertices $v_i$ will always appear (closer and closer to
the boundary of~$\re^d_+$). Moreover, the ratio $t_N$ in~(\ref{final.pmin}) may not converge to $1$, and the
algorithm becomes useless. We suggest the following modification
of the algorithm for this case, which often leads to the precise
values of LSR.
\medskip

Assume the algorithm has not terminated after $N$ steps. Take some
small $\delta > 0$ and find all vectors $v \in \cV_{N}$ such that
$\frac{v_{\min}}{v_{\max}}\, < \,\frac{\delta}{d}$, where $v_{\min}$
and $v_{\max}$ are respectively the smallest and the largest entry of~$v$.
To any such a vector~$v$ we associate a vector $h$ such that
$h^q=-\varepsilon$ if $\frac{v^{q}}{v_{\max}}\, < \,
\frac{\delta}{d}$, and $h^q = 1$ otherwise (we write $x^q$ for the
$q$th entree of the vector~$x$). The parameter $\varepsilon
> 0$ is chosen to be small and the same for all $h$. The finite set of all vectors $h$ will
be denoted as $\cH$. We also denote by $e$ the vector of ones.

For every $j=1, \ldots , m$ and for every
 $\bar h \in \cH$ we solve the following LP problem:

\begin{equation}\label{LP.m}
\left\{
\begin{array}{rcl}
\max & & t_e
\\[0.25cm]
{\rm subject \ to}
     & & \widetilde A_j\, \bar h \ge t_e e\, + \, \sum\limits_{h \in \cH}\, t_h\, h
\\[0.3cm]
{\rm and}
     & & t_h \, \ge \, 0\, , \qquad  \forall h \in \cH.
\end{array}
\right.
\end{equation}

If for some $j$ and $\bar h$ we have $t_e \le 0$, then for the
chosen values of $\delta$ and $\varepsilon$ the modification is
impossible. We can try smaller values. If $t_e > 0$ for all $j$
and $\bar h \in \cH$, then
 we restart our algorithm
with the same product $\Pi$ and with the only modification:  LP
problem~(\ref{LP.pmin}) is replaced by the following LP problem,
where $z = \widetilde A v$:

\begin{equation}\label{LP.pmin-m}
\left\{
\begin{array}{rcl}
\min & & t_0
\\[0.25cm]
{\rm subject \ to}
     & & t_0\,  z \ge \sum\limits_{x \in \cV_{k}}\, t_x\, x + \sum\limits_{h \in \cH}\, t_h\, h
\\[0.3cm]
{\rm and}
     & & \sum\limits_{x \in \cV_{k}}\, t_x \ge 1, \qquad t_x \ge 0 \quad \forall x \in \cV_k,
     \qquad t_h \ge 0 \quad \forall h \in \cH.
\end{array}
\right.
\end{equation}

{\em Explanation.} If $t_e > 0$ in LP problem~(\ref{LP.m}) for
every $j=1, \ldots , m$ and for every $\bar h \in \cH$, then the
cone $\, K_{\cH} \, = \, \bigl\{\,  x + \sum_{h \in \cH} t_h h \
\bigl| \  x \in \re^d_+\, , \ t_h \ge 0\, , \ h \in \cH   \,
\bigr\}$ is invariant for the family~$\cM$. If the algorithm
terminates after the $k$th step, then the set $P_{k-1} \, = \,
{\rm co}\, (\cV_N)\, + \, K_{\cH}$ is an extremal infinite
polytope for the family~$\widetilde \cM$, i.e., $\, \widetilde A \,
P_{k-1} \, \subset \, P_{k-1}$ for all $\widetilde A \in \widetilde M$.
Therefore, it defines an extremal  antinorm in the cone $K_{\cH}$,
and hence $\check \rho (\widetilde \cM)= 1$. Thus, we replace the
invariant cone $\re^d_+$ by a wider invariant cone $K_{\cH}$,
which covers those vertices $v\in \cV_N$ that come too close to
the boundary of $\re^d_+$. In many practical cases this trick
makes the algorithm converge within finitely many steps.
We use it in the proof of Theorem~\ref{th100} in \S \ref{subsec:OFW}.

\smallskip
We consider in the sequel several numerical examples of implementation of our
algorithms, and start with the simplest case of nonnegative
$2\times 2$-matrices. In Example~\ref{ex.jsr1}
Algorithm~\textbf{(P)} finds the JSR of two matrices, in
Example~\ref{ex.lsr1} Algorithm~\textbf{(L)} finds the LSR of
another pair of matrices. The aim of those examples is to show how the
algorithms work. Then in Section~\ref{sec:Appl} we apply our algorithms to
matrices of bigger dimensions (up to $d=50$) arising in various
problems of combinatorics and number theory. Finding the exact
values of JSR and LSR we prove, in particular, several previously
stated  conjectures in combinatorics and number theory,  and disprove one conjecture on Pascal's rhombus.
Further, in Section~\ref{sec:numer}, we show the statistics how our algorithms work for randomly
generated matrices of various dimensions. For all randomly generated families the algorithms found
the exact values of JSR and of LSR (the latter is in the case of nonnegative families).

\section{ Illustrative examples}
\label{sec:IE}

We consider here some simple examples showing the flow of the algorithms we have presented
both for the joint spectral radius and for the lower spectral radius of a nonnegative
set of matrices.

\begin{ex}\label{ex.jsr1}.{\em \textbf{Computation of the joint spectral radius.}
Consider a family ${\cM}=\{ A_1,A_2 \}$:
\begin{eqnarray*}
\displaystyle{ A_1 \, = \, \left(
\begin{array}{rr}
1 & 1 \\
0 & 1
\end{array}
\right) \qquad \mbox{and} \qquad A_2 \, = \, b\,  \left(
\begin{array}{rr}
1 & 0 \\
1 & 1
\end{array} \right), } &&
\end{eqnarray*}
with $b=9/10$. Looking over all matrix products up to some length,
we make a guess that $\Pi\, =\, A_1\, A_2\, $ is a spectrum
maximizing product. We have~$[\rho(\Pi)]^{1/2} \, = \, \sqrt{b}\,
\frac{1+\sqrt{5}}{2}$. Applying then  Algorithm~\textbf{(P)} to the
family~$\widetilde \cM\, = \, [\rho(\Pi)]^{-1/2} \cM$ we obtain at step zero
 a unique leading eigenvector $v_1$ of the product $\widetilde \Pi_1
= \widetilde \Pi$,  and $v_2 = \widetilde A_2 v_1$ -- the leading eigenvector of the cyclic permutation~$\widetilde \Pi_2$:
$\, v_1 \, = \, \bigl( 1 \ , \ \frac{\sqrt{5}-1}{2}\bigr)\, $ and
$$
v_2 \, = \widetilde A_2 v_1 \, = \, \sqrt{b} \biggl( \frac{\sqrt{5}-1}{2} \ , \ 1 \biggr) \, = \, \bigl( 0.586318522 \ , \ 0.948683298 \bigr)\,
$$
(all the values are rounded  to the ninth decimal). At the first step we get one new point
$$
v_3
\, = \, \widetilde A_1 v_1 \, = \, \frac{1}{\sqrt{b}} \biggl( 1 \ , \ \frac{\sqrt{5}-1}{\sqrt{5}+1} \biggr) \, = \, \bigl( 1.054092553 \ , \ 0.402627528 \bigr)\, ,
$$
and the other point $\, v_4 \, = \,
\widetilde A_2 v_2\, $ is ``dead'', because it belongs to the interior of $\, P_1 \, = \, {\rm co}_{-} \,
\{v_1, v_2, v_3\}$, i.e., solving LP problem~(\ref{LP.p}) for the
point $v=v_2$ and for $\widetilde A \, = \, \widetilde A_2$ we get
$t_{\{v_2, \widetilde A_2\}} \, > \, 1$. Thus, after the first step $\cV_1 \, = \,
\{v_1, v_2, v_3\}\, $ and $\, \cU_1 \, = \, \{v_3\}$.

At the second step we solve LP problem~(\ref{LP.p}) for the pairs $(v_3,
\widetilde A_1)$ and $(v_3, \widetilde A_2)$ and find that the values $\,
t_{\{v_3, \widetilde A_1\}}$ and $\, t_{\{v_3, \widetilde A_2\}}$ are both
bigger than~$1$. This means that the points $\widetilde A_1v_3$ and
$\widetilde A_2v_3$ are both internal to~$P_1$. Therefore, $\widetilde A_j P_1
\subset P_1\, , \ j = 1,2$, and so $\, P_1$ is an extremal
polytope (see Figure \ref{fig:jsr}).

\begin{figure}[ht]
\centering
\hskip 0.1cm \global\def\path{#1}\input{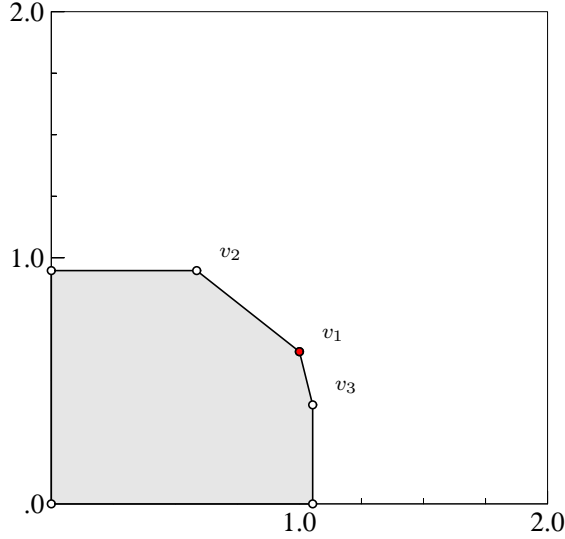} \\[2mm]
 \caption{The extremal  polytope ${P_1}$ of Example~\ref{ex.jsr1}.
 The starting leading eigenvector $v_1$ of $\Pi$ is indicated in red. \label{fig:jsr}}
\end{figure}
Thus, the algorithm terminates after the second step, and $\, \widehat \rho (\cM)\, = \,
\bigl[ \rho(\Pi) \bigr]^{1/2}\, = \, \sqrt{b}\, \frac{1+\sqrt{5}}{2}$.

The cyclic tree of this algorithm is plotted in
Figure~\ref{fig:cyc1}.
\begin{figure}[ht]
\centering
\global\def\path{#1}\input{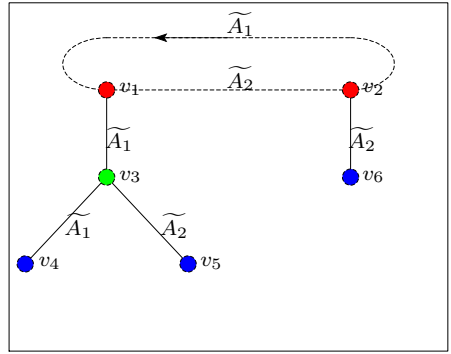}
 \caption{The cyclic tree of Example~\ref{ex.jsr1}. The root $\{v_1,v_2\}$ is in red,  the alive leaves are in green,
the dead leaves are blue. \label{fig:cyc1}}
\end{figure}
In Figure \ref{fig:jsr2} we plot the
points $\{\widetilde A_1 v_i\}_{i=1}^{3}$ (in red) and
$\{\widetilde A_2 v_i\}_{i=1}^{3}$ (in blue).
\begin{figure}[ht]
\centering
\global\def\path{#1}\input{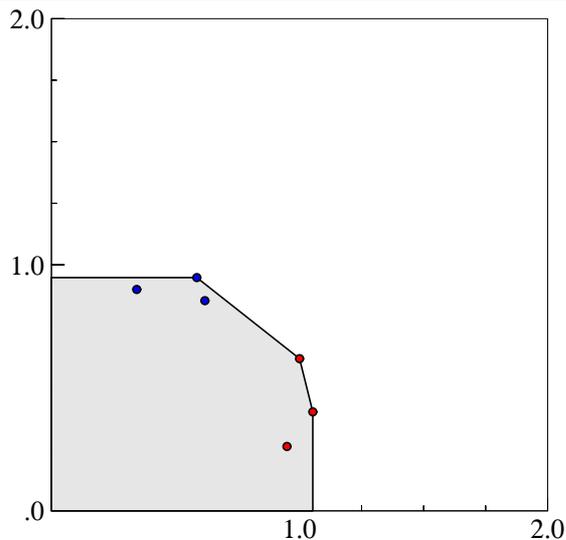} \\[2mm]
 \caption{The extremal polytope of Example~\ref{ex.jsr1} and the transformed vectors
 $\{ \widetilde A_1 v_i \}_{i=1}^{3}$ (in blue) and $\{ \widetilde A_2 v_i \}_{i=1}^{3}$
 (in red). \label{fig:jsr2}}
\end{figure}
}
\end{ex}

\begin{figure}[ht]
\centering
\global\def\path{#1}\input{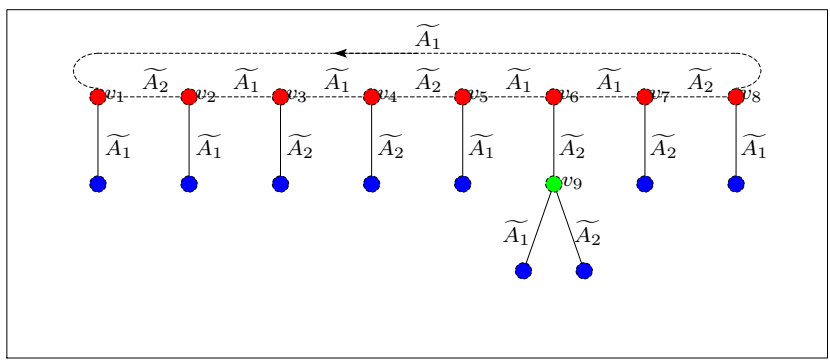}
 \caption{The cyclic tree of
  Example~\ref{ex.lsr1} (the red root,  green alive leaves and
 blue dead leaves). \label{fig:cyc2}}
\end{figure}
\begin{figure}[ht]
\centering
\global\def\path{#1}\input{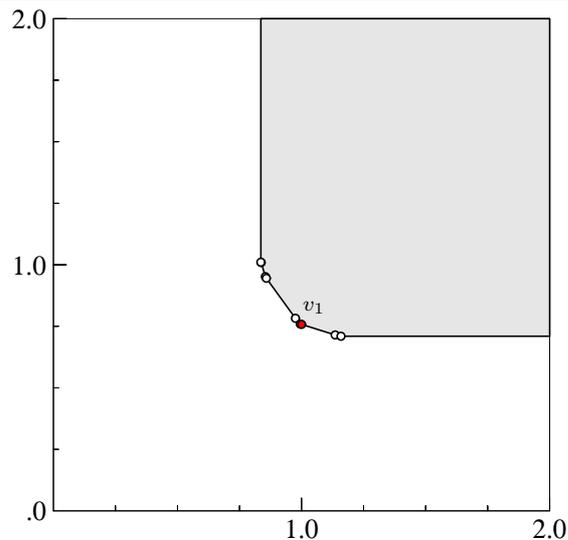} \\[2mm]
 \caption{The extremal infinite polytope in Example~\ref{ex.lsr1} (the starting eigenvector $v_1$ of $\Pi$ is red). \label{fig:pol}}
\end{figure}

\begin{ex}\label{ex.lsr1}.{\em \textbf{Computation of the lower spectral radius.}
Let ${\cM}=\{ A_1, A_2 \}$ with
\begin{eqnarray*}
A_1 \ = \ \left(\begin{array}{rr}
7 & 0
\\
2 & 3
\end{array} \right), \quad
A_2 \ = \  \left(\begin{array}{rr}
2 & 4
\\
0 & 8
\end{array} \right)\, .
\end{eqnarray*}
We prove that the product $\, \Pi \ = \ A_1\,A_2\,(A_1^2\, A_2)^2$
is  spectrum minimizing, and hence the LSR $\, \check \rho (\cM)$
equals to
$\, \rho(\Pi)^{1/8} \ = \ \left( 4 \left(213803+\sqrt{44666192953}\right) \right)^{1/8} \ = \ 6.009313489\ldots $.

Starting  Algorithm~\textbf{(L)} we define at zero step  $\widetilde M
\, = \, [\rho(\Pi)]^{-1/8}\cM\, , \, \widetilde \Pi \, = \,
[\rho(\Pi)]^{-1}\Pi$ and  get eight points $v_1, \ldots , v_8$
starting from the leading eigenvector
$$
v_1 \, = \, \biggl( 1 \ , \ \frac{97444}{\sqrt{44666192953}-82749} \biggr) \, = \, \bigl( 1 \  , \ 0.757760157 \bigr)
$$
of $\widetilde \Pi$. Thus, $\, v_2\, =\, \widetilde
A_2 v_1$, $\, v_3\, =\, \widetilde A_1v_2$, $\, v_4\, =\, \widetilde A_1
v_3$, $\, v_5\, =\, \widetilde A_2 v_4$, $\, v_6\, =\, \widetilde A_1
v_5$, $\, v_7\, =\, \widetilde A_1 v_6$, $\, v_8\, =\, \widetilde A_2
v_7$. At the first step we solve eight LP problems~(\ref{LP.pmin})
and get the only new alive vertex $v_9 \, = \, \widetilde A_2v_6$.

The other seven new vertices are ``dead leaves'': they belong to the
interior of the infinite polytope $P_1 \, = \,  {\rm co}_{+} \, \{v_i\}_{i=1}^9$,
for each of them the value $t_{\{v, \widetilde A\}}$ of problem
(\ref{LP.pmin}) is smaller than~$1$. Thus, after the first step
we have $\cV_1 = \{v_i\}_{i=1}^9$ and $\cU_1 = \{v_9\}$. The
next step produces two new vertices $\widetilde A_1v_9\, , \, \widetilde
A_2v_9$ and they are both dead. Thus, the algorithm terminates
after the second step, and $P_1$ is the extremal infinite
polytope.

\begin{figure}[ht]
\centering
\global\def\path{#1}\input{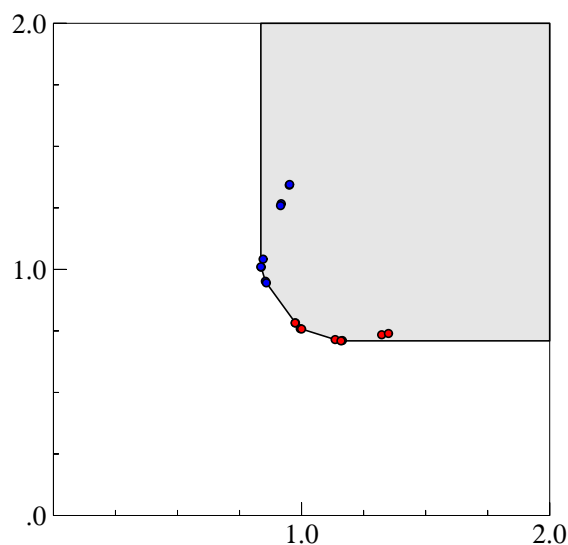} \\[2mm]
 \caption{The extremal infinite polytope $P_1$ and the transformed vectors
 $\{ \widetilde A_1 v_i \}_{i=1}^{9}$ (in red) and $\{ \widetilde A_2 v_i \}_{i=1}^{9}$ (in blue) of Example~\ref{ex.lsr1}.
 \label{fig:pol2}}
\end{figure}
}
\end{ex}

\section{ Applications}
\label{sec:Appl}

We consider four applications of Algorithms~\textbf{(P)}
and~\textbf{(L)} to various problems of combinatorics, number
theory, and theory of formal languages. Each problem is reduced to
computing JSR or LSR of some families of nonnegative matrices, which
we are able to solve exactly in the sense specified in previous sections.
\smallskip

\subsection{The  asymptotics of the number of overlap-free words}
\label{subsec:OFW}

The  problem of counting of overlap-free binary words was
intensively studied in the literature (see the recent
survey~\cite{Be}). In~\cite{C} and then in~\cite{JPB2}
this problem was reduced to computing  JSR and LSR of two
special nonnegative $20\times 20$-matrices. Those values were
computed approximately, and two conjectures were stated about
their  exact values~\cite{JPB2}. Now we prove both those conjectures by
applying Algorithms~\textbf{(P)} and \textbf{(L)}.

A binary word, i.e., a finite sequence of zeros and ones, is called
{\em overlap-free} if it does not contain a subword of the form $
xaxax, $ where $x\in \{0,1\}$ and $a$ is a word. In 1906 Thue
 proved  that there are infinitely many such
words. A natural problem, which was analyzed in many papers,  is
to estimate the total number $u_n$ of overlap-free words of
length~$n$. In~1988 Brlek showed that $u_n\, \geq \, 3\,
n\, - \, 3$, on the other hand  Restivo and Salemi in 1985  proved
the polynomial upper bound $u_n\ \leq \ C \, n^{\, r}$, where
$r=\log(15) \approx 3.906$.   This result was sharpened
successively by Kfoury (1988), Kobayashi~(1988), and
 Lepist\"o~(1995) to the value $r=1.37$. On the other hand,
Kobayashi~(1988)  showed that $ u_n \ \ge \  C\, n^{\,
1.155}$ (see~\cite{Be} for the corresponding references and historical overview). So, the number of overlap-free words grows faster than
linearly. A natural question arises, whether $\, u_n \, \asymp \, n^{\, \gamma}$
for some $\gamma \in [1.155\, , \, 1.37]$~.
Cassaigne~\cite{C} showed that the answer is negative. He
introduced the lower and the upper exponents of growth:
\begin{eqnarray}\label{alphabeta}
\alpha & = & \sup  \, \bigl\{ r \ \bigl|  \ \exists  \ C > 0\, ,  u_n  \, \ge \,  C \, n^r  \, \bigr\},\\
\nonumber \beta & = & \, \inf  \, \bigl\{ r \ \bigl|  \ \exists  \ C > 0\, ,  u_n  \, \le  \, C \, n^r  \, \bigr\},
\end{eqnarray}
and proved that $\, \alpha < \beta$. Moreover, he established that
the numbers~$u_n$ can be computed as sums of variables that are
obtained by certain linear recurrence relations.  This led  to
 the following bounds: $\alpha<1.276$ and $\beta> 1.332$.
The next improvement is due to Jungers, Protasov and Blondel~\cite{JPB2},
who,  showed that $\, \alpha\, = \, \log_2 \check \rho (A_1,
A_2)\, $ and $\, \beta\, = \, \log_2 \widehat \rho (A_1, A_2)$, where
$A_1$ and $A_2$ are special $20\times 20$-matrices with
nonnegative integer entries, which are reported in the Appendix \ref{sec:app}.

 In~\cite{JPB2} the authors introduced new algorithms for estimating the joint
and lower spectral radii based on the convex programming; by means of these algorithms
they derived the following bounds:
\begin{equation}\label{ours}
1.2690 \ <\ \alpha \ <\ 1.2736   \qquad \mbox{and} \qquad  1.3322
\ <\ \beta\ <\ 1.3326\, .
\end{equation}
This allowed the authors to make the following conjectures on the precise values:
\begin{conj}\label{conj.n1}~\cite{JPB2}

The s.m.p. for the family $\cM$ is $A_1A_2$, and $\, \beta\, =
\, \frac{1}{2}\, \log_2 \, \rho(A_1A_2)$.
\end{conj}
\begin{conj}\label{conj.n2}~\cite{JPB2}

The s.l.p. for the family $\cM$ is $A_1A_2^{10}$, and  $\,
\alpha \, = \, \frac{1}{11}\, \log_2 \, \rho(A_1A_2^{10})$.
\end{conj}
Algorithms~\textbf{(P)} and \textbf{(L)} now make it possible to prove both
these conjectures.
\begin{theorem}\label{th100}
For the upper and lower exponents of growth of the function $u_n$
one has:
\begin{eqnarray*}
\alpha & = & \frac{1}{11}\,
\log_2 \, \rho(A_1A_2^{10})\ = \  1.273553265\ldots\, .
\\
\beta  & = & \frac{1}{2}\, \log_2 \, \rho(A_1A_2)\, =
\,  \, 1.332240491\ldots .
\end{eqnarray*}
\end{theorem}
Thus, both the upper and the lower exponents of asymptotic growth of the overlap-free
words can be found precisely. To prove Theorem~\ref{th100} it suffices to present
the corresponding extremal polytopes.
Since the matrices $A_1 A_2$ are nonnegative, one can  apply
 Algorithm~\textbf{(P)} for the JSR computation. The candidate for s.m.p. is $\Pi =
A_1A_2$. The algorithm terminates having performed $k=10$ steps.
Thus, $\widehat \rho (A_1, A_2)\, = \, \sqrt{\rho(A_1A_2)}\, = \,
2.517934040\ldots$. The extremal polytope $P_9$ has~$54$ vertices.

To compute the LSR we apply Algorithm~\textbf{(L)}.
The candidate for s.l.p. is $\Pi = A_1A_2^{10}$. However,
performing $k=10$ steps we see that the algorithm does not converge.
There are two sequences of vertices~$v_i \in \cV_k$ that approach
to the boundary of the positive orthant~$\re^d_+$, i.e., those
points have very small entries on some positions. Therefore, we
apply the modified version of the algorithm (\S \ref{subsec:ML}). Taking
$\delta = 1/200$, we see that one sequence have entries of
index~$q\in I_1 = \{5,10,17,18\}$ smaller than $\delta$, the
other sequence have entries of index~$q\in I_2 =
\{7,8,15,20\}$ smaller than $\delta$. We take $\varepsilon = 1/4$
and $\cH\, = \, \{h_1, h_2\}$, where $h_j \in \re^{20}$, the $q$th
entry of $h_j$ is $-\varepsilon$ if $q \in I_j$ and $q = 1$
otherwise, $j=1,2$. Solving LP problem~(\ref{LP.m}) for all $\bar
h \in \cH\, , \, A \in \{A_1, A_2\}$  we obtain  $t_e
>0$ for each of those four problems, and hence $\, K_{\cH}\, = \, \bigl\{x + t_1 h_1 + t_2 h_2 \ \bigl| \
x \in \re^d_+\, , \ t_1 , t_2 \, \ge 0\, \bigr\}$ is a common invariant cone for $A_1,
A_2$. Now we apply again Algorithm~\textbf{(L)} with the cone
$K_{\cH}$ instead of $\re^d_+$, i.e., replacing LP
problem~(\ref{LP.pmin}) by~(\ref{LP.pmin-m}). The modified
algorithm converges: it terminates after $k=15$ steps. Thus, $\Pi
= A_1A_2^{10}$ is an s.l.p., and $\check \rho (A_1, A_2) \, = \,
\bigl[ \rho(A_1A_2^{10})\bigr]^{1/11}\, = \, 2.417562630\ldots$. The extremal
infinite polytope $P_{14}$ has $104$ vertices.
\smallskip

In order to give the formal proof of the theorem it is sufficient to simply provide the
list of vertices of the corresponding extremal polytope $P_9$ (to compute JSR) and
of the extremal infinite polytope~$P_{14}$ (to compute LSR). They are given in the Appendix \ref{sec:app}.
\smallskip

\subsection{The  density  of ones in the Pascal rhombus}
\label{subsec:PR}

The Pascal rhombus is related to the Pascal triangle, with the
only difference that each element equals to the sum of four
previous elements rather  than two (see~\cite{GKMT} for definitions and
basic properties). The elements of the Pascal rhombus arise from
linear recurrence relations on polynomials. The sequence of
polynomials $\{p_n\}$ is defined as $p_0(x)=1\, , \, p_1(x) \ = \,
x^2 + x+1\, $ and $\, p_n(x)\, = \, (x^2 + x+1)p_{n-1}(x)\, + \,
x^2 p_{n-2}(x)\, , n \ge 2\, . $ This leads to a recurrence
relation for the number~$w_n$ of odd coefficients of $p_n$~\cite{FSB}. The
asymptotic growth of~$w_n$ as $n \to \infty$ is characterized as
follows:
$$
\limsup_{n\to \infty}\, \frac{\log \, w_n}{\log\, n} \ = \ \log_2
\, \widehat \rho \ ; \qquad \liminf_{n\to \infty}\, \frac{\log \,
w_n}{\log\, n} \ = \ \log_2 \, \check \rho\, ,
$$
where $\widehat \rho$ and $\check \rho$ are respectively the JSR and LSR of matrices $A_1, A_2$ defined  as
\begin{equation}\label{pascal}
A_1 \ = \ \left(\begin{array}{rrrrr}
0 & 1 & 0 & 0 & 0
\\
1 & 0 & 2 & 0 & 0
\\
0 & 0 & 0 & 0 & 0
\\
0 & 1 & 0 & 0 & 1
\\
0 & 0 & 0 & 2 & 1
\end{array} \right), \quad
A_2 \ = \ \left(\begin{array}{rrrrr}
1 & 0 & 2 & 0 & 0
\\
0 & 0 & 0 & 2 & 1
\\
1 & 1 & 0 & 0 & 0
\\
0 & 0 & 0 & 0 & 0
\\
0 & 1 & 0 & 0 & 0
\end{array} \right)
\end{equation}
This is shown easily that $\widehat \rho = 2$, and the main difficulty
is to compute $\check \rho$.
\begin{conj}\label{conj.pascal}~\cite{F}
For the set of matrices~{\rm (\ref{pascal})} we have $\check \rho =
\frac{\sqrt{5}+1}{2}\, = \, 1.61803...$.
\end{conj}
 This conjecture is quite natural, because the
golden number appears in many problems of combinatorics.
An approximate computation of $\check \rho (A_1, A_2)$ given in~\cite{PJB}
provided the following estimate:
\begin{equation}\label{eq.pascal}
1.6180 \quad \le  \quad \check \rho \quad \le \quad  1.6376,
\end{equation}
which rather confirms Conjecture~\ref{conj.pascal}. The upper
bound in~(\ref{eq.pascal})  is obtained by the product
$A_1^3A_2^3$, for which $\rho^{1/6}(A_1^3A_2^3)\, = \, 1.6376
...$.

Applying Algorithm~\textbf{(L)} we find the precise value of
$\check \rho$. It follows that Conjecture~\ref{conj.pascal} is not
true, and the product $A_1^3\, A_2^3$ is actually an s.l.p.
\begin{theorem}\label{th110}
For the family $\cM = \{A_1, A_2\}$ one has $\check \rho (\cM)\, = \, \rho^{1/6}(A_1^3A_2^3)\, = \, 1.6376 ...$.
\end{theorem}
This solves the problem of asymptotics of the sequence $w_n$ and
disproves the golden number conjecture.

\noindent {\tt Proof of Theorem~\ref{th110}}.
To prove the theorem it is convenient to apply Algorithm~\textbf{(L)}
to the transposed family $\cM^* \, = \, \bigl\{A_1^*, A_2^* \bigr\}$
rather than to original family~(\ref{pascal}). Of course, $\,
\check \rho (\cM^*)\, = \, \check \rho (\cM)$. We take $\Pi =
(A_1^*)^3(A_2^*)^3$ as a candidate for s.l.p. The leading
eigenvector of $\widetilde \Pi_1 = \widetilde \Pi$ is
$$
v_1 \qquad = \qquad \left( 0.39925900\ ,  \quad 0.95496725\ ,
\quad 0.79851800 \ ,  \quad  0.90909090 \ ,  \quad  1 \right)\, .
$$
At step zero we get six vertices: $\, v_2\, =\, \widetilde A_2^* v_1$,
$\, v_3\, =\, \widetilde A_2^*v_2$, $\, v_4\, =\, \widetilde A_2^* v_3$,
$v_5\, =\, \widetilde A_1^* v_4$, $\, v_6\, =\, \widetilde A_1^* v_5$.
The first step gives two new alive vertices:
  $v_7\, =\, \widetilde A_1^* v_3\, $ and $\, v_8\, =\, \widetilde A_2^* v_6$.
The algorithm terminates after the second step, since we compute that the images of the new
vertices $v_7, v_8$ lie inside the infinite polytope~$P_1 = {\rm co}_+ \{v_1, \ldots , v_8\}$.
Thus, we obtain the leading eigenvectors of all cyclic permutations of $\widetilde \Pi$ plus
two additional points ($v_7$ and $v_8$).
 One can check
directly that the infinite polytope~$P_1 = {\rm co}_+ \{v_1,
\ldots , v_8\}$, where the vertices~$\{v_i\}_{i=1}^8$ are
described above, is extremal for $\widetilde \cM^*$, i.e., $\, \widetilde
A_j^* P_1 \, \subset \, P_1 \, , \ j = 1,2$.

\begin{figure}[ht]
\centering
\global\def\path{#1}\input{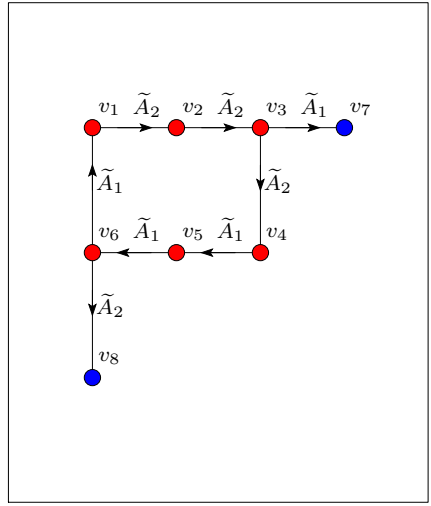}
 \caption{The vertices $\{v_i\}_{i=1}^{8}$ for LSR computation of the Pascal rhombus.
 \label{fig:plot0}}
\end{figure}

Thus, the polytope $P_1$ is extremal, and the product
$(A_1^*)^3(A_2^*)^3$ is an s.l.p. In Figure~\ref{fig:plot0} we see the cyclic tree of~$P_1$.
 Moreover, by Theorem~\ref{th.cond-min} the product $(A_1^*)^3(A_2^*)^3$ is
under-dominant.
Hence, the transpose product~$\, A_2^3\, A_1^3$ is under-dominant for~$\cM$,
and therefore so is the product~$\, A_1^3\, A_2^3\, $ being its cyclic permutation.

This concludes the proof.

{\hfill $\Box$}
\smallskip


\subsection{The Euler binary partition function}
\label{subsec:EPF}

For an arbitrary integer $r \ge 2$ the Euler binary partition function
$b(k) = b_2(r,k)$ is defined on the set of nonnegative integers
$k$ as the total number of different binary expansions $k =
\sum_{j = 0}^{\infty } d_j 2^j$, where the ''digits'' $d_j$ take
values from the set $\{0, \ldots , r-1\}$.  The asymptotic
behavior of $b(k)$ as $k \to \infty$ was studied in various
interpretations by L.~Euler, K.~Mahler, N.G.~de Bruijn, D.E.~Knuth,
B.~Reznick and others (see~\cite{P3} for the corresponding references). For even
$r = 2n$, as it was shown in~\cite{R}, one has $\, b(k)\, \asymp\,
k^{\log_2 n}$. For odd values of $r$ the asymptotic behavior of
$b(k)$ is more complicated and has been studied in~\cite{R} and~\cite{P3}. Denote
\begin{equation}\label{p1p2}
p_1 \quad =\quad \liminf_{k \to \infty }\ \log b(k)/ \log k;
\qquad p_2 \quad = \quad \limsup_{k \to \infty }\ \log b(k)/ \log
k\, .
\end{equation}

\begin{table}[thb]
\begin{center}
\begin{tabular}{|l|cccl|cccl|}\hline
 & \multicolumn{4}{c}{JSR}
  &\multicolumn{4}{c}{LSR}\\ \hline
 $r$ & $\#$ its & $\#$ vertices & $\widehat \rho$ & s.m.p. &
 $\#$ its & $\#$ vertices & $\check \rho$ & s.l.p. \\
 \hline
\rule{0pt}{9pt}\noindent
 $7$ &  5 & 8 & $3.511547$ & $A_1$ &
              6 & 14 & $3.491891$ & $A_1 A_2$ \\
 $9$ &  6 & 18 & $4.503099$ & $A_1 A_2 $ &
              5 & 17 & $4.494492$ & $A_1$ \\
 $11$ &  5 & 14 & $5.505892$ & $A_1$ &
              7 & 24 & $5.497042$ & $A_1 A_2$ \\
 $13$ &  5 & 16 & $6.502167$ & $A_1$ &
              7 & 28 & $6.498946$ & $A_1 A_2$ \\
 $15$ &  7 & 40 & $7.500106$ & $A_1 A_2$ &
              6 & 23 & $7.499841$ & $A_1$ \\
 $17$ &  7 & 40 & $8.500057$ & $A_1 A_2$ &
              6 & 30 & $8.499904$ & $A_1$ \\
 $19$ &  7 & 24 & $9.500423$ & $A_1$ &
              8 & 46 & $9.499789$ & $A_1 A_2$ \\
 $21$ &  6 & 28 & $10.500373$ & $A_1$ &
              8 & 50 & $10.499813$ & $A_1 A_2$ \\
 $23$ &  8 & 52 & $11.500053$ & $A_1 A_2$ &
              6 & 31 & $11.499894$ & $A_1$ \\
 $25$ &  9 & 34 & $12.500059$ & $A_1$ &
              8 & 58 & $12.499971$ & $A_1 A_2$ \\
 $27$ &  8 & 60 & $13.500030$ & $A_1 A_2$ &
              7 & 37 & $13.499938$ & $A_1$ \\
 $29$ &  9 & 66 & $14.500009$ & $A_1 A_2$ &
              8 & 43 & $14.499982$ & $A_1$ \\
 $31$ &  9 & 30 & $15.500001$ & $A_1$ &
              10 & 34 & $15.499999$ & $A_1 A_2$ \\
 $33$ & 11 & 36 & $16.500001$ & $A_1 A_2$ &
              10 & 55 & $16.499999$ & $A_1$ \\
 $35$ &  8 & 52 & $17.500007$ & $A_1$ &
             18 &102 & $17.499997$ & $A_1 A_2$ \\
 $37$ &  8 & 54 & $18.500012$ & $A_1$ &
              10 &113 & $18.499994$ & $A_1 A_2$ \\
 $39$ &  10 & 112 & $19.500003$ & $A_1 A_2$ &
              8 & 59 & $19.499994$ & $A_1$ \\
 $41$ &  9 & 78 & $20.500005$ & $A_1$ &
             11 &120 & $20.499997$ & $A_1 A_2$ \\
 \hline
\end{tabular} \\[2mm]
\caption{Computation of the JSR and of the LSR for the Euler partition function matrices.}\label{tab:Euler}

\end{center}
\end{table}

In~\cite{P3} it was proved that $p_{1} \, = \, \log_2 \check
\rho(A_1, A_2)\, $ and $\, p_{2} \, = \, \log_2 \widehat \rho(A_1
A_2)$, where $A_1, A_2$ are $(r-1)\times (r-1)$-matrices defined
as follows: $\, (A_s)_{ij} \, = \, 1\, $ if $\, 2-s \le 2j-i \le r-s+1$,
and $\, (A_s)_{ij} = 0$ otherwise (for $s=1,2$).
For example, for $r=7$ we have the following $6\times 6$-matrices:
$$
A_1\quad = \quad
\left(
\begin{array}{rrrrrr}
1 & 1& 1& 1& 0& 0\\
0 & 1& 1& 1& 0& 0\\
0 & 1& 1& 1& 1& 0\\
0 & 0& 1& 1& 1& 0\\
0 & 0& 1& 1& 1& 1\\
0 & 0& 0& 1& 1& 1
\end{array}
\right)\ ; \qquad
A_2\quad = \quad
\left(
\begin{array}{rrrrrr}
1 & 1& 1& 0& 0& 0\\
1 & 1& 1& 1& 0& 0\\
0 & 1& 1& 1& 0& 0\\
0 & 1& 1& 1& 1& 0\\
0 & 0& 1& 1& 1& 0\\
0 & 0& 1& 1& 1& 1
\end{array}
\right)\ .
$$
In~\cite{P3} the following conjecture was made:
\begin{conj}\label{conj.binar}
For every odd $r$ one of the two products~$A_1$ and $A_1A_2$ is an s.m.p. and the other is an s.l.p.
\end{conj}
The case $r=3$ was carried out earlier in the work~\cite{R},
for $r = 5, 7, 9, 11, 13$  Conjecture~\ref{conj.binar} was proved to hold true
in~\cite{P3}.
Algorithms~\textbf{(P)} and~\textbf{(L)} make it possible to prove
this conjecture for many more odd values (in particular we did the computation for $r \le 41$).

The results are listed in Table \ref{tab:Euler}. The first column is $r$ (where we recall that
the dimension of matrices $A_1, A_2$ is $r-1$), the second column
is the number~$k$ of iterations necessary to Algorithm~\textbf{(P)} for
terminating, the third column is the number of vertices of the extremal polytope~$P_{k-1}$, the fourth
one is the value of JSR rounded to the sixth decimal, and the fifth one is the s.m.p. The right hand
side of the table presents analogous informations for the LSR computation by
Algorithm~\textbf{(L)}.

We see that Algorithms~\textbf{(P)} and~\textbf{(L)} demonstrate a
good efficiency. Even for large dimensions of the matrices $A_1, A_2$
the total number of iterations~$k$ never exceeds 18 and the number
of vertices of the extremal polytope $P_{k-1}$ is at most 120. Let
us remark that the binary matrices $A_1, A_2$ of the partition
function are rather inconvenient for our algorithms, because of a
very small gap between JSR and LSR. For instance, for $r=33$ the
distinction between the JSR and LSR is less than $0.00002 \%$.
Therefore, all products of $A_1$ and $A_2$ of some length~$k$ have
almost the same spectral radii.
This is why we would expect Algorithms~\textbf{(P)} and~\textbf{(L)}
to need a large number of iterations.
On the contrary, they just need 11 and 10 iterations respectively.
\smallskip

Actually Algorithms~\textbf{(P)}
and~\textbf{(L)} work also for higher dimensions and
Conjecture~\ref{conj.binar} can be proved for larger~$r$.  For
instance, if~$r=51$, then Algorithm~\textbf{(L)} needs~$k=15$
iterations and produces an extremal infinite polytope~$P_{14}$
with~$135$ vertices.

\begin{figure}[t]
\centering
\global\def\path{#1}\input{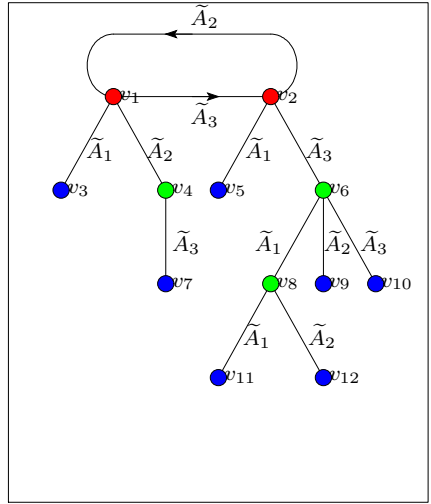}
 \caption{The vertices $\{v_i\}_{i=1}^{12}$ for JSR computation of the ternary Euler partition function.
 \label{fig:plot1}}
\end{figure}

In the next section, as further example, we consider ternary expansions and show
that also in this case we can compute the significant measures.

\smallskip

\subsection{The Euler ternary  partition function}
\label{subsec:ETF}

The LSR and JSR appear  in the problem of asymptotics of the Euler
partition function on the arbitrary base, not only for binary
expansions. For instance, the ternary partition function $b(k) =
b_3(r,k)$ is the total number of different ternary expansions $k =
\sum_{j = 0}^{\infty } d_j 3^j$, where the ''digits'' $d_j$ take
values from the set $\{0, \ldots , r-1\}$.  The largest and the
smallest exponents of growth of $b(k)$ as $k\to \infty$ are
defined  by formulas similar to~(\ref{p1p2}) through the  LSR and the JSR
of three special binary matrices $A_1, A_2, A_3$ (see~\cite{P3} for
details). In~\cite{PJB} the authors analyze the example with
$r=14$, where the matrices of $\cM = \{A_1,A_2, A_3\}$ are
\begin{eqnarray*}
A_1 & = & \left( \begin{array}{ccccccc}
     1  &   1  &   1  &   1  &   1  &   0  &   0 \\
     0  &   1  &   1  &   1  &   1  &   0  &   0 \\
     0  &   1  &   1  &   1  &   1  &   1  &   0 \\
     0  &   1  &   1  &   1  &   1  &   1  &   0 \\
     0  &   0  &   1  &   1  &   1  &   1  &   0 \\
     0  &   0  &   1  &   1  &   1  &   1  &   1 \\
     0  &   0  &   1  &   1  &   1  &   1  &   1
\end{array}  \right), \quad
A_2 \, = \, \left( \begin{array}{ccccccc}
     1  &   1  &   1  &   1  &   1  &   0  &   0 \\
     1  &   1  &   1  &   1  &   1  &   0  &   0 \\
     0  &   1  &   1  &   1  &   1  &   0  &   0 \\
     0  &   1  &   1  &   1  &   1  &   1  &   0 \\
     0  &   1  &   1  &   1  &   1  &   1  &   0 \\
     0  &   0  &   1  &   1  &   1  &   1  &   0 \\
     0  &   0  &   1  &   1  &   1  &   1  &   1
\end{array}  \right)
\\
A_3 & = & \left( \begin{array}{ccccccc}
     1  &   1  &   1  &   1  &   0  &   0  &   0 \\
     1  &   1  &   1  &   1  &   1  &   0  &   0 \\
     1  &   1  &   1  &   1  &   1  &   0  &   0 \\
     0  &   1  &   1  &   1  &   1  &   0  &   0 \\
     0  &   1  &   1  &   1  &   1  &   1  &   0 \\
     0  &   1  &   1  &   1  &   1  &   1  &   0 \\
     0  &   0  &   1  &   1  &   1  &   1  &   0
\end{array}  \right)
\end{eqnarray*}
In~\cite{PJB} the values $\check \rho (\cM)$ and $\widehat \rho(\cM)$ were
computed  approximately to the following accuracy:
$$
4.525 \quad \le \quad   \check \rho (\cM)\quad  \le \quad 4.6105 \ ; \qquad
4.72 \quad \le \quad   \widehat \rho (\cM)\quad  \le \quad 4.8\, .
$$
Algorithms~\textbf{(P)} and~\textbf{(L)}  determine their precise
values:
\begin{eqnarray*}
\check \rho ({\cM})  & = & \bigl[\rho( A_1\, A_2 )\bigr]^{1/2} \quad = \quad 4.61047781\ldots
\\
\widehat \rho({\cM}) & = & \bigl[ \rho ( A_2\,A_3)\bigr]^{1/2} \quad = \quad 4.72204513\ldots .
\end{eqnarray*}

{\em The JSR computation.} Algorithm~\textbf{(L)} starting with
the product $\Pi\, =\, A_2\,A_3$ terminates after 4
steps producing the extremal infinite polytope~$P_3 = {\rm co}_{-}\{ v_i
\}_{i=1}^{12}$, where $v_1$ is the leading eigenvector of $\widetilde
\Pi$, $v_2\, =\, \widetilde A_3 v_1\, $ (step zero); $\, v_3\, =\, \widetilde
A_1v_1$, $\, v_4\, =\, \widetilde A_2 v_1$,  $\, v_5 \, = \, \widetilde
A_1 v_2$, and $\, v_6 \, = \, \widetilde A_3 v_2\, $ (first step);
$\, v_7\, =\, \widetilde A_3 v_4$, $\, v_8\, =\, \widetilde A_1 v_6$,
$\, v_9\, =\, \widetilde A_2 v_6$, and $\, v_{10}\, =\, \widetilde A_3 v_6\, $
(second step); $\, v_{11}\, =\, \widetilde A_1 v_8\, $ and $\, v_{12}\,
=\, \widetilde A_2 v_8\, $ (third step). See the corresponding cyclic tree in Figure~\ref{fig:plot1}.
Thus $\Pi\, =\, A_2\,A_3$ is
an s.m.p.

\begin{figure}[ht]
\centering
\global\def\path{#1}\input{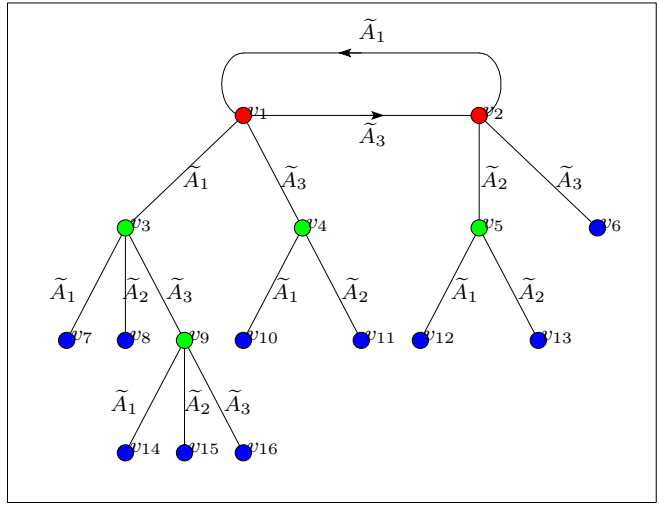}
 \caption{The vertices $\{v_i\}_{i=1}^{16}$ for LSR computation of the ternary Euler partition function
 \label{fig:plot2}}
\end{figure}

{\em The LSR computation.} Algorithm~\textbf{(P)} starting with
the product $\Pi\, =\, A_1\,A_2$ terminates after 4 steps,
the extremal  polytope~$\, P_3 = {\rm co}_{+}\{ v_i
\}_{i=1}^{16}$. The vertex $v_1$ is the leading eigenvector of $\widetilde
\Pi$, $v_2\, =\, \widetilde A_2 v_1\, $ (step zero); $\, v_3\, =\, \widetilde
A_1v_1$, $\, v_4\, =\, \widetilde A_3 v_1$,  $\, v_5 \, = \, \widetilde
A_2 v_2$, and $\, v_6 \, = \, \widetilde A_3 v_2\, $ (first step);
$\, v_7\, =\, \widetilde A_1 v_3$, $\, v_8\, =\, \widetilde A_2 v_3$,
$\, v_9\, =\, \widetilde A_3 v_3$, $\, v_{10}\, =\, \widetilde A_1 v_4$,
$\, v_{11}\, =\, \widetilde A_2 v_4$, $\, v_{12}\, =\, \widetilde A_1 v_5$,
and $\, v_{13}\, =\, \widetilde A_2 v_5\, $
(second step); $\, v_{14}\, =\, \widetilde A_1 v_9$, $\, v_{15}\, =\, \widetilde A_2 v_9$, and $\, v_{16}\,
=\, \widetilde A_3 v_9\, $ (third step). See the corresponding cyclic tree in Figure~\ref{fig:plot2}.
Therefore $\Pi\, =\, A_1\,A_2$ is an s.l.p.


\section{ Numerical results for randomly generated matrices}
\label{sec:numer}

\bigskip
\begin{table}[thb]
\begin{center}
\begin{tabular}{|l|ccl||l|ccl|}\hline
 & \multicolumn{3}{c}{JSR}
 & \multicolumn{3}{c}{JSR}
 \\ \hline
 $d$ & $\#$ its & $\#$ vertices & s.m.p. &
 $d$ & $\#$ its & $\#$ vertices & s.m.p.
 \\
 \hline
\rule{0pt}{9pt}\noindent
 $5$ &  3  & 14 & $A_1 A_2$ &
 $6$ &  4  & 26 & $A_1$
\\
\rule{0pt}{9pt}\noindent
 $5$ &  7  & 23 & $A_1 A_2^2$ &
 $6$ &  9  & 51 & $A_1 A_2$
\\
\rule{0pt}{9pt}\noindent
 $5$ &  12 & 37 & $A_1$ &
 $6$ &  5  & 38 & $A_1^2 A_2$
\\
 \hline
\rule{0pt}{9pt}\noindent
 $7$ &  17  & 100 & $A_1$ &
 $8$ &  19  & 117 & $A_1^3 A_2 A_1^4 A_2$
\\
\rule{0pt}{9pt}\noindent
 $7$ &  12  & 140 & $A_1^3 A_2 A_1 A_2$ &
 $8$ &  8   & 49  & $A_1$
\\
\rule{0pt}{9pt}\noindent
 $7$ &  24  & 223 & $A_1^3 A_2^2$ &
 $8$ &  12  &  75 & $A_1 A_2^3$
\\
 \hline
\rule{0pt}{9pt}\noindent
  $9$ &  18  & 177 & $A_1^8 A_2$ &
 $10$ &  16  & 239 & $A_1 A_2^4$
\\
\rule{0pt}{9pt}\noindent
  $9$ &  13  & 172 & $A_1^3 A_2 A_1 A_2$ &
 $10$ &  9   & 109  & $A_1$
\\
\rule{0pt}{9pt}\noindent
  $9$ &  10  & 129 & $A_2$ &
 $10$ &  24  & 408 & $(A_1^3 A_2)^2 A_2$
\\
 \hline
\rule{0pt}{9pt}\noindent
 $11$ &  20  & 707 & $A_1^3 A_2^2$ &
 $12$ &  31  & 1539 & $A_1 A_2 A_1^2 A_2^2$
\\
\rule{0pt}{9pt}\noindent
 $11$ &  14  & 340 & $A_1^2 A_2 A_1 A_2$ &
 $12$ &  9   & 211  & $A_1 A_2$
\\
\rule{0pt}{9pt}\noindent
 $11$ &  12  & 183 & $A_1^3 A_2$ &
 $12$ &  13  & 215 & $A_1 A_2^3$
\\
 \hline
\rule{0pt}{9pt}\noindent
 $15$ &  18  & 715  & $A_1^2 A_2 A_1 A_2^4$ &
 $20$ &  21  & 1539 & $A_1 A_2$
\\
\rule{0pt}{9pt}\noindent
 $15$ &  14  & 570   & $A_1^4 A_2$ &
 $20$ &  16  & 1219  & $A_1 A_2^2$
\\
\rule{0pt}{9pt}\noindent
 $15$ &  14  & 390  & $A_2$ &
 $20$ &  16  & 1247 & $A_1^2 A_2^2$
 \\
\hline
\end{tabular} \\[2mm]
\caption{Computation of the JSR for random pairs of matrices with equal norm.}\label{tab:randgen}
\end{center}
\end{table}

In this section we report some results obtained for families consisting of a
pair of random matrices of variable dimensions $d$.

The results show that the computation complexity increases significantly
as the dimension increases but also confirm the effectiveness of the method
for computing the joint and the lower spectral radius of nonnegative matrices.
We expect in general that the reachable dimension for a computation in a
reasonable time might be quite high for a set of operators sharing an invariant
cone.

First we consider the general case of two random matrices with normally
distributed entries. The generated random matrices are scaled to have
equal spectral norm. This aims to reduce the number of cases were the s.m.p is the
matrix with larger spectral radius.
The first column of Table \ref{tab:randgen} gives the dimension, the second column the number
of iterations for Algorirhm~\textbf{(R)} to converge, the third column provides the number
of vertices of the extremal polytope and the last column gives the correspondent s.m.p.; we
immediately observe that the complexity (in terms of iterations and number of vertices)
rapidly increases with the dimension. Dealing with two $20 \times 20$ matrices can be
considered a challenging computational problem.

Then we consider in Table \ref{tab:randpos} randomly generated nonnegative matrices still scaled to
have the same norm.
\begin{table}[htb]
\begin{center}
\begin{tabular}{|l|ccl|ccl|}\hline
 & \multicolumn{3}{c}{JSR}
 & \multicolumn{3}{c}{LSR}\\ \hline
 $d$ & $\#$ its & $\#$ vertices & s.m.p. & $\#$ its & $\#$ vertices & s.l.p. \\
 \hline
\rule{0pt}{9pt}\noindent
 $10$ &  3 & 6 & $A_1 A_2$ & 4 & 6 & $A_1 A_2^2$ \\
\rule{0pt}{9pt}\noindent
 $10$ &  3 & 4 & $A_1$ & 4 & 5 & $A_2$ \\
\rule{0pt}{9pt}\noindent
 $10$ &  4 & 6 & $A_1 A_2$ & 7 & 15 & $A_1^2 A_2^2$ \\
\rule{0pt}{9pt}\noindent
 $10$ &  6 & 11 & $A_1^2 A_2^2$ & 3 & 6 & $A_1 A_2$ \\
\rule{0pt}{9pt}\noindent
 $10$ &  4 & 8 & $A_1 A_2^2$ & 5 & 9 & $A_1^2 A_2$ \\
\hline
\rule{0pt}{9pt}\noindent
 $20$ &  4 & 7 & $A_2$ & 4 & 6 & $A_1$ \\
\rule{0pt}{9pt}\noindent
 $20$ &  4 & 6 & $A_1 A_2$ & 5 & 9 & $A_1 A_2^2 $ \\
\rule{0pt}{9pt}\noindent
 $20$ &  6 & 14 & $A_1^2 A_2$ & 3 & 4 & $A_1 A_2$ \\
\rule{0pt}{9pt}\noindent
 $20$ &  5 & 11 & $A_1 A_2^2$ & 6 & 14 & $A_1^2 A_2$ \\
\rule{0pt}{9pt}\noindent
 $20$ &  5 & 9 & $A_1 A_2$ & 3 & 4 & $A_1$ \\
\hline
\end{tabular} \\[2mm]
\caption{Computation of the JSR and of the LSR for random nonnegative pairs of matrices.}\label{tab:randpos}
\end{center}
\end{table}

Finally we consider binary matrices and vary the density of the number of
zero entries. We scale the pairs of matrices to have the same spectral radius;
note that in some cases either the s.m.p. or the s.l.p. are the starting matrices so that there
is no guarantee of the convergence of the algorithm we propose. Nevertheless,
when we report the number of iterations and of vertices we imply that the algorithm
has converged in a finite number of steps.

Tables \ref{tab:randbin50} and \ref{tab:randbin100} report the results obtained for pairs of matrices
respectively of dimension $d=50$ and $d=100$.
Whenever both $A_1$ and $A_2$ are either s.m.p.'s or s.l.p.'s we indicate both the numbers of iterations/vertices
taking either $A_1$ or $A_2$ as optimal product.

\begin{table}[htb]
\begin{center}
\begin{tabular}{|l|ccl|ccl|}\hline
 & \multicolumn{3}{c}{JSR}
 & \multicolumn{3}{c}{LSR}\\ \hline
 density & $\#$ its & $\#$ vertices & s.m.p. & $\#$ its & $\#$ vertices & s.l.p. \\
 \hline
\rule{0pt}{9pt}\noindent
 $0.2$ &  9 & 55 & $A_1 A_2^2$ & 4 & 8 & $A_1 \ \mbox{and} \ A_2$ \\
\rule{0pt}{9pt}\noindent
 $0.2$ &  5 & 17 & $A_1 A_2$ & 5 & 10 & $A_1^2 A_2$ \\
\rule{0pt}{9pt}\noindent
 $0.2$ &  8 & 24 & $A_1^2 A_2^2$ & 4 (4) & 6 (6)& $A_1 \ \mbox{and} \ A_2$ \\
\rule{0pt}{9pt}\noindent
 $0.2$ &  5 & 16 & $A_1^2 A_2$ & 4 (5) & 6 (8) & $A_1 \ \mbox{and} \ A_2$ \\
\rule{0pt}{9pt}\noindent
 $0.2$ & 14 & 59 & $A_1 A_2^3$ & 5 & 10 & $A_1 A_2$ \\
\hline
\rule{0pt}{9pt}\noindent
 $0.5$ &  4 & 8 & $A_1 A_2$ & 4 & 10 & $A_1 A_2$ \\
\rule{0pt}{9pt}\noindent
 $0.5$ &  5 & 14 & $A_1^2 A_2$ & 4 (3) & 5 (4) & $A_1 \ \mbox{and} \ A_2$ \\
\rule{0pt}{9pt}\noindent
 $0.5$ &  6 & 15 & $A_1 A_2^2$ & 6 & 17 & $A_1^2 A_2^2$ \\
\rule{0pt}{9pt}\noindent
 $0.5$ &  5 & 16 & $A_1 A_2$ & 4 (4) & 6 (5) & $A_1 \ \mbox{and} \ A_2$ \\
\rule{0pt}{9pt}\noindent
 $0.5$ &  6 & 20 & $A_1^3 A_2$ & 5 & 9 & $A_1 A_2$ \\
\hline
\rule{0pt}{9pt}\noindent
 $0.75$ &  5 & 16 & $A_1 A_2^2$ & 5 (7) & 12 (14) & $A_1 \ \mbox{and} \ A_2$ \\
\rule{0pt}{9pt}\noindent
 $0.75$ &  4 &  8 & $A_1 A_2$ & 6 & 16 & $A_1 A_2$ \\
\rule{0pt}{9pt}\noindent
 $0.75$ &  5 & 11 & $A_1 A_2^3$ & 6 & 19 & $A_1 A_2^2$ \\
\rule{0pt}{9pt}\noindent
 $0.75$ &  5 & 16 & $A_1^2 A_2^2$ & 11 & 170 & $A_1 A_2^5$ \\
\rule{0pt}{9pt}\noindent
 $0.75$ &  5 & 12 & $A_1^2 A_2$ & 5 (6) & 13 (12) & $A_1 \ \mbox{and} \ A_2$ \\
\hline
\rule{0pt}{9pt}\noindent
 $0.9$ &  4 (5) & 8 (9) & $A_1 \ \mbox{and} \ A_2$ & 4 & 8 & $A_1 A_2$ \\
\rule{0pt}{9pt}\noindent
 $0.9$ &  5 & 9 & $A_1^2 A_2$ & 6 & 4 & $A_1 A_2$ \\
\rule{0pt}{9pt}\noindent
 $0.9$ &  3 & 4 & $A_1 A_2$ & 7 (8) & 11 (12) & $A_1 \ \mbox{and} \ A_2$ \\
\rule{0pt}{9pt}\noindent
 $0.9$ &  4 & 11 & $A_1 A_2^2$ & 4 & 7 & $A_1 A_2$ \\
\rule{0pt}{9pt}\noindent
 $0.9$ &  7 & 14 & $A_1^3 A_2^2$ & 8 (8) & 13 (11) & $A_1 \ \mbox{and} \ A_2$ \\
\hline
\end{tabular} \\[2mm]
\caption{Computation of the JSR and of the LSR for random pairs of binary matrices of dimension $d=50$.}\label{tab:randbin50}
\end{center}
\end{table}

\begin{table}[htb]
\begin{center}
\begin{tabular}{|l|ccl|ccl|}\hline
 & \multicolumn{3}{c}{JSR}
 & \multicolumn{3}{c}{LSR}\\ \hline
 density & $\#$ its & $\#$ vertices & s.m.p. & $\#$ its & $\#$ vertices & s.l.p. \\
 \hline
\rule{0pt}{9pt}\noindent
 $0.2$ &  6 & 24 & $A_1^2 A_2$ & 6 & 31 & $A_1 A_2^2$ \\
\rule{0pt}{9pt}\noindent
 $0.2$ &  6 & 23 & $A_1 A_2$ & 6 & 28 & $A_1^2 A_2^2$ \\
\rule{0pt}{9pt}\noindent
 $0.2$ &  7 & 27 & $A_1 A_2^3$ & 6 & 20 & $A_1 A_2$ \\
\rule{0pt}{9pt}\noindent
 $0.2$ &  5 & 21 & $A_1 A_2^2$ & 7 & 24 & $A_1^2 A_2$ \\
\hline
\rule{0pt}{9pt}\noindent
 $0.5$ &  5 & 10 & $A_1 A_2$ & 5 & 15 & $A_1 A_2^2$ \\
\rule{0pt}{9pt}\noindent
 $0.5$ &  6 & 17 & $A_1^2 A_2$ & 4 & 8 & $A_1 A_2$ \\
\rule{0pt}{9pt}\noindent
 $0.5$ &  6 & 18 & $A_1^2 A_2^2$ & 5 & 16 & $A_1^2 A_2$ \\
\rule{0pt}{9pt}\noindent
 $0.5$ &  6 & 22 & $A_1 A_2^3$ & 4 (6) & 9 (14) & $A_1 \ \mbox{and} \ A_2$ \\
\hline
\rule{0pt}{9pt}\noindent
 $0.8$ &  4 & 7 & $A_1 A_2$ & 4 & 7 & $A_1 A_2$ \\
\rule{0pt}{9pt}\noindent
 $0.8$ &  7 & 18 & $A_1^2 A_2$ & 6 & 14 & $A_1^2 A_2^2$ \\
\rule{0pt}{9pt}\noindent
 $0.8$ &  5 & 14 & $A_1 A_2^2$ & 9 (7) & 14 (16) & $A_1 \ \mbox{and} \ A_2$ \\
\rule{0pt}{9pt}\noindent
 $0.8$ &  5 & 12 & $A_1^3 A_2$ & 5 & 12 & $A_1 A_2^2$ \\
\hline
\end{tabular} \\[2mm]
\caption{Computation of the JSR and of the LSR for random pairs of binary matrices of dimension $d=100$.}\label{tab:randbin100}
\end{center}
\end{table}

Some comments are necessary. The computations for the general case with $d=10$ need usually
a few minutes. The computations for the general case with $d=20$ need usually
between half an hour and one hour of computation but for some examples till $8$ hours
(in a standard laptop with i5 processor).
We hypothesize that the overall computation depend on several factors, not only the
length of the spectrum maximizing product but also the ratio between the leading
eigenvalues of the family and the closer ones that is eigenvalues of products
which have modulus close to $1$ and on the distribution of the vertices of the
extremal polytope in all the orthants.

In the nonnegative case all the vertices lie in the nonnegative orthant and this determine
a much lower complexity. The presence of quasi-optimal products that is products with
eigenvalues close to $1$ is a factor of slowdown also in this case.

For binary matrices we observe from the experiments that the behavior of the
algorithm slightly depends on the density of the zero entries and also on the
dimension. This implies that we are able to compute the joint spectral characteristics
of possibly large binary matrices with any density.

In some applications the families of the matrices have the same spectral radii and it may
happen that the optimal products are exactly the matrices themselves. In our experience there
are cases where we have been able to compute an extremal polytope invariant set by starting from
any optimal matrix of the family (see also the tables of results) but we have also encountered
cases where the algorithm has not terminated finitely. We think that an interesting open problem
is that of balancing leading eigenvectors associated to different products (which
are not powers or cyclic permutation one of the other).

Although it is true that this situation is not generic there are some applications where it naturally
occurs. We leave this topic to a future investigation.


\section{ Appendix}
\label{sec:app}

We give here the proofs of two main results of this paper
(Theorems \ref{th.cond-r} and \ref{th.cond-min})
and details about the proof of Theorem~\ref{th100}.

\smallskip

\subsection{Proof of Theorem~\ref{th.cond-r}.}

We give the proof for Algorithm~\textbf{(R)}, the proofs for
Algorithms~\textbf{(C)} and~\textbf{(P)} are analogous. We use two
auxiliary results.
\begin{lemma}\label{l.tree1r}
Let us have a cyclic tree $\cT$ with a root $\mathbf{B}$ generated
by an irreducible  word $\mathbf{b} = d_n\ldots d_1$; then for any
word $\mathbf{a}$, which is not a power of $\mathbf{b}$, we have
$\mathbf{a}^n \notin \mathbf{B}$.
\end{lemma}
{\tt Proof.} Let $l$ be the length of $\mathbf{a}$ and $p$ be and
the greatest common divisor of $l$ and $n$. If $l=kn$ for some
integer $k$, then the words $\mathbf{a}$ and $\mathbf{b}^k$ must
have different letters at some position, otherwise $\mathbf{a} =
\mathbf{b}^k$. Therefore $\mathbf{a} \notin \mathbf{B}$, and so
$\mathbf{a}^{n}\notin \mathbf{B}$.

If $l$ is not divisible by $n$, then $p < n$, and there exists an
index $j$ such that $d_{j+p}\ne d_j$, otherwise~$\mathbf{b}$ is a
power of the word $d_p\ldots d_1$, which contradicts the
irreducibility. The Diophantine equation $lx - n y = p$ has a
solution $(x, y) \in \z^2$ such that $0 \le x \le n-1$. Since the
words $\mathbf{a}^{x+1}$ and $\mathbf{b}^{y+2}$ have different
letters at the position $lx + j$, we have $\mathbf{a}^{x+1}\notin
\mathbf{B}$, and hence $\mathbf{a}^{n}\notin \mathbf{B}$, because
$n \ge x+1$.

{\hfill $\Box$}
\smallskip

\begin{lemma}\label{l.tree2r}
Let $\cT$ be the cyclic tree generated by the product $\widetilde
\Pi$. If Algorithm~\textbf{(R)} terminates within~$k$ steps, then
there is $\varepsilon > 0$ such that all nodes of $\cT$ of level
$\ge k$ are in the polytope $(1-\varepsilon)\, P_{k-1}$.
\end{lemma}
{\tt Proof.} If the algorithm terminates after $k$ steps, then
$\widetilde A_jP_{k-1} \subset P_{k-1}$ for all $j = 1, \ldots , m$.
Moreover we have $\, \cU_k = \emptyset$, which means that every node
$v\in \cT$ of the level $k$ belongs to a dead branch, and
therefore   $v \in {\rm int}\, P_{k-1}$. The total number of nodes
of level $k$ is finite, hence all of them are in $\,
(1-\varepsilon)\, P_{k-1}$ for some $\varepsilon > 0$. If $v$ is a
node of a bigger level $r > k$, then $v = R\, v_0 $, where $R \in
\widetilde \cM^{\, r-k}$ and~$v_0$ is some node of the $k$th level.
Since $v_0 \in (1-\varepsilon)\, P_{k-1}$, we have $v \in
(1-\varepsilon)\, R P_{k-1} \subset (1-\varepsilon)\, P_{k-1} $,
because $RP_{k-1} \subset P_{k-1}$.

{\hfill $\Box$}
\smallskip

{\tt Proof of Theorem~\ref{th.cond-r}}.

{\em Necessity.} Consider the cyclic tree $\cT$ generated by the product~$\widetilde \Pi$.
Assume the algorithm terminates after $k$ steps. By
Lemma~\ref{l.tree2r} all nodes of levels at least $k$ belong to
$(1-\varepsilon)\, P_{k-1}$, where $\varepsilon > 0$ is fixed. For
every product $C$, which is not a power of $\widetilde \Pi_i$,  the node
$C^{\, n}v_i$ does not belong to the root (Lemma~\ref{l.tree1r}). Hence
for each $v_i$ from the root and for every product~$C\in {\widetilde
M}^l$ that is not a power of a cyclic permutation of $\widetilde \Pi$, the
level of the node $C^{\, n+k}v_i$ is bigger than $k$. If $v$ is not
in the root, then this level is bigger than $ln+k > k$. Thus,
$C^{\, n+k}v \, \in \, (1-\varepsilon)\, P_{k-1}$ for each node $v \in
\cT$, and hence, for each vertex~$v$ of the polytope $P_{k-1}$.
This yields that~$C^{\, n+k}P_{k-1} \, \subset \, (1-\varepsilon)\,
P_{k-1}$. Therefore, $\rho(C^{\, n+k}) < 1-\varepsilon$, and so $\,
\rho(C) < \bigl(1-\varepsilon \bigr)^{1/(n+k)}$.  Consequently $\widetilde
\Pi$ is dominant.


Let us now show that $1$ is its unique and
simple leading eigenvalue. Since for $i\ne 1$ the product $\widetilde
\Pi_1$ is not a power of $\widetilde \Pi_i$, it follows that the node
$\widetilde \Pi_1^nv_i$ does not belong to the root
(Lemma~\ref{l.tree2r}). Hence, the level of the node $\widetilde
\Pi_1^{\, n+k}v_i$ is bigger than $k$. If $v$ is not in the root,
then the level of $\widetilde \Pi_1^{\, n+k}v$ is bigger than~$k$ as
well. Thus, $\widetilde \Pi_1^{\, n+k}v \, \in \, (1-\varepsilon)\,
P_{k-1}$ for all vertices~$v$ of~$P_{k-1}$, except for $v \, = \,
\pm \, v_1$. For any eigenvector $u\ne v_1$ of the operator
$\widetilde \Pi_1$ take the one-dimensional subspace $U \subset \re^d$
spanned by~$u$ (if~$u$ is complex, then $U$ is the two-dimensional
subspace spanned by $u$ and by its conjugate).  Since $v_1 \notin
U $, it follows that $\widetilde \Pi_1 (P_{k-1}\cap U) \, \subset \,
{\rm int}\, (P_{k-1}\cap U)$, where the interior is taken in $U$.
This implies that the spectral radius of $\widetilde \Pi_1$ on the
subspace~$U$ is smaller than~$1$. Thus, all eigenvalues of~$\widetilde
\Pi_1$ different from~$1$ are smaller than $1$ by modulo, and the
eigenvalue~$1$ has a unique eigenvector. Hence, the leading
eigenvalue~$1$ is unique and has only one Jordan block. The
dimension of this block cannot exceed one, otherwise
$\|\widetilde \Pi_1^k\|\to \infty$ as $k \to \infty$, which
contradicts the nondefectivity of the family~$\widetilde \cM$.
Therefore, the eigenvalue~$1$ is simple.
\smallskip

{\em Sufficiency}. The proof uses similar arguments as this given for the
{\em Small CPE Theorem} in \cite{GWZ}, to which we refer the reader.

Assume $\widetilde \Pi$ is dominant and its
leading eigenvalue is unique and simple. If the algorithm does not
terminate, then the tree $\cT$ has an infinite path of alive
leaves
 $v^{(0)} \to v^{(1)} \to \ldots $
(the node $v^{(i)}$ is on the $i$th level) starting at a node
$v^{(0)}= v_p$ from the root. For every $r$ we have $v^{(r)} \notin
{\rm int}\, P_{r-1}$. Hence $v^{(r)} \notin {\rm int}\, P_k$ for
all $k<r$. Since the family $\widetilde \cM$ is irreducible, it
follows that $0 \in {\rm int}\, P_k$ for some~$k$, and hence
the polytope $P_k$ defines a  norm $\|\cdot \|_k$ in~$\re^d$. For
this norm $\|v^{(r)}\|\ge 1$ for all $r>k$. On the other hand,
$\widetilde \cM$ is nondefective, hence the sequence $\{v^{(r)}\}$ is
bounded. Thus, there is a subsequence $\{v^{(r_i)}\}_{i\in \n}\, ,
\, r_1 \ge k,\, $
 that converges to some point $v \in \re^d$.
Clearly, $\|v\|_k\ge 1$. For every $i$ we have $v^{(r_i)} =
R_iv^{(r_1)}$ and $v^{(r_{i+1})}\, = \, C_i\, v^{(r_i)}$, where $R_i
\in \widetilde \cM^{\, r_{i}-r_1}\, , \, C_i \in \widetilde \cM^{\,
r_{i+1}-r_i}$. Denote by $\cl[\widetilde M]$ the closure of the semigroup
of all products of operators from~$\widetilde \cM$. Since this
semigroup is bounded, after possible passage to a subsequence, it
may be assumed that $R_i$ and $C_i$ converge to some $R, C \in
\cl[\widetilde M]$ respectively as $i \to \infty$. We have $Cv = v$,
hence $\rho(C) \ge 1$, which, by the domination assumption,
implies that there is $j\in \{1, \ldots , n\}$ such that $C$
belongs to $\cl[\widetilde \Pi_j]$, which is the closure of the semigroup
$\{(\Pi_j)^q\}_{q\in \n}$. Moreover, since the leading eigenvalue
of $\widetilde \Pi$ is unique and simple, we see that $v=\lambda v_j$,
where $\lambda \in \re$. We have, $\|v_j\|_k = 1$ and $\|v\|_k\ge
1$, hence $|\lambda|\ge 1$. Thus, $Rv^{(r_1)} = \lambda v_j$. The
nodes $v_j$ and $v^{(0)} = v_p$ are both from the root, hence
there is a product $S$ such that $Sv_{j} = v^{(1)}$. Taking into
account that $v^{(r_1)} = R_1v^{(1)} $, we obtain $R_1SR v^{(1)} =
\lambda v^{(1)}$. Hence $\rho(R_1SR)\ge |\lambda |$, and we
conclude that $\lambda = \pm 1$ and that $R_1SR \in  \cl[\widetilde
\Pi_i]$ for some~$i$. This yields $v^{(1)}\, = \, \mu v_i$,
where~$v_i$ is the corresponding vector from the root, $\mu \in
\re$. Since $\|v_i\| = 1$ and $\|v^{(1)}\| \ge 1$, we have $|\mu|
\ge 1$. The elements $v_i$ and $v^{(0)}=v_p$ are both from the
root, hence $\widetilde A_{p}\cdots \widetilde A_i v_i = v_p$, and
consequently $\widetilde A_s\widetilde A_{p}\cdots  \widetilde A_i v_i =
v^{(1)}$ for some~$\widetilde A_s \in \cM$. Note that $d_s \ne
d_{p+1}$, because the node $v^{(1)}$ is not in the root.
Therefore, the product $Q = \widetilde A_s \widetilde A_{p}\cdots \widetilde
A_i$ does not coincide with $\Pi_i$, and its length is at
most~$n$. Thus, $Q\, v_i = \mu v_i$, hence $\rho(Q) = 1$ and $\mu =
\pm 1$. Thus, $Q$ has spectral radius~$1$ and the leading
eigenvector~$v_i$, therefore $Q \in \cl[\widetilde \Pi_i]$. On the other
hand the length of $Q$ does not exceed $n$, hence $Q = \widetilde
\Pi_i$, which is a contradiction. Hence, the algorithm terminates
within finitely many steps.

{\hfill $\Box$}
\smallskip

\subsection{Proof of Theorem~\ref{th.cond-min}.}

We use several auxiliary results. The proof of the following lemma
is similar to the proof of Lemma~\ref{l.tree2r}.
\begin{lemma}\label{l.tree2p}
Let $\cT$ be the cyclic tree generated by the product $\widetilde
\Pi$. If the algorithm terminates within $k$ steps, then there is
$\varepsilon > 0$ such that all  vertices of the tree of level
$\ge k$ belong to the infinite polytope $(1+\varepsilon)\,
P_{k-1}$.
\end{lemma}
\begin{lemma}\label{l.tree4p}~\cite{Va}
If an operator $B$ has an invariant cone~$K$, then for every its
eigenvector from  ${\rm int}\, K$ the corresponding eigenvalue
equals to~$\rho(B)$.
\end{lemma}
\smallskip

 {\tt Proof of Theorem~\ref{th.cond-min}}.

{\em Necessity.} Consider the cyclic tree $\cT$ generated by the
product~$\widetilde \Pi$. If the algorithm terminates after $k$ steps, then
by Lemma~\ref{l.tree2p} all vertices of levels at least $k$ belong to $(1+\varepsilon)\,
P_{k-1}$. For an arbitrary  product $C$, which is not a power of $\widetilde
\Pi_i$, the node $C^{\, n}v_i$ does not belong to the root
(Lemma~\ref{l.tree1r}). Hence for every $v_i$ from the root the
level of the node $C^{\, n+k}v_i$ is bigger than $k$. If $v$ is not
in the root, then the level of $C^{\, n+k}v$ is bigger than $k$ as
well, consequently $C^{\, n+k}v \, \in \, (1+\varepsilon)\, P_{k-1}$
for each node $v \in \cT$, and hence for every vertex of  $P_k$. This
yields that~$C^{\, n+k}P_{k-1} \, \subset \, (1+\varepsilon)\, P_{k-1}$,
therefore $\rho(C^{\, n+k}) > 1+\varepsilon$, and so $\rho(C) >
\bigl(1+\varepsilon \bigr)^{1/(n+k)}$. This holds for every
product $C$ that is not a power of $\widetilde \Pi$ or of its cyclic
permutations, which completes the proof.
\smallskip

{\em Sufficiency.} Assume the converse: the product $\widetilde \Pi$
is under-dominant, but the algorithm does not produce an extremal
infinite polytope. This means that the tree $\cT$ has an infinite
path of alive leaves $v^{(0)} \to v^{(1)} \to \ldots $ starting at
a vertex $v^{(0)}= v_p$ from the root. Since the family~$\cM$ is
eventually positive, it follows from Lemma~\ref{l.pmin10} that
there exists an internal  invariant cone~$\widetilde K$, which,
moreover, contains all leading eigenvectors of products of
operators from~$\cM$. Hence, $\widetilde K$ contains the root
of~$\cT$, and therefore, it contains all the nodes~$v^{(k)}$. For
every $r$ we have $v^{(r)} \notin {\rm int}\, P_{r-1}$. Hence $v^{(r)}
\notin {\rm int}\, P_k$ for all $k < r$. Let $g(\cdot)$ be the
antinorm defined by the infinite polytope~$P_k$: $g(x)  = \sup\,
\bigl\{ \, \lambda \,  \bigl| \, \lambda^{-1} x \, \in \, P_k \,
\bigr\}$. Since a concave function is continuous at every interior
point of its domain (see, for instance,~\cite{MT}), it follows that $g$ is equivalent to
every norm and to every antinorm on the interior cone~$\widetilde K$.  In
particular, there are positive constants $c_1, c_2$ such that \\[1mm]
\centerline{$c_1f(x) \le g(x)\le c_2f(x)\, , \ x \in \widetilde K$}, \\[1mm]
where~$f$ is an invariant antinorm for $\widetilde \cM$ (Theorem~\ref{th20.5}).
For arbitrary $r$ we have $g(v^{(r)})\, \le \, 1$. On the other hand,
since~$f$ is invariant and $f(v_k)=1$ for all $k$, we have
$f(v)\ge 1$ for every node~$v$ of the tree. In particular, $\,
f(v^{(r)})\, \ge \, 1$. Thus, $c_1 \le g(v^{(r)}) \le 1$ for all
$r$. Since $g$ is equivalent to each norm on~$\widetilde K$, we see
that  the sequence $\{v^{(r)}\}$ is bounded, and hence there is a
subsequence $\{v^{(r_i)}\}_{i\in \n}\, , \, r_1 \ge 1\, , $
 that converges to some point $v \in \re^d$.
Clearly, $v \in \widetilde K$ and $c_1 \le g(v)\le 1$. For every $i$ we
have $v^{(r_i)} = R_iv^{(r_1)}$ and $v^{(r_{i+1})}\, = \,
C_iv^{(r_i)}$, where $R_i \in \widetilde \cM^{\, r_{i}-r_1}\, , \, C_i
\in \widetilde \cM^{\, r_{i+1}-r_i}$. The sequence $\{v^{(r_i)}\}$ is
contained in~$\widetilde K$, bounded, and separated from zero,
hence by Lemma~\ref{l.tree3p} the sequences of operators $\{R_i\}$
and $\{C_i\}$ are both bounded. Therefore, after a passage to
subsequences it may be assumed that these two sequences converge to
some $R, C \in \cl[\widetilde M]$ respectively as $i \to \infty$ (see the
proof of Theorem~\ref{th.cond-r} for the definition of $\cl[\widetilde
{\cM}]$ and $\cl[\widetilde \Pi]$). We have $Cv = v$. Since $v \in {\rm
int}\, \re^d_+$, if follows from Lemma~\ref{l.tree4p} that $v$ is
the leading eigenvector of~$C$. Consequently, $\, \rho(C) = 1$,
which, by the domination assumption, implies $C \in \cl[\widetilde
\Pi_j]$ and $v=\lambda v_j$ for some $j=1, \ldots , n$, and
$\lambda > 0$. Since $g(v_j) = 1$ and $g(v) \le 1$, it follows
that $\lambda \le 1$. Thus, $Rv^{(r_1)} = \lambda v_j$. The
elements $v_j$ and $v^{(0)} = v_p$ are both from the root, hence
there is a product $S$ such that $Sv_{j} = v^{(1)}$. Taking into
account that $v^{(r_1)} = R_1v^{(1)} $, we obtain $R_1SR v^{(1)} =
\lambda v^{(1)}$. Again invoking Lemma~\ref{l.tree4p}, we conclude
that $v^{(1)}$ is the leading eigenvector of $R_1SR$. Hence
$\rho(R_1SR) =  \lambda $, and we see that $\lambda =  1$,  hence
$R_1SR \in  \cl[\widetilde \Pi_i]$ for some~$i$. This yields $v^{(1)}\, =
\, \mu v_i$, where~$v_i$ is the corresponding vector from the
root, $\mu > 0$. Since $g(v_i) = 1$ and $g(v^{(1)}) \, \le \, 1$,
we have $\mu \le 1$. Elements $v_i$ and $v^{(0)}=v_p$ are both
from the root, hence $\widetilde A_{p}\cdots \widetilde A_i v_i = v_p$,
and consequently $\widetilde A_s\widetilde A_{p}\cdots  \widetilde A_i v_i =
v^{(1)}$ for some~$A_s \in \cM$. Note that $d_s \ne d_{p+1}$,
because the vertex $v^{(1)}$ is not in the root. Therefore, the
product $Q = \widetilde A_s \widetilde A_{p}\cdots \widetilde A_i$ does not
coincide with $\Pi_i$, and its length is at most~$n$. Thus, $Qv_i
= \mu v_i$, and by Lemma~\ref{l.tree4p} $\rho(Q) = \mu$.
Consequently, $\mu = 1$  and $v_i$ is the leading eigenvector of
the operator $Q \in \cl[\widetilde \Pi_i]$. On the other hand the length
of the product $Q$ does not exceed $n$, therefore $Q = \widetilde
\Pi_i$, which is a contradiction.

{\hfill $\Box$}
\smallskip

\subsection{The $20\times 20$-matrices $A_1, A_2$ for the problem of overlap-free words of \S \ref{subsec:OFW}
and the proof of Theorem~\ref{th100}.}
We write the two matrices $A_1$, $A_2$, associated to the problem discussed in
\S \ref{subsec:OFW}, \small
\begin{eqnarray*}
\left(\begin{array}{rrrrrrrrrrrrrrrrrrrr}
     0  &  0  &  0  &  0  &  0  &  0  &  0  &  2  &  4  &  2  &  0  &  0  &  0  &  0  &  0  &  0  &  0  &  0  &  0  &  0 \\
     0  &  0  &  1  &  1  &  0  &  1  &  1  &  0  &  0  &  0  &  0  &  0  &  0  &  0  &  0  &  0  &  0  &  0  &  0  &  0 \\
     0  &  0  &  0  &  0  &  0  &  1  &  1  &  1  &  1  &  0  &  0  &  0  &  0  &  0  &  0  &  0  &  0  &  0  &  0  &  0 \\
     0  &  0  &  1  &  1  &  0  &  0  &  0  &  0  &  0  &  0  &  0  &  0  &  0  &  0  &  0  &  0  &  0  &  0  &  0  &  0 \\
     0  &  0  &  0  &  0  &  0  &  0  &  0  &  0  &  0  &  0  &  0  &  0  &  0  &  0  &  0  &  0  &  0  &  0  &  0  &  0 \\
     0  &  1  &  0  &  0  &  1  &  0  &  0  &  0  &  0  &  0  &  0  &  0  &  0  &  0  &  0  &  0  &  0  &  0  &  0  &  0 \\
     1  &  1  &  0  &  0  &  0  &  0  &  0  &  0  &  0  &  0  &  0  &  0  &  0  &  0  &  0  &  0  &  0  &  0  &  0  &  0 \\
     0  &  0  &  0  &  0  &  0  &  2  &  0  &  0  &  0  &  0  &  0  &  0  &  0  &  0  &  0  &  0  &  0  &  0  &  0  &  0 \\
     0  &  0  &  1  &  0  &  0  &  0  &  0  &  0  &  0  &  0  &  0  &  0  &  0  &  0  &  0  &  0  &  0  &  0  &  0  &  0 \\
     0  &  0  &  0  &  0  &  0  &  0  &  0  &  0  &  0  &  0  &  0  &  0  &  0  &  0  &  0  &  0  &  0  &  0  &  0  &  0 \\
     0  &  0  &  0  &  0  &  0  &  0  &  0  &  1  &  2  &  1  &  0  &  0  &  0  &  0  &  0  &  0  &  0  &  1  &  2  &  1 \\
     0  &  0  &  1  &  1  &  0  &  1  &  1  &  0  &  0  &  0  &  0  &  0  &  0  &  0  &  0  &  0  &  0  &  0  &  0  &  0 \\
     0  &  0  &  0  &  0  &  0  &  0  &  0  &  1  &  1  &  0  &  0  &  0  &  0  &  0  &  0  &  1  &  1  &  0  &  0  &  0 \\
     0  &  0  &  0  &  0  &  0  &  0  &  0  &  0  &  0  &  0  &  0  &  0  &  1  &  1  &  0  &  0  &  0  &  0  &  0  &  0 \\
     1  &  2  &  0  &  0  &  1  &  0  &  0  &  0  &  0  &  0  &  0  &  0  &  0  &  0  &  0  &  0  &  0  &  0  &  0  &  0 \\
     0  &  0  &  1  &  0  &  0  &  1  &  0  &  0  &  0  &  0  &  0  &  0  &  0  &  0  &  0  &  0  &  0  &  0  &  0  &  0 \\
     0  &  0  &  0  &  0  &  0  &  0  &  0  &  0  &  0  &  0  &  0  &  0  &  0  &  0  &  0  &  0  &  0  &  0  &  0  &  0 \\
     0  &  0  &  0  &  0  &  0  &  0  &  0  &  1  &  0  &  0  &  0  &  0  &  0  &  0  &  1  &  0  &  0  &  0  &  0  &  0 \\
     0  &  0  &  0  &  0  &  0  &  0  &  0  &  0  &  0  &  0  &  0  &  1  &  0  &  0  &  0  &  0  &  0  &  0  &  0  &  0 \\
     0  &  0  &  0  &  0  &  0  &  0  &  0  &  0  &  0  &  0  &  1  &  0  &  0  &  0  &  0  &  0  &  0  &  0  &  0  &  0
\end{array} \right), \quad 
\left(\begin{array}{rrrrrrrrrrrrrrrrrrrr}
     0  &  0  &  0  &  0  &  0  &  0  &  0  &  1  &  2  &  1  &  0  &  0  &  0  &  0  &  0  &  0  &  0  &  1  &  2  &  1 \\
     0  &  0  &  1  &  1  &  0  &  1  &  1  &  0  &  0  &  0  &  0  &  0  &  0  &  0  &  0  &  0  &  0  &  0  &  0  &  0 \\
     0  &  0  &  0  &  0  &  0  &  0  &  0  &  1  &  1  &  0  &  0  &  0  &  0  &  0  &  0  &  1  &  1  &  0  &  0  &  0 \\
     0  &  0  &  0  &  0  &  0  &  0  &  0  &  0  &  0  &  0  &  0  &  0  &  1  &  1  &  0  &  0  &  0  &  0  &  0  &  0 \\
     1  &  2  &  0  &  0  &  1  &  0  &  0  &  0  &  0  &  0  &  0  &  0  &  0  &  0  &  0  &  0  &  0  &  0  &  0  &  0 \\
     0  &  0  &  1  &  0  &  0  &  1  &  0  &  0  &  0  &  0  &  0  &  0  &  0  &  0  &  0  &  0  &  0  &  0  &  0  &  0 \\
     0  &  0  &  0  &  0  &  0  &  0  &  0  &  0  &  0  &  0  &  0  &  0  &  0  &  0  &  0  &  0  &  0  &  0  &  0  &  0 \\
     0  &  0  &  0  &  0  &  0  &  0  &  0  &  1  &  0  &  0  &  0  &  0  &  0  &  0  &  1  &  0  &  0  &  0  &  0  &  0 \\
     0  &  0  &  0  &  0  &  0  &  0  &  0  &  0  &  0  &  0  &  0  &  1  &  0  &  0  &  0  &  0  &  0  &  0  &  0  &  0 \\
     0  &  0  &  0  &  0  &  0  &  0  &  0  &  0  &  0  &  0  &  1  &  0  &  0  &  0  &  0  &  0  &  0  &  0  &  0  &  0 \\
     0  &  0  &  0  &  0  &  0  &  0  &  0  &  0  &  0  &  0  &  0  &  0  &  0  &  0  &  0  &  0  &  0  &  2  &  4  &  2 \\
     0  &  0  &  0  &  0  &  0  &  0  &  0  &  0  &  0  &  0  &  0  &  0  &  1  &  1  &  0  &  1  &  1  &  0  &  0  &  0 \\
     0  &  0  &  0  &  0  &  0  &  0  &  0  &  0  &  0  &  0  &  0  &  0  &  0  &  0  &  0  &  1  &  1  &  1  &  1  &  0 \\
     0  &  0  &  0  &  0  &  0  &  0  &  0  &  0  &  0  &  0  &  0  &  0  &  1  &  1  &  0  &  0  &  0  &  0  &  0  &  0 \\
     0  &  0  &  0  &  0  &  0  &  0  &  0  &  0  &  0  &  0  &  0  &  0  &  0  &  0  &  0  &  0  &  0  &  0  &  0  &  0 \\
     0  &  0  &  0  &  0  &  0  &  0  &  0  &  0  &  0  &  0  &  0  &  1  &  0  &  0  &  1  &  0  &  0  &  0  &  0  &  0 \\
     0  &  0  &  0  &  0  &  0  &  0  &  0  &  0  &  0  &  0  &  1  &  1  &  0  &  0  &  0  &  0  &  0  &  0  &  0  &  0 \\
     0  &  0  &  0  &  0  &  0  &  0  &  0  &  0  &  0  &  0  &  0  &  0  &  0  &  0  &  0  &  2  &  0  &  0  &  0  &  0 \\
     0  &  0  &  0  &  0  &  0  &  0  &  0  &  0  &  0  &  0  &  0  &  0  &  1  &  0  &  0  &  0  &  0  &  0  &  0  &  0 \\
     0  &  0  &  0  &  0  &  0  &  0  &  0  &  0  &  0  &  0  &  0  &  0  &  0  &  0  &  0  &  0  &  0  &  0  &  0  &  0
\end{array}
\right)
\end{eqnarray*}
\normalsize respectively.



To give the rigorous proof of the theorem it now suffices to list all the vertices of the extremal polytopes obtained by applying
Algorithms~\textbf{(P)} and~\textbf{(L)}.


\smallskip
{\tt Proof of Theorem~\ref{th100}}. Denote $\Pi = A_1A_2$. To show that
$$\, \widehat \rho (A_1, A_2)\, = \,\bigl[ \rho(\Pi)\bigr]^{1/2}$$
it suffices to present an extremal polytope~$P$
for the operators $\widetilde A_1 = \bigl[ \rho(\Pi)\bigr]^{-1/2}A_1\, $ and $\, \widetilde A_2 = \bigl[ \rho(\Pi)\bigr]^{-1/2}A_2$.
This polytope is $\, P = {\rm co}_{-}\, \bigl( \{v_i\}_{i=1}^{54}\bigr)$, where the first vertex
$v_1$ is the leading eigenvector of $\Pi$, and the other vertices are
\small
$$
\begin{array}{llllll}
\, v_2 = Ш\widetilde A_2Ш v_1,  &
\, v_3 = Ш\widetilde A_1Ш v_1, &
\, v_4 = Ш\widetilde A_1Ш v_2,  &
\, v_5 = Ш\widetilde A_2Ш v_2, &
\, v_6 = Ш\widetilde A_2Ш v_3, &
\, v_7 = Ш\widetilde A_1Ш v_4, \\
\, v_8 = Ш\widetilde A_1Ш v_5, &
\, v_9 = Ш\widetilde A_1Ш v_6,  &
\, v_{10} = Ш\widetilde A_2ШШv_4, &
\, v_{11} = Ш\widetilde A_2Ш v_5, &
\, v_{12} = Ш\widetilde A_2Ш v_6, &
\, v_{13} = Ш\widetilde A_1Ш v_8, \\
\, v_{14} = Ш\widetilde A_1Ш v_9, &
\, v_{15} = Ш\widetilde A_1Ш v_{10}, &
\, v_{16} = Ш\widetilde A_1Ш v_{11}, &
\, v_{17} = Ш\widetilde A_1Ш v_{12}, &
\, v_{18} = Ш\widetilde A_2Ш v_7,    &
\, v_{19} = Ш\widetilde A_2Ш v_8,    \\
\, v_{20} = Ш\widetilde A_2Ш v_9,    &
\, v_{21} = Ш\widetilde A_2Ш v_{10}, &
\, v_{22} = Ш\widetilde A_2Ш v_{11}, &
\, v_{23} = Ш\widetilde A_2Ш v_{12}, &
\, v_{24} = Ш\widetilde A_1Ш v_{13}, &
\, v_{25} = Ш\widetilde A_1Ш v_{14}, \\
\, v_{26} = Ш\widetilde A_1Ш v_{16}, &
\, v_{27} = Ш\widetilde A_1Ш v_{17}, &
\, v_{28} = Ш\widetilde A_1Ш v_{19}, &
\, v_{29} = Ш\widetilde A_1Ш v_{21}, &
\, v_{30} = Ш\widetilde A_1Ш v_{22}, &
\, v_{31} = Ш\widetilde A_2Ш v_{13}, \\
\, v_{32} = Ш\widetilde A_2 Шv_{14}, &
\, v_{33} = Ш\widetilde A_2 Шv_{16}, &
\, v_{34} = Ш\widetilde A_2 Шv_{17}, &
\, v_{35} = Ш\widetilde A_2 Шv_{19}, &
\, v_{36} = Ш\widetilde A_2 Шv_{20}, &
\, v_{37} = Ш\widetilde A_2 Шv_{21}, \\
\, v_{38} = Ш\widetilde A_2 Шv_{22}, &
\, v_{39} = Ш\widetilde A_1 Шv_{25}, &
\, v_{40} = Ш\widetilde A_1 Шv_{26}, &
\, v_{41} = Ш\widetilde A_1 Шv_{27}, &
\, v_{42} = Ш\widetilde A_1 Шv_{30}, &
\, v_{43} = Ш\widetilde A_1 Шv_{32}, \\
\, v_{44} = Ш\widetilde A_1 Шv_{33}, &
\, v_{45} = Ш\widetilde A_2 Шv_{25}, &
\, v_{46} = Ш\widetilde A_2 Шv_{26}, &
\, v_{47} = Ш\widetilde A_2 Шv_{30}, &
\, v_{48} = Ш\widetilde A_2 Шv_{33}, &
\, v_{49} = Ш\widetilde A_1 Шv_{39}, \\
\, v_{50} = Ш\widetilde A_1 Шv_{43}, &
\, v_{51} = Ш\widetilde A_1 Шv_{45}, &
\, v_{52} = Ш\widetilde A_1 Шv_{46}, &
\, v_{53} = Ш\widetilde A_2 Шv_{39}, &
\, v_{54} = Ш\widetilde A_2 Шv_{45}.
\end{array}
$$
\normalsize

\smallskip

Now let $\Pi = A_1A_2^{10}$. To prove that
$$\, \check \rho (A_1, A_2)\, = \,\bigl[ \rho(\Pi)\bigr]^{1/11}$$
it suffices to present an extremal infinite polytope~$P$
for the operators $\widetilde A_1 = \bigl[ \rho(\Pi)\bigr]^{-1/11}A_1\, $ and $\, \widetilde A_2 = \bigl[ \rho(\Pi)\bigr]^{-1/11}A_2$.
This polytope is $\, P = {\rm co}\, \bigl( \{v_i\}_{i=1}^{104}\bigr)\, + \, K_{\cH}$, where the first vertex
$v_1$ is the leading eigenvector of $\Pi$, and the other vertices are
\small
$$
\begin{array}{llllll}
\,  v_{2} = Ш\widetilde A_2 v_1,&\,  v_3 = Ш\widetilde A_2v_2,&\,  v_4 =
Ш\widetilde A_2v_3,&\,  v_5 = Ш\widetilde A_2 v_4, &\,  v_6 = Ш\widetilde
A_2v_5,& \, v_7 = Ш\widetilde A_2v_6,\\
\,  v_8 = Ш\widetilde A_2v_7,& \, v_9 = Ш\widetilde A_2 v_8,& \,  v_{10} =
Ш\widetilde A_2v_9,& \, v_{11} = Ш\widetilde A_2 v_{10},& \,  v_{12} =
\widetilde A_1Шv_1,& \, v_{13} = \widetilde A_1Шv_2,\\
  \,  v_{14} =
\widetilde A_1Шv_3,& \,  v_{15} = \widetilde A_1Шv_4,&
 \, v_{16} = \widetilde A_1Шv_6,& \,  v_{17} = \widetilde
A_1Шv_7,& \, v_{18} = \widetilde A_1Шv_8,& \, v_{19} = Ш \widetilde
A_1v_9,\\
 \, v_{20} = Ш \widetilde A_1v_{10},& \,  v_{21} = Ш\widetilde
A_2v_{11},& \, v_{22} = \widetilde A_1Шv_{12}, & \,  v_{23} = \widetilde
A_1Шv_{13},& \, v_{24} = \widetilde A_1Шv_{18},& \, v_{25} = Ш \widetilde
A_1v_{19}, \\
\, v_{26} = Ш \widetilde A_1v_{20}, &
 \, v_{27} = Ш \widetilde
A_1v_{21},& \, v_{28} = Ш\widetilde A_2 v_{12},& v_{29} = Ш\widetilde A_2
v_{18}, & \, v_{30} = \widetilde A_2v_{19},&\, v_{31} = \widetilde
A_2Шv_{20}, \\
\, v_{32} = \widetilde A_2Шv_{21},&\, v_{33} = \widetilde
A_1Шv_{22}, &\,
 v_{34} = \widetilde A_1Шv_{27},&\,  v_{35} = \widetilde A_1Шv_{30},&\,  v_{36}    = \widetilde
 A_1Шv_{31},&
 \, v_{37} = Ш \widetilde A_1 v_{32},\\
  \,  v_{38} = \widetilde A_2Ш v_{27}, & \,
   v_{39} = \widetilde A_2Шv_{30}, &\,  v_{40} = \widetilde A_2Шv_{31}, & \,
 v_{41} = \widetilde A_2 v_{32}, & \,  v_{42} = \widetilde A_1Шv_{34}, &
 \, v_{43} = Ш \widetilde A_1 v_{37}, \\
 \,  v_{44} = Ш\widetilde  A_1 v_{38}, & \,
  v_{45} = Ш\widetilde  A_1v_{40}, & \,   v_{46} = Ш \widetilde A_1 v_{41}, & \,  v_{47} = \widetilde A_2
v_{34}, & \,  v_{48} = \widetilde A_2 v_{38}, & \, v_{49} = \widetilde
A_2Шv_{40}, \\
\, v_{50} = \widetilde A_1 v_{42}, &
 \,  v_{51} = \widetilde A_1 v_{43}, & \,  v_{52} = Ш\widetilde  A_1 v_{44}, & \, v_{53}
= Ш\widetilde A_1 v_{48}, & \,
 v_{54} = \widetilde A_2 v_{42}, & \,
v_{55} = \widetilde A_2Шv_{43}, \\
 \,  v_{56} = \widetilde A_2Шv_{44}, &
\, v_{57} = \widetilde A_2Ш v_{47},&
 \,  v_{58} = \widetilde A_2 v_{48}, & \,
v_{59} = \widetilde A_1Шv_{50}, & \,
 v_{60} = \widetilde A_1Шv_{52}, &\,
v_{61} = Ш\widetilde  A_1 v_{54}, \\
\,  v_{62} = \widetilde A_2Ш v_{50}, &\, v_{63} = \widetilde A_2Шv_{52}, &
\,  v_{64} = \widetilde A_2Шv_{57}, & \, v_{65} = \widetilde A_1Шv_{59}, &
\,  v_{66} = \widetilde A_1Шv_{60}, & \, v_{67} = Ш\widetilde  A_1 v_{62},
\\
 \,  v_{68} = Ш\widetilde A_2 v_{59}, & \, v_{69} = \widetilde A_2
v_{60}, & \, v_{70} = Ш\widetilde A_2 v_{62}, & \, v_{71} = \widetilde
A_1Шv_{65}, & \,  v_{72} = \widetilde A_1Шv_{67}, & \, v_{73} =
Ш\widetilde  A_1 v_{68}, \\
 \, v_{74} = Ш\widetilde A_1 v_{69}, & \,
v_{75} = \widetilde A_1 v_{70}, & \, v_{76} = \widetilde A_2 v_{65}, & \,
v_{77} = \widetilde A_2 v_{67}, & \, v_{78} = \widetilde A_2Шv_{68}, & \,
v_{79} = \widetilde A_1Шv_{71}, \\
 \, v_{80} = \widetilde A_1Шv_{72}, &\,
v_{81} = \widetilde A_1Шv_{73}, & \, v_{82} = \widetilde A_1Шv_{76}, & \,
v_{83} = \widetilde A_1Шv_{77}, & \,
 v_{84} = \widetilde A_1Шv_{78}, & \, v_{85} = Ш\widetilde A_2 v_{71}, \\
 \,  v_{86} = \widetilde A_2Шv_{72}, & \,  v_{87} = \widetilde A_2Ш v_{73}, & \,  v_{88} = Ш\widetilde A_2 v_{76}, &
 \,  v_{89} = \widetilde A_2 v_{78}, & \,  v_{90} = \widetilde A_1Шv_{79}, & \, v_{91} = \widetilde A_1Шv_{81}, \\
  \,  v_{92} = \widetilde A_1Шv_{82},&
  \,  v_{93} = \widetilde A_1Шv_{84}, & \,  v_{94} = \widetilde A_1Шv_{85}, & \,  v_{95} = \widetilde A_1Шv_{87}, &
  \,  v_{96} = \widetilde A_1 Шv_{88}, & \, v_{97} = \widetilde A_2 v_{79}, \\
   \,  v_{98} = \widetilde A_2 v_{81}, & \, v_{99} = \widetilde A_2 v_{82}, & \,  v_{100} = \widetilde A_2  v_{85}, &
   \, v_{101} = \widetilde A_2 v_{87},
& \,  v_{102} = \widetilde A_2  v_{88}, & \, v_{103} = Ш\widetilde  A_1
v_{92}, \\
 \,  v_{104} = \widetilde A_2 v_{100}\, .
\end{array}
$$
\normalsize

The proof is completed by routine computations.

{\hfill $\Box$}

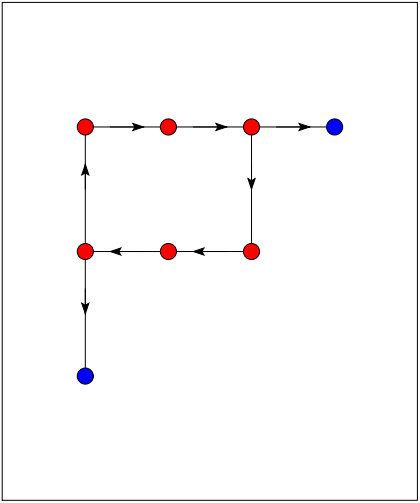
\end{document}